% TeX-Command: tex
\magnification 1200
\parskip 4pt plus 2pt

%%USE THE FOLLOWING FOUR LINES IF USING THE CURRENT ('97) VERSION
%%OF AMSPPT.
\input amstex
\documentstyle{amsppt}
\pagewidth{391.46176pt}
\pageheight{536.00175pt}

%% autonum.tex
%% Version 1.2, 11/1/2003
%%
%% Copyright (C) 1996, 2003 Michael A. Mandell
%% All rights reserved.
%% Copying of this file is authorized only if either:
%%    (1) No changes are made and the file name ``autonum.tex'' is
%%        retained.
%%    (2) The file name is changed.
%%
%% (The copyright is to ensure compatibility and control versions.)
%%
%% Please send comments or errors to mandell@math.uchicago.edu
%%
%% This file provides the macros for doing automatic number of
%% sections, theorems, and labels in AMS-TeX.  It should be loaded
%% at most once.

\ifx\AutoNum\undefined \def\AutoNum{Version 1.2, 8/8/2003}
\else \message{[already loaded]} \fi
\message{^^JAutoNum, \AutoNum.^^J}

%% The file should be \input after the style file

\ifx\styname\undefined
   \def\next{%
   \let\AutoNum\undefined
   \errhelp={It will probably fix it if you type:^^J
   I \documentstyle{amsppt}  }
   \errmessage{AutoNum needs to be loaded after the style file.^^J}
    }
\else
   \let\next\relax
\fi
\next
\def\styletest{AMSPPT}
\ifx\styname\styletest \relax
\else \message{Not expecting style \styname.  Hope this works...^^J^^J}
\fi

%% Numbered sections and subsections use the \section and \subsection
%% macros as in the following examples:
%% \section{Outline of the paper}
%% \subsection{Reduction to basic example}
%%
%% For unnumbered sections and subsections, use the \section* and
%% \subsection* macros as in the following examples:
%% \section*{Introduction}
%% \subsection*{Acknowledgments}
%%
%% To change the formatting styles of sections and subsections, \def
%% \thesection and \thesubsection macros in your TeX file.  The variants
%% \thesectiondot and \thesubsectiondot include a trailing separator.

\newcount\sectionnum
\newcount\subsectionnum

\def\thesection{%
\ifnum\sectionnum>0
\number\sectionnum%
\fi}

\def\thesectiondot{%
\ifnum\sectionnum>0
\number\sectionnum.%
\fi}

\def\thesubsection{%
\thesectiondot%
\ifnum\subsectionnum>0
\number\subsectionnum%
\fi}

\def\thesubsectiondot{%
\thesectiondot%
\ifnum\subsectionnum>0
\number\subsectionnum.%
\fi}

%% \proclaim, \definition, \example, and \remark now automatically
%% number.  For unnumbered theorems, etc., use \oldproclaim,
%% \oldexample, and \oldremark.
%%
%% By default these all use the same number scheme.

\newcount\proclaimnum
\let\definitionnum=\proclaimnum
\let\examplenum=\proclaimnum
\let\remarknum=\proclaimnum

%% To use different numbering for definitions, examples, or remarks,
%% comment out the appropriate ``let'' above put the appropriate
%% ``newcounter'' macro in your TeX file:
%\newcount\definitionnum
%\newcount\examplenum
%\newcount\remarknum

%% To change the formatting of the numbers, \def the macros below
%% in your TeX file.

\def\theproclaim{\thesubsectiondot\number\proclaimnum}
\def\thedefinition{\thesubsectiondot\number\definitionnum}
\def\theexample{\thesubsectiondot\number\examplenum}
\def\theremark{\thesubsectiondot\number\remarknum}

%% The macros \resetcounters and \resetsubcounters
%% reset the counters that should be reset at the start of a numbered
%% section and at the start of a numbered subsection (resp.).
%% (Unnumbered sections and subsections do not reset counters.)
%% Change these as appropriate in your TeX file.
%% Don't forget the ``\global'' !

\def\resetcounters{%
\global\subsectionnum=0
\global\proclaimnum=0
\global\definitionnum=0
\global\examplenum=0
\global\remarknum=0
\relax}

\def\resetsubcounters{%
\global\proclaimnum=0
\global\definitionnum=0
\global\examplenum=0
\global\remarknum=0
\relax}

\def\pagestyling#1{\ifnum #1<0 \romannumeral-#1\else\number #1\fi}

%% The following commands are also provided:
%%    \refonsamepage{lbl}{text1}{text2}
%%    \refonearlierpage{lbl}{text1}{text2}
%%    \refonlaterpage{lbl}{text1}{text2}
%% The first, refonsamepage, typesets text1 if lbl is defined and on
%% same page and otherwise typesets text2 (and gives a warning if lbl
%% is undefined).  The other two are similar but check for earlier or
%% later page, resp.  This check knows nothing about how your document
%% is formatted: it assumes pages are numbers occur linearly and that
%% page -2 comes before page -1 comes before page 0, etc.  (Typically
%% an actual document with a page -2 has page -1 BEFORE page -2 and so
%% \refonearlierpage will do the wrong thing (i.e., text1) if called
%% on page -1 about a reference on page -2.)

%% The remaining data for the last section, subsection, proclaim,
%% definition, example, or remark are available in the following
%% macros.  (Which start as empty.)

\def\therefnum{}
\def\therefname{}
\def\thereftype{}
\def\thereftext{}

%% In the current version \therefname and \thereftype are always text
%% versions of the \refname and \reftype -- this may not be what you
%% want if the proclaim, definition, example, or remark name contains
%% a control sequence.

%%
%% This is all the information needed for standard use.
%%

%% The following ``hooks'' allow for extensions by other packages.
%% They can be used, for example, to write a compatible package that
%% automatically generates a table of contents or a package that
%% prints out the labels and reference calls for debugging.
%% The hooks are called after \therefnum, \therefname, \thereftype, and
%% \thereftext are set but before any typesetting occurs.

%% These macros are called when a section and subsection macro is
%% hit, the first set in the regular numbered form and the second set
%% in the * form.

\def\xsectionhook{}
\def\xsubsectionhook{}
\def\xsectionstarhook{}

%% These are called when proclaim, definition, example, and remark
%% are hit

\def\xproclaimhook{}
\def\xdefinitionhook{}
\def\xexamplehook{}
\def\xremarkhook{}

%% Hooks are not provided for \oldproclaim, \olddefinition,
%% \oldexample, and \oldremark because no data is analyzed.

%% These hooks are called when a label is set and called.  Unlike the
%% previous hooks these have arguments.  For \xlabelhook, the argument
%% is the name of the label, i.e., in the call \label{lbl}, the
%% argument is passed as lbl. For \xrefcallhook #1 is the csname of
%% the reference call and #2 is the label name.  For example,
%% \refnum{lbl} calls \xrefcallhook{\refnum}{lbl}.

\def\xlabelhook#1{}
\def\xrefcallhook#1#2{}

%% The remaining comments are technical documentation.

%% LABELS AND REFERENCES

%% We must keep track of the page number and also make sure that the
%% reference is added to the file when references are saved.  In tex
%% the only way to keep track of the page numbers is in macros (or
%% token #'s but we would quickly run out of those).  Obviously it
%% would be nice for these macros not to conflict with the user's
%% macros.  We make conflict extremely unlikely by having our macros
%% include the character ``%''.  On the other hand to keep the
%% documentation going, we will need a new comment character; we
%% choose ``;''. (Note that there should be no problem for the user to
%% include ``;'' as a charcter in a label name.)

\catcode`\%=11
\catcode`\;=14

;; We will keep track of data in the macros ??ref%lbl
;; where lbl is the label name as in \label{lbl} and ?? correspoonds
;; to the type of data.  Of course, before trying to expand ??ref%lbl,
;; need to make sure that it is defined

\def\refnum#1{\xrefcallhook{\refnum}{#1}\relax\expandafter\ifx\csname nref%#1\endcsname\relax\ref%%n%f{#1}{num}{?}\else\csname nref%#1\endcsname\fi}

\def\refname#1{\xrefcallhook{\refname}{#1}\relax\expandafter\ifx\csname naref%#1\endcsname\relax\ref%%n%f{#1}{name}{??????}\else\csname naref%#1\endcsname\fi}

\def\reftype#1{\xrefcallhook{\reftype}{#1}\relax\expandafter\ifx\csname tyref%#1\endcsname\relax\ref%%n%f{#1}{type}{??????}\else\csname tyref%#1\endcsname\fi}

\def\reftext#1{\xrefcallhook{\reftext}{#1}\relax\expandafter\ifx\csname tref%#1\endcsname\relax\ref%%n%f{#1}{text}{??????}\else\csname tref%#1\endcsname\fi}

\def\refpage#1{\xrefcallhook{\refpage}{#1}\relax\expandafter\ifx\csname pref%#1\endcsname\relax\ref%%n%d{#1}?\else\pagestyling{\csname pref%#1\endcsname}\fi}

\def\refonsamepage#1#2#3{\expandafter\ifx\csname pref%#1\endcsname\relax\ref%%n%d{#1}{#3}\else\expandafter\ifnum\csname pref%#1\endcsname=\pageno #2\else #3\fi\fi}

\def\refonearlierpage#1#2#3{\expandafter\ifx\csname pref%#1\endcsname\relax\ref%%n%d{#1}{#3}\else\expandafter\ifnum\csname pref%#1\endcsname<\pageno #2\else #3\fi\fi}

\def\refonlaterpage#1#2#3{\expandafter\ifx\csname pref%#1\endcsname\relax\ref%%n%d{#1}{#3}\else\expandafter\ifnum\csname pref%#1\endcsname>\pageno #2\else #3\fi\fi}

\def\ref%%n%f#1#2{\expandafter\ifx\csname pref%#1\endcsname\relax\ref%%n%d{#1}\else\immediate\write16{Reference #1 has no #2 data (while proc. p. \pagenumber).}\fi}
\def\ref%%n%d#1{\immediate\write16{Reference #1 undefined (while proc. p. \pagenumber).}}

\def\setrefnum#1#2{\expandafter\gdef\csname nref%#1\endcsname{#2}}
\def\setrefname#1#2{\expandafter\gdef\csname naref%#1\endcsname{#2}}
\def\setreftype#1#2{\expandafter\gdef\csname tyref%#1\endcsname{#2}}
\def\setreftext#1#2{\expandafter\gdef\csname tref%#1\endcsname{#2}}
\def\setrefpage#1#2{\expandafter\gdef\csname pref%#1\endcsname{#2}}

;; Now we load the saved reference table

\newread\%reference%fileno
\openin\%reference%fileno=\jobname.ref
{\def\next{\ifeof\%reference%fileno\else\read\%reference%fileno to \nextref\nextref\next\fi}\next}
\closein\%reference%fileno

;; We save the reference table for this run as it is built.  Changes
;; to page numbers so not appear until the next run.

\newwrite\%reference%fileno
\immediate\openout\%reference%fileno=\jobname.ref

;; By changing the definition of \label and adding appropriate ``set''
;; other kinds of references can be kept.  Since the setpageref
;; is inside the write macro, the page is typeset before the
;; \the\pageno is evaluated and so we obtain correct results even in
;; paragraphs in which the page breaks.  To ensure that the other
;;  parts of reference are correct, they are written immediately.
;;
;; Have to be careful that inclusion of \label does not affect
;; spacing. (The only way the debugging message about vertical mode
;; can occur is if \xlabelhook does something funny.  If \xlabelhook
;; leaves vertical mode intentionally, have it call \labelvmodefalse.)
;;
;; The macro \fl%exp gives the first-level expansion.  This keeps (for
;; example) ties from being expanded.
;;
;; The %ref%def keeps track of which labels have been defined (in the
;; macros def%lbl) and on what page and returns warnings when labels
;; are multiply defined. 

\newdimen\saveskip
\newif\iflabelvmode
\newif\iflabelpageonly
\def\label#1{\labelpageonlyfalse\do%label{#1}}
\def\labelpage#1{\labelpageonlytrue\do%label{#1}}
\def\do%label#1{\labelvmodefalse;
\ifvmode\labelvmodetrue\saveskip=\lastskip\vskip-\saveskip\fi;
\%ref%def{#1}{\pagenumber}\xlabelhook{#1}\relax;
\iflabelpageonly\relax;
\else;
\immediate\write\%reference%fileno{\string\setrefnum{#1}{\fl%exp\therefnum}};
\immediate\write\%reference%fileno{\string\setrefname{#1}{\fl%exp\therefname}};
\immediate\write\%reference%fileno{\string\setreftype{#1}{\fl%exp\thereftype}};
\immediate\write\%reference%fileno{\string\setreftext{#1}{\fl%exp\thereftext}};
\fi;
\write\%reference%fileno{\string\setrefpage{#1}{\the\pageno}};
\iflabelvmode
   \ifvmode\vskip\saveskip
   \else\immediate\write16{Unexpectedly forced out of
      vertical mode while creating label #1.}
   \fi
\fi\ignorespaces}

\def\%ref%def#1#2{\expandafter\ifx\csname def%#1\endcsname\relax\expandafter\xdef\csname def%#1\endcsname{#2}\else\immediate\write16{Reference #1 multiply defined (while proc. pp. \csname def%#1\endcsname, \pagenumber).}\fi}

\def\fl%exp#1{\expandafter\strip%macro\meaning#1\end%macro}
\def\strip%macro#1>#2\end%macro{#2}

;; This completes the label and reference macros.  Now for the
;; automatic numbering macros.

;; PROCLAIM

;;    Save old definitions
;;    (Need to erase \proclaim here because it is ``\outer''.
;;     Eventually \old%proclaim will contain the old \proclaim
;;     macro, but (since it is ``\outer''), can't ``\let'' it yet,
;;     so we save the old macro in \%temp instead.)

\let\%temp\proclaim\let\proclaim\relax
\let\old%endproclaim\endproclaim

;;    \oldproclaim calls the old \proclaim
;;    Since the old \proclaim calls ``\proclaim'' we need to
;;    restore the definition to the macro ``\proclaim''.
;;    We use \endproclaim to restore the new \proclaim.
;;    Logically, ``\oldproclaim'' should be ended by
;;    ``\endoldproclaim'', so we define this macro.
;;    (It is harmless to end ``\oldproclaim'' with ``\endproclaim''.)

\outer\def\oldproclaim{\let\proclaim\old%proclaim\proclaim}
\def\endproclaim{\old%endproclaim\let\proclaim\new%proclaim}
\let\endoldproclaim\endproclaim

;;    New \proclaim puts in number and then calls the old \proclaim

\outer\def\new%proclaim#1{\advance\proclaimnum by 1
\xdef\therefnum{\theproclaim};
\gdef\thereftext{#1}
\bgroup\edef\next{\gdef\noexpand\thereftype{\fl%exp\thereftext}};
\lowercase\expandafter{\next}\egroup
\xdef\therefname{\fl%exp\thereftext\noexpand~\therefnum};
\xproclaimhook\relax;
\let\proclaim\old%proclaim\proclaim{#1~\therefnum}}

;;    Finally, we can define ``\old%proclaim'' and ``\proclaim''

\let\old%proclaim\%temp
\let\proclaim\new%proclaim

;; DEFINITION

\let\old%definition\definition
\let\old%enddefinition\enddefinition

\def\olddefinition{\let\definition\old%definition\definition}
\def\enddefinition{\old%enddefinition\let\definition\new%definition}
\let\endolddefinition\enddefinition

\def\new%definition#1{\advance\definitionnum by 1
\xdef\therefnum{\thedefinition};
\gdef\thereftext{#1};
\bgroup\edef\next{\gdef\noexpand\thereftype{\fl%exp\thereftext}};
\lowercase\expandafter{\next}\egroup
\xdef\therefname{\fl%exp\thereftext\noexpand~\therefnum};
\xdefinitionhook\relax;
\let\definition\old%definition\definition{#1~\therefnum}}

\let\definition\new%definition

;; EXAMPLE

\let\old%example\example
\let\old%endexample\endexample

\def\oldexample{\let\definition\old%definition\let\example\old%example\example}
\def\endexample{\old%endexample\let\example\new%example\let\definition\new%definition}
\let\endoldexample\endexample

\def\new%example#1{\advance\examplenum by 1
\xdef\therefnum{\theexample};
\gdef\thereftext{#1};
\bgroup\edef\next{\gdef\noexpand\thereftype{\fl%exp\thereftext}};
\lowercase\expandafter{\next}\egroup
\xdef\therefname{\fl%exp\thereftext\noexpand~\therefnum};
\xexamplehook\relax;
\oldexample{#1~\therefnum}}

\let\example\new%example

;; REMARK
;;   (Also requires fixing ``\demo'' to use oldremark)

\let\old%remark\remark
\let\old%endremark\endremark

\let\olddemo\demo
\let\old%enddemo\enddemo
\def\demo{\let\remark\old%remark\old%demo}
\def\enddemo{\old%enddemo\let\remark\new%remark}
\let\endolddemo\enddemo
\let\old%demo\olddemo

\def\oldremark{\let\remark\old%remark\remark}
\def\endremark{\old%endremark\let\remark\new%remark}
\let\endoldremark\endremark

\def\new%remark#1{\advance\remarknum by 1
\xdef\therefnum{\theremark};
\gdef\thereftext{#1};
\bgroup\edef\next{\gdef\noexpand\thereftype{\fl%exp\thereftext}};
\lowercase\expandafter{\next}\egroup
\xdef\therefname{\fl%exp\thereftext\noexpand~\therefnum};
\xremarkhook\relax;
\oldremark{#1~\therefnum}}

\let\remark\new%remark

;; SECTIONS

;;   The section macro needs to check whether or not the next
;;   character is * -- this is the standard method using \futurelet

\let\%temp\relax
\outer\def\section{\futurelet\%temp\open%section}
\def\open%section#1{\gdef\thereftype{section};
\if\%temp*
\def\%temp##1{\gdef\thereftext{##1}\xsectionstarhook\relax\heading{\thereftext}\endheading}
\gdef\therefnum{}
\gdef\therefname{Section};
\else
\global\advance\sectionnum by 1;
\xdef\therefnum{\thesection}
\xdef\therefname{Section \therefnum};
\gdef\thereftext{#1};
\xsectionhook\relax;
\def\%temp{\heading{\thesection. \thereftext}\endheading}
\resetcounters\fi\%temp}

;;   ``\subsection'' is essentially the same as ``\section'' m.m., but
;;   for some reason the ``\subheading'' macro has a diffent usage
;;   than ``\heading''.  Go figure.

\outer\def\subsection{\futurelet\%temp\open%subsection}
\def\open%subsection#1{\gdef\thereftype{subsection};
\if\%temp*
\def\%temp##1{\gdef\thereftext{##1}\xsubsectionstarhook\relax\subheading{\thereftext}}
\gdef\therefnum{}
\gdef\therefname{Section};
\else
\global\advance\subsectionnum by 1;
\xdef\therefnum{\thesubsection}
\xdef\therefname{Subsection \therefnum};
\gdef\thereftext{#1};
\xsubsectionhook\relax;
\def\%temp{\subheading{\thesubsection. \thereftext}}
\resetsubcounters\fi\%temp}

;; That's it!

;; Restore the characters ``%'' and ``;'' to their usual meanings

\catcode`\%=14
\catcode`\;=12

%% end of autonum.tex

\input xy
\xyoption{matrix}
\xyoption{arrow}
\xyoption{2cell}
\UseTwocells

\def\id{{\mathop{\roman{id}}\nolimits}}
\def\op{{\roman{op}}}
\def\A{{\Cal{A}}}
\def\Al{{\bold{A}}}
\def\AA{\Al^{\Sigma_*}}
\def\AE{\Al^{E\Sigma_*}}
\def\B{{\Cal{B}}}
\def\Bi{{\bold{B}}^{\Sigma_*}}
\def\bMsig{b\M^{\Sigma_*}}
\def\bEMsig{b_E\M^{\Sigma_*}}
\def\C{{\Cal{C}}}
\def\D{{\Cal{D}}}
\def\E{{\Cal{E}}}
\def\F{{\Cal{F}}}
\def\IO#1{\iota_{#1}}
\def\IP{I^{+}}
\def\KP{J^{+}}
\def\M{{\bold{M}}}
\def\MI{\Bbb\IP}
\def\MK{\Bbb\KP}
\def\MM{\Bbb{M}}
\def\MMp{\MM'}
\def\Mp{{\M'}}
\def\Msig{\ell\M^{\Sigma_*}}
\def\N{{\Bbb{N}}}
\def\O{{\Cal{O}}}
\def\OB{O}
\def\OBp{{\OB'}}
\def\P{{\bold{P}}}
\def\PP{{\Bbb{P}}}
\def\Q{{\bold{Q}}}
\def\R{{\Cal{R}}}
\def\rMsig{r\M^{\Sigma_*}}
\def\S{{\Cal{S}}}
\def\SS{{\Cal{SS}}}
\def\T{{\Cal{T}}}
\def\U#1{\Bbb{U}_{#1}}
\def\Cat{{\bold{Cat}}}
\def\ACat{{\bold{C}}}
\def\obj#1{\mathop{\roman{Ob}}\left(#1\right)}
\def\b#1{\bold{#1}}
\def\Deltaop{\Delta^\op}
\def\br#1{\langle#1\rangle}
\def\brn#1{\langle \b{n}\rangle^{#1}}
\def\blank{\underline{\phantom{m}}}
\def\klin{\P_k}
\def\kmap{{\hbox{$k${\rm-map}}}}
\def\Aut{{\mathop{\roman{Aut}}}}
\def\ttheta{{\tilde{\theta}}}
\def\u#1{\underline{#1}}
\def\overto#1{\buildrel{#1}\over{\longrightarrow}}
\let\iso\cong
\let\sma\wedge
\def\barC{\overline{\C}}
\def\SSps#1{S^{#1}_{\bullet}}
\def\SSp{\SSps1}
\def\SSpn{\SSps{n}}
\def\SSpd#1#2{S^{#1}_{#2}}
\def\KSeg{K^{\text{Seg}}}
\def\Knew{K^{\text{new}}}
\def\subst#1{\lceil_{\!#1}}
\def\Colim{\mathop{\roman{Colim}}}
\def\Times{\mathop{\times}\limits}
\def\Prod{\prod\limits}
\def\multprod{\Gamma}
\def\SimpAC{\ACat^{\Deltaop}}
\def\SpecAC{\S_{\ACat}}

\topmatter
\title Rings, modules, and algebras in infinite loop space
theory
\endtitle
\author A.\ D.\ Elmendorf and\\
M.\ A.\ Mandell
\endauthor
\affil Purdue University Calumet\\
The University of Chicago
\endaffil
\address Department of Mathematics,
Purdue University Calumet,
Hammond, IN 46323
\endaddress
\email aelmendo\@math.purdue.edu
\endemail
\address Department of Mathematics,
University of Chicago,
Chicago, IL 60637
\endaddress
\email mandell\@math.uchicago.edu
\endemail
\thanks
The second author was supported in part by NSF grant
DMS-0203980
\endthanks
\date February 6, 2004
\enddate
\abstract We give a new construction of the algebraic $K$-theory
of small permutative categories that preserves multiplicative
structure, and therefore allows us to give a unified treatment of
rings, modules, and algebras in both the input and output.  This
requires us to define multiplicative structure on the category of
small permutative categories.  The framework we use is the concept
of multicategory, a generalization of symmetric monoidal category
that precisely captures the multiplicative structure we have
present at all stages of the construction. Our method ends up in
Smith's category of symmetric spectra, with an intermediate stop
at a new category that may be of interest in its own right, whose
objects we call symmetric functors.
\endabstract
\endtopmatter

\document

\section{Introduction}

This paper offers a new treatment of multiplicative infinite loop
space theory that expands and improves on the account in the
literature.  The motivation comes from the new tools provided by
the modern categories of spectra such as those of \cite{5} and
\cite{7}, which provide cleaner versions of old questions as well
as new ones that could not be asked before.  We now know that any
$E_\infty$ ring spectrum is equivalent to a strictly commutative
ring in any of the new categories of spectra.   It has been known
since the 1980's that the $K$-theory of a bipermutative category
is an $E_\infty$ ring spectrum, although there are gaps in the
proof in the literature which we describe below, and circumvent by
our new methods. The next natural question, asked by Gunnar
Carlsson, is: What structure on a permutative category makes its
$K$-theory into a module over this commutative ring?  We give a
full answer to this question, as well as corresponding ones about
rings, modules, and algebras of all sorts in the context of
permutative categories and their $K$-theory spectra.

Our treatment of multiplicative structures relies on the concept
of {\sl multicategory}, which is an old, familiar friend to
category theorists and computer scientists, but likely foreign to
topologists and $K$-theorists.  It was introduced by Lambek in
1969 in \cite{10}, although without the symmetric group actions we
require. A multicategory is a simultaneous generalization of an
operad and a symmetric monoidal category, and can be thought of as
an ``operad with many objects'' in precisely the same way that a
category can be thought of as a ``monoid with many objects.''
Indeed, an operad is precisely a multicategory with one object.
Any symmetric monoidal category has an underlying multicategory
(more accurately, one for each choice of associating sums, all of
which are canonically isomorphic), but there are many other
multicategories besides these. In particular, restricting to a
subclass of objects in a multicategory again results in a
multicategory, in contrast to what happens with a symmetric
monoidal category.  The natural structure-preserving maps between
multicategories are called {\sl multifunctors}. Every
multicategory has an underlying category, and a multifunctor gives
a functor between underlying categories.

Just as it is often fruitful to consider categories enriched over
a symmetric monoidal category other than sets, so too with
multicategories.  The multicategories we study are all enriched
over either small categories or simplicial sets, and these
enrichments play a crucial role in our theory.  If a multicategory
is enriched over small categories, we also consider it as enriched
over simplicial sets via the nerve construction with no further
comment.

Our use of multicategories in this paper is structural: we
construct a multicategory enriched over small categories whose
objects are the small permutative categories -- we could do so
more generally for symmetric monoidal categories, but to no
additional advantage.  We give a new construction of the
$K$-theory of a small permutative category which gives us an
enriched multifunctor to the symmetric monoidal category of
symmetric spectra constructed in
\cite{7}.  The proof of the following theorem occupies
Sections~\refnum{secperm}--\refnum{sftoss}.

\proclaim{Theorem}\label{num6}
The category of small permutative categories forms a multicategory
enriched over the category of small categories. There is a
multifunctor $K$ from small permutative categories to symmetric
spectra, equivalent to the usual $K$-theory functor, respecting
the enrichment over simplicial sets.
\endproclaim

As a consequence of this theorem, any structure on small
permutative categories captured by a map out of a ``parameter''
multicategory passes directly to $K$-theory spectra.  In the case
of ring structure, the parameter multicategories have only one
object, i.e., they are operads.

We define ring structures on permutative categories in
\refname{secperm} in terms of a second monoidal product and
distributivity maps that satisfy certain coherence relations.  The
noncommutative version we call ``associative'' categories, and the
$E_\infty$ version we call bipermutative categories.  The second of these is the 
generalization for lax morphisms of the usual definition (for 
example, in May \cite{14}); see the discussion preceding 
\refname{num11}, below. We prove the following theorem in
\refname{secassoc} that interprets these structures in terms
operads.

\proclaim{Theorem}\label{num911}
There is an operad $\Sigma_*$ for which an associative structure
(\refname{num9}) on a small permutative category $\A$ determines
and is determined by a multifunctor $\Sigma_*\to\P$ sending the
single object of $\Sigma_*$ to $\A$.  There is an $E_\infty$
operad $E\Sigma_*$ for which a bipermutative structure
(\refname{num11}) on a small permutative category $\R$ determines
and is determined by a multifunctor $E\Sigma_*\to\P$ sending the
single object of $E\Sigma_*$ to $\R$.
\endproclaim

We will see that, as an immediate
consequence of these two theorems, our $K$-theory functor sends
associative categories to ring symmetric spectra and bipermutative
categories are sent to $E_{\infty}$ ring symmetric spectra.  In
\refname{secmodules}, we prove analogous theorems that give
parameter multicategory interpretations of various types of module
structures, defined in terms of a pairing of an associative or
bipermutative category with a small permutative category, and also
algebra structures, defined in terms of certain maps from a
bipermutative category to an associative category.  Again as
immediate consequences of
\refname{num6}, all such ring, module, and algebra structures pass
via $K$-theory to the corresponding structures in the category of
symmetric spectra.

Since we wish our output structures to be as rigid as possible, we
prove a theorem comparing $E_\infty$ versions of rings, modules,
and algebras with their strictly commutative analogues. We do this
by studying model category structures on categories of
multifunctors into the category $\S$ of symmetric spectra.  We
prove the following theorem in \refname{secpmssmodel}.

\proclaim{Theorem}\label{num5}
Suppose $\M$ is a small multicategory enriched over simplicial
sets, and let $\S^\M$ be the category of multifunctors from $\M$
to the category $\S$ of symmetric spectra. There is a simplicial
model structure on $\S^\M$ whose weak equivalences are the
objectwise stable equivalences and whose fibrations are the
objectwise positive stable fibrations of symmetric spectra.
\endproclaim

The map of operads from the $E_{\infty}$ operad $E\Sigma_{*}$
describing bipermutative categories to the one point operad
%$\Cal{C}om$
describing commutative monoids or commutative ring
symmetric spectra is an example of a ``weak equivalence'' of
multicategories, as is the multifunctor from the multicategory
describing modules over $E\Sigma_{*}$ algebras to the
multicategory describing modules over a commutative monoid. (See
\refname{wkequiv} for the general definition of weak equivalence
of multicategories.) We prove the following theorem in
\refname{secpmssfunc}.

\proclaim{Theorem}\label{num5.5}
Let $\M$ and $\Mp$ be small multicategories enriched over
simplicial sets. If $f\colon\M\to\Mp$ is a simplicial
multifunctor, then the induced functor $f^*:\S^{\Mp}\to\S^\M$ is
the right adjoint in a Quillen adjunction.  If in addition $f$ is
a weak equivalence, then the Quillen adjunction is a Quillen
equivalence and therefore induces an equivalence on homotopy
categories.
\endproclaim

As a corollary of this general rectification result, we conclude
that any $E_\infty$ ring in symmetric spectra is equivalent to a
strictly commutative ring spectrum (as was already well-known),
but also that any $E_\infty$ module over an $E_\infty$ ring is
equivalent to a strict module over an equivalent commutative ring,
as well as a wide range of similar results for many other
structures.

The need to use a multicategory structure on small permutative
categories rather than a symmetric monoidal structure seems
intrinsic: contrary to Thomason's claim in the introduction to
\cite{18}, small permutative categories appear not to support a
symmetric monoidal structure consistent with a reasonable notion
of multiplicative structure. We will explain in a later paper how
this problem can be resolved by embedding into a larger symmetric
monoidal category (whose objects are, ironically,
multicategories), but the necessary complications are irrelevant
to the present paper.

On a technical note, our construction of the $K$-theory
multifunctor is actually a two step process, with an intermediate
stop at a new multicategory which may be of interest in its own
right; we call the objects {\sl symmetric functors}.  They are
described in \refname{secsym}.

Historically, the question of what additional structure to impose on a
permutative, or more generally a symmetric monoidal category in
order to give its $K$-theory some sort of ring structure was first
investigated by Peter May in \cite{14}.  He defined bipermutative
categories, and offered a proof that their $K$-theory spectra are
$E_\infty$ ring spectra.  Unfortunately, he made a combinatorial
error (found by Steinberger), as explained in Appendix A of
\cite{16}. This led May to write
\cite{16}, whose main results are entirely correct. However, there is a
further combinatorial error in \cite{16}, Section 7, which was
patched by Uwe Hommel; unfortunately, the patch was never
published. Gerry Dunn also found an error in the category theory
in Section 4 of \cite{16}, which he described and attempted to
patch in
\cite{2}. However, there is a critical error in
\cite{2}, Section 2 (the evaluation $\xi$ of Lemma~2.2(ii) is not
well-defined).  The categorical error in \cite{16} can apparently
be fixed by making a correction to the left adjoints, although a
detailed check has yet to be made.  One benefit of the current paper            
is to give a new proof of
this theorem. Since there were no reasonable concepts of module
and algebra spectra available at the time \cite{16} was written,
the question of which permutative categories give rise to these
sorts of $K$-theory spectra was not addressed; we do so now.

The paper is organized as follows:  Section~2 contains a precise
definition of multicategory and a description of types of
parameter multicategories giving ring, module, and algebra
structures. Section~3 constructs the multicategory structure on
the category of small permutative categories and describes our
results on ring structure in greater detail. In Section~4, we
recall the construction of the $K$-theory of a permutative
category in the literature, give our new construction as a functor
(as opposed to a multifunctor), and prove that our construction is
equivalent to the old one. Section~5 is devoted to the description
of
%the intermediate stop in our construction,
the multicategory of symmetric functors.
Section~6 constructs the multifunctor from
%then gives the part of our $K$-theory multifunctor that goes from
permutative categories to symmetric functors, and Section~7 constructs
the multifunctor from symmetric functors to symmetric spectra; the
composite of these two is our $K$-theory multifunctor.
Section~8 proves \refname{num911}, describing associative categories
and bipermutative categories in terms of actions of the operads
$\Sigma_*$ and $E\Sigma_*$. Section~9 describes
the various sorts of modules and algebras in
permutative categories
in terms of parameter
multicategories. Section~10 describes various ways in which
free permutative categories can have associative or bipermutative
structure. Finally, Sections 11 and 12 contain the proofs of our model
category results, Theorems~\refnum{num5} and~\refnum{num5.5}.

The first author would like to thank Gunnar Carlsson for asking some
very interesting questions, and Peter May for both encouragement and criticism.

\section{Multicategories}\label{secmult}

\definition{Definition}\label{num1}  A {\bf multicategory} $\M$
consists of the following:
\roster
\item  A collection of objects (which may form a
proper class)
\item For each $k\ge
0$, $k$-tuple of objects $(a_1,\ldots,a_k)$ (the ``source'') and
single object $b$ (the ``target''), a set $\M_k(a_1,\ldots,a_k;b)$
(the ``$k$-morphisms'')
\item
A right action of $\Sigma_k$ on the collection of all
$k$-morphisms, where for $\sigma\in\Sigma_k$,
$$\sigma^*:\M_k(a_1,\ldots,a_k;b)
\to\M_k(a_{\sigma(1)},\ldots,a_{\sigma(k)};b)
$$
\item
A distinguished ``unit'' element $1_a\in\M_1(a;a)$
for each object $a$, and
\item A composition ``multiproduct''
$$\eqalign{\multprod:&\;\M_n(b_1,\ldots,b_n;c)\times
\M_{k_1}(a_{11},\ldots,a_{1k_1};b_1)\times
\cdots\times\M_{k_n}(a_{n1},\ldots,a_{nk_n};b_n)\cr
&\longrightarrow\M_{k_1+\cdots+k_n}(a_{11},\ldots,a_{nk_n};c).}
$$
\endroster
subject to the identities for an operad listed on pages 1--2 in
\cite{13}, which still make perfect sense in this context.   In
greater detail, we require the diagrams (1)--(4) below to commute
for all nonnegative integers $k$, $j_s$ for $1\le s\le k$, and
$i_{sq}$ for $1\le q\le j_s$, and all objects $d$,  $c_s$ for
$1\le s\le k$,  $b_{sq}$ for $1\le s\le k$ and $1\le q\le j_s$,
and  $a_{sqp}$ for $1\le s\le k$, $1\le q\le j_s$, and $1\le p\le
i_{sq}$. In these diagrams,  we write $i_s$ for
$\sum_{q=1}^{j_s}i_{sq}$, $i$ for $\sum_{s=1}^k i_s$, and $j$ for
$\sum_{s=1}^k j_s$, and to compress the diagrams to fit on the
page, we write lists like $c_{1},\dotsc,c_{k}$ as $\br{c}$ or as
$\br{c_{s}}_{s=1}^{k}$ when the index is ambiguous.
\roster
\item We require the following
multiassociativity diagram to commute.
$$\xymatrix@C=-122pt @R-10pt{
&\M_k(\br{c};d)\times\Prod_{s=1}^k
\M_{i_s}(\br{\br{a_{sqp}}_{p=1}^{i_{sq}}}_{q=1}^{j_s};c_s)
\ar[ddr]^-{\multprod}
\\ \M_k(\br{c};d)\times\Prod_{s=1}^k
\left(\M_{j_s}(\br{b_{sq}}_{q=1}^{j_s};c_s)
\times\Prod_{q=1}^{j_s}\M_{i_{sq}}
(\br{a_{sqp}}_{p=1}^{i_{sq}};b_{sq})\right)
\ar[ur]^-{\id\times\multprod}\ar[dd]_-{\cong}
\\&&\M_i(\br{\br{\br{a_{sqp}}_{p=1}^{i_{sq}}}_{
q=1}^{j_s}}_{s=1}^k;d).
\\ \M_k(\br{c};d)\times\Prod_{s=1}^k
\M_{j_s}(\br{b_{sq}}_{q=1}^{j_s};c_s)
\times\Prod_{s=1}^k\ar[dr]_-{\multprod\times1}
\Prod_{q=1}^{j_s}\M_{i_{sq}}
(\br{a_{sqp}}_{p=1}^{i_{sq}};b_{sq})
\\& \M_j(\br{\br{b_{sq}}_{q=1}^{j_s}}_{s=1}^k;d)
\times\Prod_{s=1}^k
\Prod_{q=1}^{j_s}\M_{i_{sq}}
(\br{a_{sqp}}_{p=1}^{i_{sq}};b_{sq})
\ar[uur]_-{\multprod}}
$$
\item We require the following unit diagrams to commute:
$$\xymatrix@C=1pt{
\M_k(\br{c};d)\times\{1\}^k\ar[r]^-{\cong}
\ar[d]_-{\id\times1^k}
&\M_k(\br{c};d),
\\ \M_k(\br{c};d)\times\Prod_{s=1}^k\M_1(c_s;c_s)
\ar[ur]_-{\multprod}}
\xymatrix@C=7pt{
\{1\}\times\M_k(\br{c};d)\ar[r]^-{\cong}
\ar[d]_-{1\times\id}
&\M_k(\br{c};d).
\\ \M_1(d;d)\times\M_k(\br{c};d)\ar[ur]_-{\multprod}}
$$
\item Given $\sigma\in\Sigma_k$, we require the following
equivariance
diagram to commute:
$$\xymatrix{
\M_k(\br{c};d)\times\Prod_{s=1}^k
\M_{j_s}(\br{b_{sq}}_{q=1}^{j_s};c_s)
\ar[r]^-{\multprod}\ar[d]_-{\sigma\times\sigma^{-1}}
&\M_j(\br{\br{b_{sq}}_{q=1}^{j_s}}_{s=1}^k;d)
\ar[d]^-{\sigma\br{j_{\sigma(1)},\ldots,j_{\sigma(k)}}}
\\ \M_k(\br{c_{\sigma(s)}}_{s=1}^{k};d)\times\Prod_{s=1}^k
\M_{j_{\sigma(s)}}
(\br{b_{{\sigma(s)}q}}_{q=1}^{j_{\sigma(s)}};c_{\sigma(s)})
\ar[r]_-{\multprod}
&\M_j(\br{\br{b_{\sigma(s)q}}_{
q=1}^{j_{\sigma(s)}}}_{s=1}^k;d),}
$$
where $\sigma\br{j_{\sigma(1)},\ldots,j_{\sigma(k)}}$ permutes
blocks as indicated.
\item Given $\tau_s\in\Sigma_{j_s}$ for $1\le s\le
k$, we require the following equivariance diagram to commute:
$$\xymatrix{
\M_k(\br{c};d)\times\Prod_{s=1}^k
\M_{j_s}(\br{b_{sq}}_{q=1}^{j_s};c_s)
\ar[r]^-{\multprod}\ar[d]_-{\id\times\Prod\tau_s}
&\M_j(\br{\br{b_{sq}}_{q=1}^{j_s}}_{s=1}^k;d)
\ar[d]^-{\tau_1\oplus\cdots\oplus\tau_k}
\\ \M_k(\br{c};d)\times\Prod_{s=1}^k
\M_{j_s}(\br{b_{s\tau(q)}}_{q=1}^{j_s};c_s)
\ar[r]_-{\multprod}
&\M_j(\br{\br{b_{s\tau(q)}}_{q=1}^{j_s}}_{s=1}^k;d).}
$$
\endroster
This concludes the definition of a multicategory. However, we may
also ask that the $k$-morphisms $\M_k(a_1,\ldots,a_k;b)$ take
values in a symmetric monoidal category other than sets; the
examples we are interested in take values in either categories or
simplicial sets. This gives the concept of an {\bf enriched}
multicategory.  Note that a multicategory enriched over small
categories can be considered enriched over simplicial sets by
applying the nerve functor to the $k$-morphisms, since the nerve
functor preserves categorical products.
\enddefinition

\definition{Definition}
For multicategories $\M$ and $\Mp$, a {\bf multifunctor} from $\M$
to $\Mp$ consists of a function $f$ from the objects of $\M$ to
the objects of $\Mp$, and for all objects $b$ and $k$-tuples of
objects $a_{1},\dotsc, a_{k}$, a function
$\M_k(a_{1},\dotsc,a_{k};b)\to
\M'_k(f(a_{1}),\dotsc,f(a_{k});f(b))$ which preserves the
$\sigma_{k}$ action on the collection of all $k$-morphisms,
preserves the units, and preserves the multiproduct.  When $\M$
and $\Mp$ are enriched over simplicial sets or small categories,
the multifunctor is enriched when the maps on $k$-morphisms
preserve the enrichment; in this context, ``multifunctor'' always
means enriched multifunctor.
\enddefinition

\oldexample{Example}  In any symmetric monoidal
category $(\M,\oplus,0)$,
we can define $k$-morphisms as
$\M_k(a_1,\ldots,a_k;b):=\M(a_1\oplus\cdots\oplus a_k,b)$, with
the sums associated in any fixed order.
\endoldexample

\oldexample{Example} An operad is simply
a multicategory with one
object.
\endoldexample

\oldremark{\bf Remark} If we restrict our attention just to the
objects and 1-morphisms of a multicategory, we get a category.
\endoldremark

A major theme of this paper is that rings, modules, and algebras
can be described in any multicategory, and as we shall see in
\refname{secassoc}, the enrichments present in our examples
of interest allow for $E_\infty$ versions of these concepts as
well.  These are all described by means of maps out of what we
call {\sl parameter multicategories}, which are simply specific,
very small examples of multicategories.  Since our construction of
the $K$-theory of a small permutative category is a multifunctor,
it follows automatically that ring, module, and algebra structures
on small permutative categories are preserved in their $K$-theory
spectra.  We turn next to descriptions of our basic classes of
parameter multicategories.

\definition{Definition}\label{num2}  Let $\O$ be an operad (a multicategory
with only one object) and $\Q$ a multicategory. An $\O${\bf-ring}
in $\Q$ is a multifunctor from $\O$ to $\Q$. Usually we speak of
the target object in $\Q$ as being the ring.  If the morphism
spaces of $\O$ are all contractible, then we say that the target
object is an $E_\infty$ ring.
\enddefinition

For example, if $\O$ is the final operad with $\O_k=*$ for all
$k$, then an $\O$-ring in a symmetric monoidal category is simply
a commutative monoid in that category.  In particular, if the target category is
abelian groups under tensor product, an $\O$-ring is simply a
commutative ring.  As another example, if $\O=\Sigma_*$ is the
``associative'' operad with $\O_k=\Sigma_k$ (described in greater
detail below), then an $\O$-ring in a symmetric monoidal category
is a monoid in the underlying monoidal category. In the case of
abelian groups, we get a ring.

We also define parameter multicategories for modules and algebras.

\definition{Definition}\label{num3}  Let $\M$ be a multicategory with two
objects, $R$ (the ``ring'') and $M$ (the ``module'').  We say that
$\M$ is a {\bf parameter multicategory for modules} if we have
$\M_k(B_1,\ldots,B_k;C)=\emptyset$ unless all variables are $R$, or else $C$ and exactly one of the $B$'s are $M$.  If all the nonempty morphism
spaces are contractible, then we say that $\M$ is a parameter
multicategory for $E_\infty$ modules.
\enddefinition

In the special case where $\M_k(B_1,\ldots,B_k;C)=*$ whenever it
is not required to be empty, we find that a multifunctor into a
symmetric monoidal category consists of a commutative monoid (the
image of $R$) and an action of that monoid on another object (the
image of $M$).  In the special case of abelian groups, we get a
commutative ring and a module over it.

As another example, if $\O$ is an operad, we can let
$\M_k(B_1,\ldots,B_k;C)=\O_k$ whenever it is not required to be
empty.  This recovers the notion of $\O$-module defined by
Ginzburg and Kapranov in \cite{6} and discussed by Kriz and May in
Section~I.4 of \cite{9}.  In particular, if $\O=\Sigma_*$, we get a
monoid and a ``bimodule'' (which has commuting left and right
actions).

For a third example, we let
$\M_k(R^{j-1},M,R^{k-j};M)=\{\sigma\in\Sigma_k:\sigma(j)=k\}$ for
all $j$, $\M_k(R^k;R)=\Sigma_k$, and we make all other morphism
sets empty.  Then a multifunctor from $\M$ to a symmetric monoidal
category is a (noncommutative) monoid and a left action on another
object of the category.  If instead we make
$\M_k(R^{j-1},M,R^{k-j};M)=\{\sigma\in\Sigma_k:\sigma(j)=1\}$,
then we get a right action.

Next we turn to algebra structures.

\definition{Definition}\label{num4}  A parameter multicategory for algebras
is a multicategory $\Al$ with two objects, $R$ (the ``ring'') and
$A$ (the ``algebra''), subject to the following condition. Suppose
given inputs $B_1,\ldots,B_k$ with at least one of the $B_j$'s being
equal to $A$.  Then we require that
$\Al_k(B_1,\ldots,B_k;R)=\emptyset$. If all the other $k$-morphism
spaces are contractible, then we say that $\Al$ is a parameter
multicategory for $E_\infty$ algebras.
\enddefinition

Again, we can look at the example in which all the nonempty
$k$-morphism spaces are a single point, and we map to a symmetric
monoidal category. Then the images of both $R$ and $A$ are
commutative monoids, and the rest of the structure is induced by a
strict map of monoids from $R$ to $A$ given by the single element
of $\Al_1(R;A)$.

A more interesting example is given by letting $S=\{j:B_j=A\}$ in
the expression $\Al_k(B_1,\ldots,B_k;C)$ and, if not required to
be empty, setting this $k$-morphism space equal to
$\Sigma_k/\sim$, where $\sim$ is the equivalence relation on
$\Sigma_k$ given by requiring $\sigma\sim\sigma'$ if and only if,
for all elements $i$ and $j$ of $S$,
$\sigma(i)<\sigma(j)\Leftrightarrow\sigma'(i)<\sigma'(j)$.  Then a
multifunctor to a symmetric monoidal category makes the image of
$R$ again a commutative monoid, the image of $A$ is now a
noncommutative monoid, and the map induced by the single element
of $\Al_1(R;A)$ is central in the obvious sense.

For a third example, let $\O$ be an operad.  Then we can let
$\Al_k(B_1,\ldots,B_k;C)=\O_k$ whenever it is not required to be
empty.  Then the images of both $R$ and $A$ are $\O$-rings, and
there is a map of $\O$-rings given by the identity element of
$\O_1=\Al_1(R;A)$ which determines the entire algebra structure.

We describe further variants of module and algebra structures and
their applications to permutative categories in
\refname{secmodules}.

\section{The Multicategory of Permutative Categories}\label{secperm}

In this section we describe the multicategory of permutative
categories.  We begin by recalling the definition of permutative
category.

\definition{Definition}\label{num7}  A {\bf permutative category}\/
is a category $\C$ with a functor $\oplus\colon \C\times\C\to\C$,
an object $0\in\obj\C$, and a natural isomorphism $\gamma\colon
a\oplus b\cong b\oplus a$ satisfying:
\roster
\item $(a\oplus b)\oplus c=a\oplus(b\oplus c)$ (strict associativity),
\item $a\oplus0=a=0\oplus a$ (strict unit),
\item The following three diagrams must commute:
$$\xymatrix{a\oplus 0\ar[rr]^\gamma_\cong\ar[dr]_=&&0\oplus a\ar[dl]^=\\
&a,\\}
\qquad\qquad
\xymatrix{a\oplus b\ar[rr]^=\ar[dr]^\cong_\gamma&&a\oplus b\\
&b\oplus a,\ar[ur]^\cong_\gamma\\}
$$
$$\xymatrix{a\oplus b\oplus c\ar[rr]^\gamma\ar[dr]_{1\oplus\gamma}&&c\oplus
a\oplus b\\ &a\oplus c\oplus b.\ar[ur]_{\gamma\oplus 1}\\}
$$
\endroster
A permutative category is {\bf small}\/ if its underlying category
is small.
\enddefinition

Any symmetric monoidal category is naturally equivalent to a
permutative category by a well-known theorem of Isbell \cite{8}.
We also have the following examples of small permutative
categories from $K$-theory.

\oldexample{Examples} Let $A$ be a ring and let $\roman{GL} A$ be the
category whose objects are the standard free modules $A^{n}$ and whose
morphisms are the (left) $A$-module isomorphisms.  Direct sum makes
$\roman{GL}A$ into a small permutative category, whose $K$-theory is the
``free module'' algebraic $K$-theory of $A$.  More generally, let
$\roman{Pr}A$ be the following category.  An object is a pair
$(A^{n},i)$ where $i\colon A^{n}\to A^{n}$ is an idempotent left
$A$-module endomorphism.  A map from $(A^{m},i)$ to $(A^{n},j)$ is a
left $A$-module isomorphism from $\mathop{\roman{Im}}(i)$ to
$\mathop{\roman{Im}}(j)$. Again, direct sum (of modules and
idempotents) makes $\roman{Pr}A$ a small permutative category.  The
$K$-theory of $\roman{Pr}A$ is the algebraic $K$-theory of the ring
$A$. The functor $\roman{GL}A\to \roman{Pr}A$ that sends $A^{n}$ to
$(A^{n},\id)$ induces a map on $K$-theory that is an isomorphism on
homotopy groups in all degrees except (possibly) degree
zero. \endoldexample

The following definition describes the multicategory we study
whose objects are the small permutative categories.

\definition{Definition}\label{num8} Let $\C_1,\ldots,\C_k$ and $\D$
be small permutative categories.  We define categories
$\klin(\C_1,\ldots,\C_k;\D)$ that provide the categories of
$k$-morphisms for the
multicategory $\P$ of permutative categories.
The objects of $\klin(\C_1,\ldots,\C_k;\D)$ consist of
functors
$$f\colon \C_1\times\cdots\times\C_k\to\D
$$
which we think of as $k$-linear maps, satisfying
$f(c_1,\ldots,c_k)=0$ if any of the $c_i$ are 0, together with
natural transformations, which we think of as distributivity maps,
$$\delta_i\colon f(c_1,\ldots,c_i,\ldots,c_k)\oplus
f(c_1,\ldots,c_i',\ldots,c_k)\to f(c_1,\ldots,c_i\oplus
c_i',\ldots,c_k)
$$
for $1\le i\le k$.  We conventionally suppress the variables
that do not change, writing
$$\delta_i\colon f(c_i)\oplus f(c_i')\to f(c_i\oplus c_i').
$$
We require $\delta_i=\id$ if either $c_i$ or $c_i'$ is 0, or if
any of the other $c_j$'s are 0. These natural transformations are
subject to the commutativity of the following diagrams:
$$ \xymatrix{
f(c_i)\oplus f(c_i')\ar[d]_-{\gamma}^-{\cong} \ar[r]^-{\delta_i}
&f(c_i\oplus c_i') \ar[d]^-{f(\gamma)}_-{\cong}
\\
f(c_i')\oplus f(c_i)\ar[r]_-{\delta_i} &f(c_i'\oplus c_i),
\\}
\qquad\qquad
\xymatrix{
f(c_i)\oplus f(c_i')\oplus f(c_i'')\ar[d]_-{\delta_i\oplus 1}
\ar[r]^-{1\oplus\delta_i} &f(c_i)\oplus f(c_i'\oplus c_i'')
\ar[d]^-{\delta_i}
\\
f(c_i\oplus c_i')\oplus f(c_i'')\ar[r]_-{\delta_i} &f(c_i\oplus
c_i'\oplus c_i''),\\}
$$
and for $i\ne j$,
$$\xymatrix@C=-2pc{&f(c_i\oplus c_i',c_j)\oplus f(c_i\oplus c_i',c_j')
\ar[ddr]^-{\delta_j}
\\
f(c_i,c_j)\oplus f(c_i',c_j)\oplus f(c_i,c_j')\oplus f(c_i',c_j')
\ar[ur]^-{\delta_i\oplus\delta_i}
\ar[dd]_-{1\oplus\gamma\oplus1}^{\cong}
\\
&&f(c_i\oplus c_i',c_j\oplus c_j').
\\
f(c_i,c_j)\oplus f(c_i,c_j')\oplus f(c_i',c_j)\oplus f(c_i',c_j')
\ar[dr]_{\delta_j\oplus\delta_j}
\\
&f(c_i,c_j\oplus c_j')\oplus f(c_i',c_j\oplus c_j')
\ar[uur]_-{\delta_i}}
$$
This completes the definition of the objects of
$\klin(\C_1,\ldots,\C_k;\D)$.  To specify its morphisms,
given two objects $f$ and $g$, a morphism $\phi\colon f\to g$ is a
natural transformation commuting with all the $\delta_i$'s, in the
sense that all the diagrams
$$\xymatrix{f(c_i)\oplus f(c_i')\ar[d]_-{\phi\oplus\phi}
\ar[r]^-{\delta_i^f} & f(c_i\oplus c_i') \ar[d]^-{\phi}
\\
g(c_i)\oplus g(c_i') \ar[r]_-{\delta_i^g} &g(c_i\oplus c_i')
\\}$$
commute.

In order to make the $\klin(\C_1,\ldots,\C_k;\D)$'s the
$k$-morphisms of a multicategory, we must specify a $\Sigma_k$
action and a multiproduct.  The $\Sigma_k$ action
$$\sigma^*f\colon \C_{\sigma(1)}\times\cdots\times\C_{\sigma(k)}
\to\D
$$
is specified by
$$\sigma^*f(c_{\sigma(1)},\ldots,c_{\sigma(k)})
=f(c_1,\ldots,c_k),
$$
with the structure maps $\delta_i$ inherited from $f$ (with the
appropriate permutation of the indices).  We define the
multiproduct as follows: Given $f_j\colon
\C_{j1}\times\cdots\times\C_{jk_j}\to \D_j$ for $1\le j\le n$ and
$g\colon \D_1\times\cdots\times\D_n\to\E$, we define
$$\multprod(g;f_1,\ldots,f_n):=g\circ(f_1\times\cdots\times f_n).$$
To specify the structure maps, suppose $k_1+\cdots+k_{j-1}<s\le
k_1+\cdots+k_j$, and let $i=s-(k_1+\cdots+k_{j-1})$.  Then
$\delta_s$ is given by the composite
$$\xymatrix{
g(f_j(c_{ji}))\oplus g(f_j(c'_{ji}))\ar[r]^-{\delta^g_j}
&g(f_j(c_{ji})\oplus f_j(c'_{ji}))\ar[r]^-{g(\delta^{f_j}_i)}
&g(f_j(c_{ji}\oplus c'_{ji})).
}$$
The authors have checked that these definitions satisfy the
required properties of the structure maps $\delta_s$, and the
diligent reader will do so as well; the pentagonal diagram for the
last structure map has two cases. These definitions extend easily to
morphisms, and we leave to the reader the
straightforward task of checking that the necessary identities for
a multicategory are satisfied.
\enddefinition

\oldremark{\bf Remark} The morphisms of the category of permutative categories that we get by remembering
only the 1-morphisms are called {\bf lax} maps.  To describe them
explicitly, suppose $\C$ and $\D$ are permutative categories. Then
a lax map $f\colon \C\to\D$ consists of a functor on the underlying
categories for which $f(0)=0$, together with a natural
transformation $\lambda\colon f(c)\oplus f(c')\to f(c\oplus c')$.  We
require $\lambda=\id$ if either $c$ or $c'$ is $0$, together with
the commutativity of the first two diagrams in \refname{num8}; the
third diagram does not apply in this situation.  The reader can
now supply the definition of composition of lax maps.
\endoldremark

\oldremark{\bf Variant} A {\bf strong} map of permutative categories
is a lax map for which the natural transformation $\lambda$ of the
previous remark is a natural isomorphism.  When we require the
distributivity transformations $\delta_{i}$ of the previous
definition to be isomorphisms, we obtain a multicategory structure
whose underlying category is the category of strong maps of
small permutative categories. \endoldremark

In the rest of this section, we describe the analogues of rings
and commutative rings that appear to be most useful in the context
of permutative categories, and give some examples.
We begin with the definition of associative
category.  This is the analogue in permutative categories of an
associative ring with unit.

\definition{Definition}\label{num9} An {\bf associative} category is a
permutative category $\A$ together with a functor
$\otimes\colon \A\times\A\to\A$ that is strictly associative with a
strict unit object 1, and natural {\bf distributivity} maps
$$d_l\colon (a\otimes b)\oplus(a'\otimes b)\to(a\oplus
a')\otimes b
$$
and
$$d_r\colon (a\otimes b)\oplus(a\otimes b')\to a\otimes(b\oplus b'),
$$
subject to the following requirements:
\roster
\item"(a)" $a\otimes0=0\otimes a=0$ for all $a$.

\item"(b)" The following diagram commutes, as does an analogous one for $d_r$:
$$\xymatrix{
(a\otimes b)\oplus(a'\otimes b)\oplus(a''\otimes
b)\ar[r]^-{d_l\oplus1}\ar[d]_-{1\oplus d_l}
&((a\oplus a')\otimes b)\oplus(a''\otimes b)\ar[d]^-{d_l}\\
(a\otimes b)\oplus((a'\oplus a'')\otimes b)\ar[r]_-{d_l}
&(a\oplus a'\oplus a'')\otimes b.\\
}$$

\item"(c)" The following diagram commutes, as does an analogous one for $d_r$:
$$\xymatrix{
(a\otimes b)\oplus(a'\otimes
b)\ar[r]^-{d_l}\ar[d]_-{\gamma^\oplus}
&(a\oplus a')\otimes b\ar[d]^-{\gamma^\oplus\otimes1}\\
(a'\otimes b)\oplus(a\otimes b)\ar[r]_-{d_l}&(a'\oplus a)\otimes
b. }$$

\item"(d)" The following diagram commutes, as does an analogous one for $d_r$:
$$\xymatrix{
(a\otimes b\otimes c)\oplus(a'\otimes b\otimes c)\ar[dr]^-{d_l}
\ar[d]_-{d_l}\\
((a\otimes b)\oplus(a'\otimes b))\otimes c\ar[r]_-{d_l\otimes1}
&(a\oplus a')\otimes b\otimes c\\
}$$

\item"(e)" The following diagram commutes:
$$\xymatrix{
(a\otimes b\otimes c)\oplus(a\otimes b'\otimes c)
\ar[r]^-{d_l}\ar[d]_-{d_r}
&((a\otimes b)\oplus(a\otimes b'))\otimes c\ar[d]^-{d_r\otimes1}
\\a\otimes((b\otimes c)\oplus(b'\otimes c))\ar[r]_-{1\otimes d_l}
&a\otimes(b\oplus b')\otimes c.}
$$

\item"(f)" The following diagram commutes:
$$\xymatrix@C-65pt @R-6pt{
&(a\otimes(b\oplus b'))\oplus(a'\otimes(b\oplus b'))\ar[ddr]^-{d_l}\\
(a\otimes b)\oplus(a\otimes b')\oplus(a'\otimes b)\oplus(a'\otimes
b')\ar[ur]^-{d_r\oplus d_r}
\ar[dd]_-{1\oplus\gamma\oplus1}\\
&&(a\oplus a')\otimes(b\oplus b').\\
(a\otimes b)\oplus(a'\otimes b)\oplus(a\otimes b')\oplus(a'\otimes
b')\ar[dr]_-{d_l\oplus
d_l}\\
&((a\oplus a')\otimes b)\oplus((a\oplus a')\otimes b')\ar[uur]_-{d_r}}
$$
\endroster
\enddefinition

\oldexample{Example}  The primary examples of associative categories are
categories of endomorphisms of small permutative categories.  Let $\C$
be a small permutative category.  Then we can give the category of lax
maps $\P_1(\C;\C)$ the structure of an associative category as
follows.  Suppose we have two objects $f$ and $g$, i.e.,
lax maps from $\C$ to itself.  We define $f\oplus g$ as the lax
map for which $(f\oplus g)(c):=fc\oplus gc$, with lax structure
map given by the composite
$$\xymatrix@R=8pt @C=-15pt{
(f\oplus g)(c)\oplus(f\oplus g)(c')\ar[r]^-{=}&fc\oplus gc\oplus
fc'\oplus gc'
\ar[r]^-{\gamma}_-{\cong}
&fc\oplus fc'\oplus gc\oplus gc'
\\ \relax\ar[r]^-{\lambda_f\oplus\lambda_g}
&f(c\oplus c')\oplus g(c\oplus c')=(f\oplus g)(c\oplus c').}
$$
(Notice that even if both lax structure maps $\lambda_f$ and
$\lambda_g$ were the identity, the lax structure map for $f\oplus
g$ would still involve the transposition isomorphism.)  This gives
us permutative structure on $\P_1(\C;\C)$.  The associative structure
is given by composition of lax maps; we leave the necessary
verifications to the reader.
\endoldexample

\oldexample{Example}
If $\C$ is a small monoidal category with a strictly associative and
unital monoidal
product, then the ``free permutative category'' on $\C$ is
functorially an associative category, in fact, in uncountably many
ways.  See \refname{freepcs} for details.
\endoldexample

As further motivation for the definition of the multicategory
structure on permutative categories, we offer the following
theorem, proved in \refname{secassoc}.  The operad $\Sigma_*$ mentioned in the theorem is discussed immediately below.

\proclaim{Theorem}\label{num10} An associative structure on a
small permutative category
$\A$ determines and is determined by a multifunctor
$\Sigma_*\to\P$ sending the single object of $\Sigma_*$ to $\A$.
\endproclaim

Here, as above, $\Sigma_*$ denotes the fundamental
``associative'' operad of sets whose algebras are the associative
monoids.  For convenience, we recall the definition.  The component
sets of $\Sigma_{*}$ are the symmetric groups $\Sigma_k$ and the
multiproduct is described as follows: Let $\sigma\in\Sigma_k$,
$\phi_i\in\Sigma_{j_i}$ for $1\le i\le k$. Then we must have
$\multprod(\sigma;\phi_1,\ldots,\phi_k)\in\Sigma_j$, where
$j=j_1+\cdots+j_k$.  This is specified as the composite
$$\xymatrix{j_1\coprod\cdots\coprod j_k\ar[r]^-{\coprod_i\phi_i}
&j_1\coprod\cdots\coprod j_k\ar[rr]^-{\sigma\langle j_1,\ldots,j_k\rangle}
&&j_{\sigma^{-1}(1)}\coprod\cdots\coprod j_{\sigma^{-1}(k)},}
$$
where $\sigma\langle j_1,\ldots,j_k\rangle$ permutes the blocks
$j_1,\ldots,j_k$ as indicated.  The right action of $\Sigma_k$ is
simply right multiplication.

Since the algebras for the operad $\Sigma_*$ in any symmetric
monoidal category are simply the monoids in the underlying
monoidal category, \refname{num6} now implies the following
corollary.

\proclaim{Corollary} If $\A$ is an associative category, then
$K\A$ is a strict ring symmetric spectrum.
\endproclaim

We next consider commutativity in multiplication, which cannot be
strict in our context; we must settle for $E_{\infty}$.  To
describe the relevant $E_\infty$ operad, we need the following
construction.  Consider the forgetful functor from small
categories to sets that forgets the morphisms and remembers only
the objects. This functor has a right adjoint $E$ that takes a set
$X$ and produces the category $EX$ with $X$ as its set of objects,
and with exactly one morphism between each pair of objects;
formally, the morphism set is $X\times X$. We use $E$ for this
construction because if the set is actually a group $G$, the
classifying space of the category $EG$ is the usual construction
of the universal principal $G$-bundle.  Since $E$ is a right
adjoint, it preserves products, and therefore if $\O$ is any
operad of sets, $E\O$ is an operad of categories. Applying $E$ to
the operad $\Sigma_{*}$ defines the categorical Barratt-Eccles
operad $E\Sigma_*$.  Since $\Sigma_*$ is $\Sigma$-free, so is
$E\Sigma_*$, and $EX$ is always contractible. The structures in
$\P$ induced by $E\Sigma_*$ turn out to be bipermutative
categories, as defined below.
We note that our
bipermutative categories are more general than
May's (\cite{14}, p.~154) both in requiring only distributivity morphisms
rather than isomorphisms, and in deleting the requirement that one
of the distributivity morphisms be the identity.  Laplaza's
symmetric bimonoidal categories \cite{11} are more general even
than our bipermutative categories, and since they can be rectified
to equivalent bipermutative categories in May's sense, so can
ours.  Our
explicit definition is as follows:

\definition{Definition}\label{num11} A {\bf bipermutative category} is a
permutative category $(\R,\oplus,0)$ together with a second
permutative structure $(\R,\otimes,1)$ with symmetry isomorphism
$\gamma^\otimes:a\otimes b\cong b\otimes a$, and natural
distributivity maps
$$d_l\colon (a\otimes b)\oplus(a'\otimes b)\to(a\oplus a')\otimes b
$$
and
$$d_r\colon (a\otimes b)\oplus(a\otimes b')\to a\otimes(b\oplus b').
$$
These are subject to the requirement that the diagrams for an
associative category given in \refname{num9} commute, except with
diagram (e) replaced with the following diagram (e$'$):
$$\xymatrix{
(a\otimes b)\oplus(a'\otimes
b)\ar[r]^-{d_l}\ar[d]_-{\gamma^\otimes\oplus\gamma^\otimes}
&(a\oplus a')\otimes
b\ar[d]^-{\gamma^\otimes}\\
(b\otimes a)\oplus(b\otimes a')\ar[r]_-{d_r}&b\otimes(a\oplus a').\\
}
$$
\enddefinition

\oldexample{Example} Let $A$ be a commutative ring. The categories
$\roman{GL}A$ and $\roman{Pr}A$ described above become bipermutative
categories using the tensor
product $\otimes_{A}$, when we identify
$A^{m}\otimes_{A} A^{n}$ with $A^{mn}$ using lexicographical order on
the standard
basis. \endoldexample

We prove the following result in \refname{secassoc}.

\proclaim{Theorem}\label{num12} Bipermutative structure on a
small permutative category $\R$ determines and is determined by a
multifunctor $E\Sigma_*\to\P$ sending the single object of
$E\Sigma_*$ to $\R$.
\endproclaim

\proclaim{Corollary}\label{num13} Any small bipermutative category is an
associative category.
\endproclaim

\demo{Proof} Compose the given multifunctor $E\Sigma_{*}\to \P$ with
the map of operads $\Sigma_*\to E\Sigma_*$ that is the inclusion of
objects.
\enddemo

The ``small'' hypothesis is not really necessary: A component of the
argument for \refname{num12} is a direct verification that a
bipermutative category satisfies diagram~(e) of \refname{num9}
(see Figure~1 on page~\refpage{figure}).

Since the map $E\Sigma_*\to*$ of operads is a weak equivalence,
and the algebras for the one-point operad in any symmetric
monoidal category are the commutative monoids in that category,
Theorems \refnum{num6} and \refnum{num5.5} now give the following
corollary.

\proclaim{Corollary} If $\R$ is a bipermutative category, then
$K\R$ is equivalent to a strictly commutative ring symmetric
spectrum.
\endproclaim

\section{The $K$-Theory of Permutative Categories}\label{secK}

In this section, we construct the underlying functor of our $K$-theory
multifunctor from permutative categories to symmetric spectra, and
show that it is equivalent to the $K$-theory functor in the literature.
Since our functor is a modification of the usual Segal
construction of the $K$-theory spectrum of a small permutative
category, we describe that first, using the construction from
\cite{15}.

\definition{Construction} For a small permutative category $\C$ and a finite
based set $A$, let $\barC_{A}$ denote the category whose objects
are the systems $\{C_{S}, \rho_{S,T}\}$, where
\roster
\item $S$ runs through
the subsets of $A$ {\bf not} containing the basepoint,
\item $S,T$ runs
through the pairs of such subsets with $S\cap T=\emptyset$,
\item the
$C_{S}$ are objects of $\C$ and the $\rho_{S,T}$ are isomorphisms
$C_{S}\oplus C_{T}\to C_{S\cup T}$
\endroster
such that $\C_{S}=0$ and
$\rho_{S,T}=\id_{C_{T}}$ when $S=\emptyset$, and the following
diagrams commute for all $S,T,U$:
$$ \xymatrix{
C_{S}\oplus C_{T}\ar[d]_{\gamma}\ar[r]^{\rho_{S,T}}
&C_{S\cup T}\ar @{=}[d]\\
C_{T}\oplus C_{S}\ar[r]_{\rho_{T,S}} &C_{T\cup S}
}\quad
\xymatrix@C+15pt{
C_{S}\oplus C_{T}\oplus C_{U}
\ar[d]_{\id_{C_{S}}\oplus \rho_{T,U}}
\ar[r]^{\rho_{S,T}\oplus \id_{U}}
&C_{S\cup T}\oplus C_{U}
\ar[d]^{\rho_{S\cup T,U}}\\
C_{S}\oplus C_{T\cup U}
\ar[r]_{\rho_{S,T\cup U}}
&C_{S\cup T\cup U}.
}
$$
A morphism $f\colon \{C_{S},\rho_{S,T}\}\to
\{C'_{S},\rho'_{S,T}\}$ consists of morphisms $f_{S}\colon
C_{S}\to C'_{S}$ in $\C$ for all $S$, such that
$f_{\emptyset}=\id_{0}$, and the following diagram commutes for all
$S,T$:
$$ \xymatrix{
C_{S}\oplus C_{T}\ar[r]^{\rho_{S,T}}\ar[d]_{f_{S}\oplus f_{T}}
&C_{S\cup T}\ar[d]^{f_{S\cup T}}\\
C'_{S}\oplus C'_{T}\ar[r]_{\rho'_{S,T}}
&C'_{S\cup T}.
} $$
\enddefinition

\oldremark{Remark}
The construction is described in \cite{15} in terms of based subsets
of a based set as indices.  This leads to
some awkwardness in defining functoriality which the formalism
above avoids.  The description in \cite{15} can be recovered simply by
reattaching the basepoint to all indexing subsets.
\endoldremark

\proclaim{Theorem}\label{barcgamma}
The assignment $A\mapsto\barC_A$ defines a functor $\barC$ from
the category of finite based sets to the category of small
categories.
\endproclaim

\demo{Proof}
A map of finite based sets $\alpha \colon A\to A'$ induces the
functor $\barC_\alpha$ that sends the object
$\{C_{S},\rho_{S,T}\}$ of $\barC_{A}$ to the object
$\{C^{\alpha}_{S},\rho^{\alpha}_{S,T} \}$ of $\barC_{A'}$ where
$C^{\alpha}_{S}=C_{\alpha^{-1}S}$ and
$\rho^{\alpha}_{S,T}=\rho_{\alpha^{-1}S,\alpha^{-1}T}$.  Note that
since $\alpha$ is basepoint-preserving, $\alpha^{-1}(S)$ does
{not} contain the basepoint.  Likewise, $\barC_\alpha$ sends the
map $\{f_{S}\}$ to the map $\{f^{\alpha}_{S}\}$ where
$f^{\alpha}_{S}=f_{\alpha^{-1}S}$. Clearly, when $\alpha$ is the
identity, $\barC_\alpha$ is the identity, and for $\alpha'\colon
A'\to A''$,
$\barC_{\alpha'\circ\alpha}=\barC_{\alpha'}\circ\barC_\alpha$.
\enddemo

In the conventions of \cite{1}, a ``$\Gamma$-space'' is a functor
from the category of finite based sets to the category of
simplicial sets that takes the trivial based set (consisting of
only the base point) to a constant one point simplicial set. It
follows that $N\barC$ is a $\Gamma$-space, where $N$ denotes the
nerve functor. Standard notation is to use $\b n$ to denote the
finite based set $\{0,1,2,\ldots,n\}$ with $0$ serving as the
basepoint.   The category $\barC_{\b1}$ is then canonically
isomorphic to the original category $\C$.  For $n>0$, the based
maps $\b n\to \b1$ that send all but one of the non-basepoint
elements to the basepoint induce a functor
$$ p_{n}\colon \barC_{\b n}\to \barC_{\b1}\times \cdots \times \barC_{\b1}
\iso  \C\times \cdots \times\C  $$
that is easily identified as the functor that sends
$\{C_{S},\rho_{S,T}\}$ to $(C_{\{1\}},\ldots,C_{\{n\}})$ and is an
equivalence of categories.  The $\Gamma$-space $N\barC$ is
therefore ``special'' in the terminology of \cite{1} in that the
map $p_{n}\colon N\barC_{\b n}\to N\barC_{\b1}\times \cdots
\times N\barC_{\b1}$ is a homotopy equivalence for each $n>0$.

The spectrum associated to a $\Gamma$-space $X$ is constructed as
follows.  Let $\SSp$ denote the following simplicial model of the
circle: The set of $n$-simplices is $\SSpd{1}{n}=\b{n}$ with face
maps $d_{i}$ the order-preserving maps that delete the element $i$
and the degeneracy maps $s_{i}$ the order-preserving maps that
skip the element $i$.  Then $\SSp$ has one $0$-simplex and one
non-degenerate $1$-simplex; all $n$-simplices are degenerate for
$n>1$. Regarding $\SSp$ as a simplicial based set and applying the
functor $X$ degreewise, we obtain a bisimplicial set $X_{\SSp}$,
which we regard as a simplicial set by taking the diagonal.
Writing $\SSpn$ for the $n$-fold smash power of $\SSp$ (with
$\SSps0$ the constant simplicial set $\b1$), we likewise get
simplicial sets $X_{\SSpn}$.  Since $S_q^{n-1}\wedge S_q^1
=S_q^n$, each $q$-simplex $x$ of $\SSp$ induces a map of based
sets
$$ \SSpd{n-1}{q}\iso \SSpd{n-1}{q}\sma \{0,x\} \to \SSpd{n}{q} $$
that assemble to a based map
$$ X_{\SSpd{n-1}{q}}\sma \SSpd1q
\iso
\bigvee_{x\in\SSpd1q\setminus\{0\}}(X_{\SSpd{n-1}{q}\sma \{0,x\}})
\to X_{\SSpd{n}{q}} $$ for each $q$.
Taking these together for all $q$ and $n$ form the ``structure
maps'' $\Sigma X_{\SSps{n-1}}\to X_{\SSpn}$ that make
$\{X_{\SSpn}\}$ into a spectrum.  In fact, $\{X_{\SSpn}\}$ forms a
symmetric spectrum, where the symmetric group action on
$X_{\SSpn}$ is induced by permuting the smash factors of $\SSpn$.
The main theorem of
\cite{17} then can be phrased as saying that when $X$ is a
special $\Gamma$-space, this spectrum is an ``almost
$\Omega$-spectrum'' in that after geometric realization, the maps
$$ |{X_{\SSpn}}|\to \Omega|{X_{\SSps{n+1}}}|
$$
adjoint to the structure maps are homotopy equivalences for all
$n\geq1$.

Although we have followed \cite{15} in constructing $\barC$ and
\cite{1} in constructing the associated (symmetric) spectrum, we
refer to this as Segal's construction.% of the $K$-theory of $\C$.

\definition{Definition}
Segal's construction of $K$-theory of the permutative category $\C$ is
the symmetric spectrum $\KSeg\C=\{N\barC_{\SSpn}\}$.
\enddefinition

Previously, the main difficulty with constructing ring and module
structures on the spectra associated to permutative categories was
the lack of a symmetric monoidal product on the target category of
spectra. Even using the category of symmetric spectra, which does
have a symmetric monoidal product, the previous definition does
not carry ring structures (e.g., associative category structures)
to ring structures.  A suitable collection of maps
$$ N\barC_{\b m}\sma N\barC_{\b n}\to N\barC_{\b{m}\sma \b{n}} $$
would give rise to a pairing $\KSeg\C \sma \KSeg\C \to \KSeg\C$,
but no reasonable definition of pairing on the permutative
category $\C$ gives rise to such a collection of maps.  We can
illustrate this by looking at just the zero simplices, or
equivalently, the objects in the categories. Given some kind of
pairing $\otimes$ on $\C$ and objects $\{C_{S},\rho_{S,T}\}$ of
$\barC_{\b m}$ and $\{C'_{S},\rho_{S,T}\}$ of $\barC_{\b n}$, we
need to construct an object $\{C''_{S},\rho_{S,T}\}$ of
$\barC_{\b{m}\sma \b{n}}$.  It seems natural to take
$$ C''_{S\times T}=C_{S}\otimes C'_{T} $$
on the subsets of the form $S\times T\subset \b{m}\sma \b{n}$, but how
do we fill in the objects $C''_{U}$ for subsets $U$ not of this form?

Our basic idea is to modify the construction of $\barC$ so the
objects correspond only to those subsets of the appropriate form.
The set of $q$-simplices $\SSpd{n}{q}$ of $\SSpn$ is
$\SSpd1q\sma\cdots\sma\SSpd1q$; instead of using
$N\barC_{\SSpd{n}{q}}$ where we choose objects $C_{T}$ for all
subsets $T$ of $\SSpd1q\sma\cdots\sma\SSpd1q$ not containing the
basepoint, we can use a variant where we only choose them for the
subsets of the form $T_{1}\times\cdots\times T_{n}$.  We make one
other alteration: Since we have defined the multicategory of
permutative categories using lax distributivity maps, we do not
require the morphisms $\rho$ to be isomorphisms. Before describing
the construction, it is useful to introduce the following
notation. Given finite basepoint-free (sub)sets
$S_{1},\ldots,S_{n}$, we write $\br{S}$ for the $n$-tuple
$(S_{1},\ldots,S_{n})$.  Given a finite basepoint-free set $T$ and
$i\in\{1,\ldots,n\}$, we write $\br{S\subst{i}{T}}$ for the
$n$-tuple $(S_{1},\ldots,S_{i-1},T,S_{i+1},\ldots,S_{n})$ obtained
by substituting $T$ in the $i$-th position.

\definition{Construction}
For a small permutative category $\C$ and finite based sets
$A_{1},\ldots, A_{n}$, let $\barC_{(A_{1},\ldots,A_{n})}$ denote
the category whose objects are the systems $\{C_{\br{S}},
\rho_{\br{S};i,T,U}\}$, where
\roster
\item $\br{S}=(S_{1},\ldots,S_{n})$ runs through all $n$-tuples of
based subsets $S_{i}\subset A_{i}$,
\item For $\rho_{\br{S};i,T,U}$, $i$ runs through $1,\ldots,n$, and
$T,U$ run through the basepoint-free subsets of $S_{i}$ with $T\cap
U=\emptyset$ and $T\cup U = S_{i}$
\item The $C_{\br{S}}$ are objects of $\C$, and
\item The $\rho_{\br{S};i,T,U}$ are morphisms $C_{\br{S\subst{i}T}}\oplus
C_{\br{S\subst{i}U}}\to
C_{\br{S}}$ in $\C$
\endroster
such that
\roster
\item $\C_{\br{S}}=0$ if $S_{k}=\emptyset$ for any $k$,
\item $\rho_{\br{S};i,T,U}=\id$ if any of the $S_{k}$ (for any $k$),
$T$, or $U$ are empty.
\item For all $\rho_{\br{S};i,T,U}$ the following
diagram commutes:
$$ \xymatrix@C+40pt @R+12pt{
C_{\br{S\subst{i}T}}\oplus C_{\br{S\subst{i}U}}
\ar[d]_{\gamma}\ar[r]^{\rho_{\br{S};i,T,U}}
&C_{\br{S}}\ar @{=}[d]\\
C_{\br{S\subst{i}U}}\oplus C_{\br{S\subst{i}T}}
\ar[r]_{\rho_{\br{S};i,U,T}}
&C_{\br{S}}
} $$
\item For all $\br{S}$, $i$, and $T,U,V\subset A_{i}$ with $T\cup
U\cup V=S_{i}$ and $T$, $U$, and $V$ all mutually disjoint, the
following diagram commutes:
$$ \xymatrix@C+50pt @R+12pt{
C_{\br{S\subst{i}T}}\oplus C_{\br{S\subst{i}U}}\oplus C_{\br{S\subst{i}V}}
\ar[d]_{\id\oplus \rho_{\br{S\subst{i}(U\cup V)};i,U,V}}
\ar[r]^{\rho_{\br{S\subst{i}(T\cup U)};i,T,U}\oplus \id}
&C_{\br{S\subst{i}(T\cup U)}}\oplus C_{\br{S\subst{i}V}}
\ar[d]^{\rho_{\br{S};i,T\cup U,V}}\\
C_{\br{S\subst{i}T}}\oplus C_{\br{S\subst{i}(U\cup V)}}
\ar[r]_{\rho_{\br{S};i,T,U\cup V}}
&C_{\br{S}},
} $$
\item For all $\rho_{\br{S};i,T,U}$ and $\rho_{\br{S};j,V,W}$ with
$i\ne j$, the following diagram commutes:
$$ \xymatrix@C-35pt @R-3pt{
&C_{\br{S\subst{j}V}}\oplus
C_{\br{S\subst{j}W}}
\ar[ddr]^{\rho_{\br{S};j,V,W}}
\\
C_{\br{S\subst{i}T\subst{j}V}}\oplus
C_{\br{S\subst{i}U\subst{j}V}}\oplus
C_{\br{S\subst{i}T\subst{j}W}}\oplus
C_{\br{S\subst{i}U\subst{j}W}}
\ar[ur]^{(\rho_{\br{S\subst{j}V};i,T,U})\oplus (\rho_{\br{S\subst{j}W};i,T,U})\qquad\qquad }
\ar[dd]_{\id\oplus \gamma \oplus \id}
\\
&&C_{\br{S}}.\\
C_{\br{S\subst{i}T\subst{j}V}}\oplus
C_{\br{S\subst{i}T\subst{j}W}}\oplus
C_{\br{S\subst{i}U\subst{j}V}}\oplus
C_{\br{S\subst{i}U\subst{j}W}}
\ar[dr]_{(\rho_{\br{S\subst{i}T};j,V,W})\oplus (\rho_{\br{S\subst{i}U};j,V,W})\qquad\qquad }
\\
&C_{\br{S\subst{i}T}}\oplus
C_{\br{S\subst{i}U}}
\ar[uur]_{\rho_{\br{S};i,T,U}}
} $$
\endroster
A morphism $f\colon \{C_{\br{S}},\rho_{\br{S};i,T,U}\}\to
\{C'_{\br{S}},\rho'_{\br{S};i,T,U}\}$
consists of morphisms $f_{S}\colon C_{\br S}\to C'_{\br S}$ in
$\C$ for all $\br S$ such that $f_{\br S}$ is the identity
$\id_{0}$ when $S_{i}=\emptyset$ for any $i$, and the following
diagram commutes for all $\rho_{\br{S};i,T,U}$:
$$ \xymatrix@C+15pt{
C_{\br{S\subst{i}T}}\oplus C_{\br{S\subst{i}U}}\ar[r]^-{\rho_{\br{S};i,T,U}}
\ar[d]_{f_{\br{S\subst{i}T}}\oplus f_{\br{S\subst{i}U}}}
&C_{\br{S}}\ar[d]^{f_{\br{{S}}}}\\
C'_{\br{S\subst{i}T}}\oplus
C'_{\br{S\subst{i}U}}\ar[r]_-{\rho'_{\br{S};i,T,U}}
&C'_{\br{S}}.
} $$
\enddefinition

We make $\barC$ into a functor from $n$-tuples of based spaces to
categories just as in \refname{barcgamma}. The categories
$\barC_{\br{A}}$ have further functoriality as well:

\olddefinition{Permutation Functors}
A permutation $\sigma$ in $\Sigma_{n}$ induces a functor
$$ \sigma_{!}\colon \barC_{(A_{1},\ldots,A_{n})}
\to\barC_{(A_{\sigma^{-1}(1)},\ldots,A_{\sigma^{-1}(n)})}, $$
which is an isomorphism of categories, as follows:  The object
$\{C_{\br{S}},\rho_{\br{S};i,T,U}\}$ is sent to the object
$\{C^{\sigma}_{\br{S'}},\rho^{\sigma}_{\br{S'};i,T}\}$ where
$$
C^{\sigma}_{\br{S'}}
 = C_{\sigma\br{S'}},
\qquad
\rho^{\sigma}_{\br{S'};i,T,U}
 = \rho_{\sigma\br{S'};\sigma(i),T,U},
\qquad
\sigma\br{S'}=(S'_{\sigma(1)},\ldots,S'_{\sigma(n)}),
$$
so if $S'_{i}=S_{\sigma^{-1}(i)}\subset A_{\sigma^{-1}(i)}$, then
$\sigma\br{S'}=\br{S}$. The morphism $\{f_{\br{S}}\}$ is sent to
the morphism $\{f^{\sigma}_{\br{S'}}\}$ where
$f^{\sigma}_{\br{S'}}=f_{\sigma
\br{S'}}$.
\endolddefinition

\olddefinition{Extension Functors}
We have an isomorphism of
categories
$$ e\colon \barC_{(A_{1},\ldots,A_{n})}\to
\barC_{(A_{1},\ldots,A_{n},\b1)} $$
defined as follows:
The object $\{C_{\br{S}},\rho_{\br{S};i,T,U}\}$ is sent to
the object $\{C^{e}_{\br{S'}},\rho^{e}_{\br{S'};i,T,U}\}$ where
$$ \alignat 3
C^{e}_{(S_{1},\ldots,S_{n},\{1\})}
 &= C_{\br{S}},
&\quad
\rho^{e}_{(S_{1},\ldots,S_{n},\{1\});i,T,U}
 &= \rho_{\br{S};i,T,U} \hbox to 0pt{{ for $i<n+1$,}\hss}
&\quad &\\
C^{e}_{(S_{1},\ldots,S_{n},\emptyset)}&=0,
&\quad
\rho^{e}_{(S_{1},\ldots,S_{n},\emptyset);i,T,U}
 &= \id,&\quad
\rho^{e}_{(S_{1},\ldots,S_{n},\{1\});n+1,T,U}&=\id.
\endalignat $$
The morphism $\{f_{\br{S}}\}$ is sent to the morphism
$\{f^{e}_{\br{S'}}\}$ where
$$f^{e}_{(S_{1},\ldots,S_{n},\{1\})}=f_{\br{S}}, \qquad
f^{e}_{(S_{1},\ldots,S_{n},\emptyset)}=\id.
$$
This description of the components of the objects and morphisms is
complete since the only two basepoint-free subsets of $\b1$ are
$\{1\}$ and $\emptyset$. The inverse of this isomorphism is
induced by dropping the $\{1\}$ from $(n+1)$-tuples of the form
$(S_{1},\ldots,S_{n},\{1\})$.  Of course, for any other set
$\{*,x\}$ with precisely one non-basepoint, we have an extension
functor $e_{x}\colon
\barC_{(A_{1},\ldots,A_{n})}\to
\barC_{(A_{1},\ldots,A_{n},\{*,x\})}$ given by the composite of $e$
and the functor induced by the unique based bijection $\b1\to\{*,x\}$.
\endolddefinition

The various functors above satisfy certain compatibility relations
that we describe implicitly in the next section, by abstracting
them into the definition of symmetric functor.  We can also extend
such functors naturally to functors from finite simplicial based
sets to simplicial categories by applying the functor degreewise.
The nerve of a simplicial category is formed by taking the nerve
degreewise and then taking the diagonal. The underlying functor of
the $K$-theory multifunctor we describe in
Sections~\refnum{pctosf} and~\refnum{sftoss} is naturally
isomorphic to the $K$-theory functor in the following definition.

\definition{Definition}
For a small permutative category $\C$, the symmetric spectrum
$\Knew\C$ is defined by $(\Knew\C)(0)=N\barC_{S^{0}}$,
$(\Knew\C)(1)=N\barC_{\SSp}$, $(\Knew\C)(2)=N\barC_{(\SSp,\SSp)}$,
and in general,
$$ (\Knew\C)(n) =
N\barC_{(\underbrace{\scriptstyle \SSp,\dotsc,\SSp}_{n})},
$$
with symmetric
action induced by the permutation functors above and structure maps
$$ N\barC_{(\SSpd1q,\ldots,\SSpd1q)}\sma \SSpd1q\iso
\bigvee_{x\in\SSpd1q\setminus\{0\}}
N\barC_{(\SSpd1q,\ldots,\SSpd1q,\{0,x\})}
\to N\barC_{(\SSpd1q,\ldots,\SSpd1q,\SSpd1q)}  $$
induced by the extension functors above.
\enddefinition

We close this section by showing that the symmetric spectra
$\KSeg\C$ and $\Knew\C$ are equivalent.   First we note that we
have a canonical functor $\barC_{A_{1}\sma \cdots \sma A_{n}}\to
\barC_{(A_{1},\ldots,A_{n})}$ that takes the object
$\{C_{S},\rho_{S,T}\}$ to the object
$$ C_{\br{S}}=C_{S_{1}\times\cdots\times S_{n}}, \qquad
\rho_{\br{S};i,T,U} =
 \rho_{S_{1}\times\cdots\times T\times\cdots \times S_{n},S_{1}\times\cdots
 \times U\times\cdots \times S_{n}}.
$$
This functor is natural in $\C$ and $A_{1},\ldots,A_{n}$ and commutes
with the permutation and extension
functors.  We therefore get a natural map of symmetric spectra
$\KSeg\to\Knew$.

\proclaim{Theorem}
The natural map of symmetric spectra $\KSeg\C\to\Knew\C$ is a level
equivalence for every $\C$.
\endproclaim

\demo{Proof}
It suffices to show that the map
$N\barC_{\b{m}_{1}\sma\cdots\sma\b{m}_{n}}\to N\barC_{\br{\b{m}}}$
is a weak equivalence for all $n$, $\br{\b m}$. Write $\sma\br{\b
m}$ as an abbreviation for $\b{m}_{1}\sma\cdots\sma \b{m}_{n}$ and
let $m=m_{1}\dotsb m_{n}$. Let $p_{m}\colon
\barC_{\sma\br{\b{m}}}\to\C^{m}$ denote the functor that takes
$\{C_{S},\rho_{S,T}\}$ to the $m$-tuple whose
$(i_1,\ldots,i_n)$-th coordinate is $C_{\{(i_1,\ldots,i_n)\}}$.
Then $p_m$ is an equivalence of categories. Let $q_{\br m}\colon
\barC_{\br{\b{m}}}\to \C^{m}$ denote the functor that takes
$\{C_{\br{S}},\rho_{\br{S};i,T,U}\}$ to the $m_1\cdots m_n$-tuple
whose $(i_1,\ldots,i_n)$'th coordinate is
$C_{(\{i_{1}\},\ldots,\{i_{n}\})}$.  Then $q_{\br m}$ is not an
equivalence of categories but does have a left adjoint, namely the
functor that sends an object with coordinates
$(X_{i_{1},\ldots,i_{n}})$ to the object with
$$ C_{\br{S}}=\bigoplus_{i_{1}\in S_{1}}\cdots
\bigoplus_{i_{n}\in S_{n}} X_{i_{1},\ldots,i_{n}}
$$
(ordered using the natural order on $S_{i}\subset \b{m}$), with the
convention that the empty sum is the unit $0$ of $\C$; the
$\rho$'s are defined by the appropriate rearrangement using the
commutativity isomorphism $\gamma$. The functor $q_{\br m}$
therefore induces a homotopy equivalence on nerves.  Since the
functor $p_{m}$ factors as the composite of the functor
$\barC_{\sma\br{\b{m}}}\to\barC_{\br{\b{m}}}$ we are interested in
and the functor $q_{\br m}$, we conclude that the map
$N\barC_{\sma\br{\b m}}\to N\barC_{\br{\b{m}}}$ is a homotopy
equivalence.  This completes the proof.
\enddemo

\section{The Multicategory of Symmetric Functors}\label{secsym}

Extending the $K$-theory functor to a multifunctor from the
multicategory of permutative categories to the multicategory of
symmetric spectra requires a detailed study of the properties of
the constructions of the previous section.  Instead of carrying
along the details, it is useful to abstract the essential
properties, and this leads us to the symmetric functors that we define in
this section. To simplify things, the definition of symmetric
functor throws away some inessential data that is recoverable up
to isomorphism. Instead of working with all finite based sets,
symmetric functors are defined just in terms of the finite based
sets $\b0,\b1,\b2,\ldots$.  Instead of keeping track of some
analogue of the extension functors of the previous section,
symmetric functors compress these away by looking at the colimit.
Finally, to avoid possible confusion as to the meaning of
``functor'' and ``natural transformation'' where they occur below,
we define symmetric functors with values in an arbitrary category
$\ACat$ that has finite (categorical) products. In our
applications $\ACat$ is always either $\Cat$, the category of
small categories, or $\SS$, the category of simplicial sets.

\definition{Definition}\label{num48} Let $\F$ be the category with
objects $\b0, \b1, \b2,\ldots$, where (as above) $\b n=\{0,1,2,\ldots,n\}$, and
with morphisms based functions with $0$ as the basepoint. Form the
direct system of categories
$$\F\to\F^2\to\F^3\to\cdots,
$$
where $\F^{p-1}$ maps to $\F^p$ by setting the last coordinate
equal to $\b1\in\obj{\F}$ for objects and $\id_{\b1}$ for
morphisms. Let $\F^\infty$ be the colimit of this system; it is a
category with objects sequences $(\b n_1,\b n_2,\ldots)$ with $\b
n_i=\b1$ for all but finitely many $i$, and morphisms sequences of
morphisms $(f_1,f_2,\ldots)$, where $f_i\colon \b n_i\to\b n'_i$ with
$f_i=\id_{\b1}$ for all but finitely many $i$.  A {\bf symmetric
functor} (with values in $\ACat$) is a functor
$$F\colon \F^\infty\to\ACat
$$
such that $F(\b n_1,\ldots)=*$ (the final object) whenever any of the
$\b n_i=\b0$,
together with an action of $\Aut(\N)$ in the following sense. Let
$\theta\in\Aut(\N)$. Then there is an induced functor
$\theta^*\colon \F^\infty\to\F^\infty$, given explicitly on
objects by $\theta^*(\b n_1,\b n_2,\ldots):=(\b n_{\theta(1)},\b
n_{\theta(2)},\ldots)$ and similarly on morphisms.  We require, as
part of the structure, a natural isomorphism $\theta_!\colon
F\theta^*\to F$ for each $\theta$. We express this
diagrammatically as
$$\xymatrix@=40pt{
\F^\infty\ar[r]^-{\theta^*}
\drtwocell<0>_{F}{<-2.5>\;\theta_!}
&\F^\infty
\ar[d]^-{F}
\\
&\ACat.
\\}$$
These must satisfy the coherence conditions that $\id_!=\id_F$,
and that the following diagrams of natural transformations
coincide:
%% The two cell arrows need to be fixed here
$$ \vcenter{
\xymatrix@=40pt{
\F^\infty\ar[r]^-{(\theta_1)^*}
\drtwocell<0>_{F}{<-2>\phantom{mm}(\theta_1)_!}
&\F^\infty
\ar[r]^-{(\theta_2)^*}
\ar[d]_(.65){F}
&\F^\infty
\dltwocell<0>^{F}{<2.5>(\theta_2)_!\quad}
\\
&\ACat
\\}
}\qquad=\qquad\vcenter{
\xymatrix@R=40pt{
\F^\infty
%\ar@{{}{}{}}[d]_(.6){=}
\ar[rr]^-{(\theta_1\theta_2)^*}
\drtwocell<0>_{F}{<-3>\qquad\,(\theta_1\theta_2)_!}
&&\F^\infty
\ar[dl]^-{F}
\\
\relax&\ACat.
\\}
} $$ In addition, we require a component of the natural
isomorphism $\theta_!$ to coincide with the identity whenever its
indexing object $\br{\b n}\in\obj{\F^\infty}$ satisfies
$\theta(i)\ne i\Rightarrow {\b n_i}={\b1}$, i.e., whenever the
only entries in $\br{\b n}$ that $\theta$ moves are $\b1$'s. This
concludes the definition of the objects of the multicategory of
symmetric functors.
\enddefinition

To give meaning to the previous definition, the reader should keep
in mind the following example, which is also the one of main
interest.

\example{Example}\label{Jdef}
Let $\C$ be a small permutative category, and let
$$ J\C\br{\b n} = \Colim_{m\geq m_{0}} \barC_{(\b n_{1},\ldots,\b n_{m})}
$$ where ${\b n}_{m_{0}}$ is the last non-$\b1$ entry, and the colimit is taken
over the extension functors $e$ (which are isomorphisms of
categories).  Then $J\C$ becomes a symmetric functor with values
in $\Cat$, where the isomorphisms $\theta_{!}$ are induced by the
permutation functors. Likewise $NJ\C$ is a symmetric functor with
values in $\SS$.
\endexample

So far we have described just the objects, and we still need to
describe the $k$-morphisms.  In the following definition, we use
$k\cdot\N$ to denote $\{1,\ldots,k\}\times\N$. For a bijection
$\beta
\colon \N\to k\cdot \N$, $\beta^{*}$ denotes the functor
$(\F^{\infty})^{k}\to
\F^{\infty}$ defined by the formula
$$ \beta^{*}\colon
 ((\b n_{11},\b n_{12}, \b n_{13},\ldots),\ldots,
  (\b n_{k1},\b n_{k2}, \b n_{k3},\ldots))
\mapsto (\b n_{\beta(1)}, \b n_{\beta(2)}, \b n_{\beta(3)},\ldots).
$$
Abstractly, $\beta^{*}$ is the induced isomorphism of categories
given by identifying $(\F^\infty)^k$ with a subcategory of
$\F^{k\cdot\N}$, pulling back along $\beta$, and observing that
this process factors through the subcategory $\F^{\infty}$ of
$\F^{\N}$.

\definition{Definition}\label{sfkmap} The $k$-morphisms of the
multicategory of symmetric functors are defined as follows:  Let
$F_1,\ldots,F_k$, and $G$ be symmetric functors.  A $k$-morphism
$f\colon (F_1,\ldots,F_k)\to G$ assigns to each choice of
bijection $\beta\colon \N\overto{\cong}k\cdot\N$ a natural
transformation $f_\beta$ as in the following (noncommutative)
diagram:
$$\xymatrix@=40pt{
(\F^\infty)^k
\ar[r]^-{F_1\times\cdots\times F_k}
\ar[d]_-{\beta^*}
\drtwocell<\omit>{\;f_\beta}
&\ACat^k
\ar[d]^-{\times}
\\
\F^\infty
\ar[r]_-{G}
&\ACat.
\\}$$
Here the right vertical arrow is the categorical product in
$\ACat$, and the left vertical arrow is as described above.  These
transformations $f_\beta$ must satisfy the following coherence
conditions, the first of which is an equivariance statement
connecting the actions of $\Aut(\N)$ on the $F$'s with the action
on $G$, and the second relates the $f_\beta$'s for different
choices of $\beta$.  For the first one, given elements
$\theta_1,\ldots,\theta_k$ of $\Aut(\N)$, let
$\ttheta=\beta^{-1}(\theta_1\coprod\cdots\coprod\theta_k)\beta$.
Abbreviate $\theta_!$ for
$(\theta_1)_!\times\cdots\times(\theta_k)_!$. Then we require the
following diagrams of natural transformations to coincide:
$$ \vcenter{
\xymatrix@R=30pt @C=55pt{
(\F^\infty)^k
\ar[r]^-{\theta_1^*\times\cdots\times\theta_k^*}
\drtwocell<0>_{F_1\times\cdots\times
F_k\phantom{mmmm}}{<-2.5>\;\theta_!}
\ar[dd]_-{\beta^*}
&(\F^\infty)^k
\ar[d]^-{F_1\times\cdots\times F_k}
\\
\relax
\drtwocell<\omit>{\;f_\beta}
&\ACat^k
\ar[d]^-{\times}
\\
\F^\infty
\ar[r]_-{G}
&\ACat
\\}
}=\qquad\vcenter{
\xymatrix@R=30pt @C=15pt{
(\F^\infty)^k
\ar[rr]^-{\theta_1^*\times\cdots\times\theta_k^*}
\ar[dd]_{\beta^*}
&&(\F^\infty)^k
\ar[d]^-{F_1\times\cdots\times F_k}
\ar[dl]_-{\beta^*}
\ddtwocell<\omit>{<4>f_\beta}
\\
&\F^\infty
\ar[dr]^-{G}
&\ACat^k
\ar[d]^{\times}
%\dtwocell<\omit>{<2> f_\beta}
\\
\F^\infty
\ar[ur]^{(\ttheta)^*}
\rrtwocell<0>_{G}{<-2.5>\quad(\ttheta)_!}
&&\ACat.
\\}
%%\xymatrix@C=40pt{
%%(\F^\infty)^k
%%\ar[r]^-{\theta_1^*\times\cdots\times\theta_k^*}
%%\ar[dd]_{\beta^*}%\ar[dd]_-{=\phantom{mmmmmmm}\beta^*}
%%&(\F^\infty)^k
%%\ar[r]^-{F_1\times\cdots\times F_k}
%%\ddrtwocell<\omit>{<-1>\;f_\beta}
%%\ar[d]_-{\beta^*}
%%&\ACat^k
%%\ar[dd]^{\times}%\ar[dd]^-{\times\phantom{mmmmmmmmmm}}
%%\\
%%&\F^\infty
%%\ar[dr]_-{G}
%%\\
%%\F^\infty
%%\ar[ur]^{(\ttheta)^*}
%%\rrtwocell<0>_{G}{<-2.5>\quad(\ttheta)_!}
%%&&\ACat.
%%\\}
} \leqno{\text{\bf \refnum{sfkmap}(a)}}
$$
For the second condition, given bijections $\beta_1$ and
$\beta_2\colon \N\overto{\cong}k\cdot\N$, there is a unique
$\theta\in\Aut(\N)$ such that $\beta_1=\beta_2\theta$. We require
the following diagrams of natural transformations to coincide:
$$ \vcenter{
\xymatrix{
(\F^\infty)^k
\ar[dd]_-{\beta_2^*}
\ar[dr]_-{\beta_1^*}
\ar[rr]^-{F_1\times\cdots\times F_k}
\ddrrtwocell<\omit>{<-4.5>\quad f_{\beta_1}}
&&\ACat^k
\ar[dd]^-{\times}
\\
&\F^\infty
\ar[dr]^-{G}
\\
\F^\infty
\ar[ur]^-{\theta^*}
\rrtwocell<0>_{G}{<-2.5>\quad\theta_!^G}
&&\ACat
\\}
}\qquad=\qquad\vcenter{
\xymatrix@=60pt{
(\F^\infty)^k
\ar[r]^-{F_1\times\cdots\times F_k}
\ar[d]_-{\beta_2^*}
%\ar[d]_-{=\phantom{mmmmmmm}\beta_1^*}
\drtwocell<\omit>{\quad f_{\beta_2}}
&\ACat^k
\ar[d]^-{\times}
%\ar[d]^-{\times\phantom{mmmmmmm}}
\\
\F^\infty
\ar[r]_-{G}
&\ACat.
\\}} \leqno{\text{\bf \refnum{sfkmap}(b)}}
$$
This condition implies that it is sufficient to specify $f_\beta$
for a single choice of $\beta$, since all other bijections from
$\N$ to $k\cdot\N$ are of the form $\beta\theta$ for
$\theta\in\Aut(\N)$.  In particular, if $k=1$, it suffices to use
$\beta=\id_\N$, in which case we find that a 1-morphism is just an
equivariant natural transformation.
\enddefinition

When $\ACat$ is enriched over small categories or simplicial sets, the
$k$-morphisms
from $F_{1},\ldots,F_{k}$ to $G$ inherit a canonical enrichment:
The natural transformations between functors into $\ACat$
form a category or simplicial set, and each requirement for
diagrams to coincide specifies an equalizer, defining the
$k$-morphisms as a category or simplicial set.  Moreover, since the
nerve functor preserves all limits, when $\ACat$ is enriched over
small categories, the nerve of the category of $k$-morphisms is
canonically isomorphic to the simplicial set of $k$-morphisms obtained
by viewing $\ACat$ as enriched over simplicial sets.

Finally, to give symmetric functors the structure of a
multicategory, we must specify the $\Sigma_k$-action on the
$k$-morphisms and the multiproduct.  We write
$\kmap(F_1,\ldots,F_k;G)$ for the $k$-morphisms from
$(F_1,\ldots,F_k)$ to $G$.

\definition{Definition}\label{sfmultcat}
Given
$\sigma\in\Sigma_k$, we define
$$\sigma^*\colon \kmap(F_1,\ldots,F_k;G)
\to\kmap(F_{\sigma(1)},\ldots,F_{\sigma(k)};G)
$$
by requiring $(\sigma^*f)_\beta$ to be the natural map
$$\xymatrix@C=50pt{
(\F^\infty)^k
\ar[r]^-{F_{\sigma(1)}\times\cdots\times F_{\sigma(k)}}
\ar[d]_-{\sigma^*}
&\ACat^k
\ddrtwocell<0>^{\times}{<2>\sigma_*}
\ar[d]_-{\sigma^*}
\\
(\F^\infty)^k
\ar[r]^-{F_1\times\cdots\times F_k}
\ar[d]_-{\beta^*}
\drrtwocell<\omit>{f_\beta}
&\ACat^k
\ar[dr]_-{\times}
\\
\F^\infty
\ar[rr]_-{G}
&&\ACat.
\\}$$
Here $\sigma_*$ is the natural isomorphism that permutes the factors
of the categorical
product.  Elementary pasting arguments show that $\sigma^*f$ is again
a $k$-morphism.

To define the multiproduct, consider symmetric functors
$F_{ij_i}$ for $1\le i\le k$ and $1\le j_i\le n_i$ for each $i$,
also $G_i$ for $1\le i\le k$, and finally $H$.   Given
$n_i$-morphisms $f_i\colon \br{F_i}\to G_i$, where
$\br{F_i}:=(F_{i1},\ldots,F_{in_i})$, and a $k$-morphism $g\colon
\br{G}\to H$, we need to define an $(n_1+\cdots+n_k)$-morphism
$$h:=\multprod(g;f_1,\ldots,f_k)\colon \br{F}\to H.$$
Write $n$ for $n_1+\cdots+n_k$.  For each $\delta\colon
\N\overto{\cong}n\cdot\N$, we must produce a natural map
$h_\delta\colon F\to H\circ \delta^*$, subject to coherence.  Pick
bijections $\alpha_i\colon \N\overto{\cong}n_i\cdot\N$
arbitrarily, and let
$\beta=(\alpha_1\coprod\cdots\coprod\alpha_k)^{-1}\delta$, so
$\delta=(\alpha_1\coprod\cdots\coprod\alpha_k)\beta$.  We now
define $h_\delta$ as the transformation obtained from the
following gluing diagram:
$$\xymatrix@=30pt{
(\F^\infty)^n
\ar[r]^-{\cong}
\ar[ddr]_-{\delta^*}
&(\F^\infty)^{n_1}\times\cdots\times(F^\infty)^{n_k}
\ar[rr]^-{F=F_{1}\times \cdots \times F_{k}}
\ar[d]_{\alpha_{1}^{*}\times \cdots \times \alpha_{k}^{*}}
\drtwocell<\omit>{\phantom{mmmmmmmmm}(f_{1})_{\alpha_{1}}\times \cdots
\times(f_{k})_{\alpha_{k}}}
&&\ACat^{n_{1}}\times \cdots \times \ACat^{n_{k}}
\ar[d]^-{(\times)^k}
\\
&(\F^{\infty})^{k}
\ar[rr]_{G=G_{1}\times \cdots \times G_{k}}
\drrtwocell<\omit>{\;\;g_\beta}
\ar[d]^-{\beta^*}
&\relax
&\ACat^{k}
\ar[d]^-{\times}
\\
&\F^\infty
\ar[rr]_-{H}
&\relax
&\ACat.
\\}$$
It is an interesting exercise to use all the previous coherence
conditions to verify that this definition does not depend on the
choices of $\alpha_i$, and that it does itself satisfy the
necessary coherence relations for an $n$-morphism.
\enddefinition

\oldremark{\bf Remark}
The definition of symmetric functor given above is geared toward
our application rather than generality.  Our definition uses the
categorical product and basepoint (final object) restrictions to
avoid introducing the poorly behaved ``smash product'' of (based)
categories.  The correct general definition, using a symmetric
monoidal product in place of the categorical product, would drop
the conditions involving the final object (rather than replacing
them with the analogous condition for the unit) by restricting to
the full subcategory of objects $\b{n}$ where none of the
$\b{n}_{i}$ are zero.  This version of the definition then further
generalizes to the case when $\ACat$ is a multicategory.
\endoldremark

\section{From Permutative Categories to Symmetric
Functors}\label{pctosf}

Now that we have given our intermediate category of symmetric
functors, we can prove the following theorem.

\proclaim{Theorem} There is a multifunctor $J$ from permutative
categories to symmetric functors extending the construction $J$ of
Example \refnum{Jdef}.
\endproclaim

\demo{Proof} We need to give functors
$$J\colon \klin(\C_1,\ldots,\C_k;\D)\to\kmap(J\C_1,\ldots,J\C_k;J\D)
$$
which preserve the multicategory structure.

We begin by giving $J$ on objects of $\klin(\C_1,\ldots,\C_k;\D)$; for
this fix a $k$-linear map $f\colon
\C_1\times\cdots\times\C_k\to\D$ with structure maps
$\delta_i\colon f(c_i)\oplus f(c_i')\to f(c_i\oplus c_i')$ for
$1\le i\le k$.   We need an induced $k$-map $Jf\colon
(J\C_1,\ldots,J\C_k)\to J\D$.  This in turn consists of natural
functors $(Jf)_\beta$ for each choice of bijection $\beta\colon
\N\to k\cdot\N$; we need to
specify a functor
$$(Jf)_\beta\colon J\C_1\br{\b n_1}\times\cdots\times J\C_k\br{\b n_k}
\to J\D\beta^*(\br{\b n_1},\ldots,\br{\b n_k}).
$$

Again, we begin by specifying $(Jf)_{\beta}$ on objects.
An object of the source of this functor is a $k$-tuple
$(A_1,\ldots,A_k)$ where $A_i$ assigns an object $A_i\br{S_i}$ of
$\C_i$ to each sequence of subsets
$$S_{ij}\subset\{1,2,\ldots,n_{ij}\}\subset \b{n}_{ij},
$$
where $\br{\b n_i}=({\b n}_{i1},{\b
n}_{i2},\ldots)\in\obj{\F^\infty}$.  An object of the target
assigns an object of $\D$ to each sequence of subsets
$\br{T}=(T_1,T_2,\ldots)$, where
$$T_s\subset\{1,\ldots,n_{\beta(s)}\}\subset \b{n}_{\beta(s)},
$$
and we must specify such an assignment for each object
$(A_1,\ldots,A_k)$ of the source.  To do so,
for $1\leq i\leq k$, let $\br{\beta_{*}\br{T}_{i}}$ be the sequence with
$j$-th entry the subset
$$
\beta_{*}\br{T}_{ij}=T_{\beta^{-1}(i,j)}\subset\{1,\ldots,n_{ij}\}
\subset \b{n}_{ij}.
$$
Then
$(A_1\br{\beta_{*}\br{T}_1},\ldots,A_k\br{\beta_{*}\br{T}_k})$ is an
object of $\C_1\times\cdots\times\C_k$, to which we may apply our
$k$-linear map $f$ to produce an object of $\D$.  We therefore
define
$$(Jf)_\beta(A_1,\ldots,A_k)\br{T}:=
f(A_1\br{\beta_{*}\br{T}_1},\ldots,A_k\br{\beta_{*}\br{T}_k}).
$$
To complete the description of $(Jf)_\beta(A_1,\ldots,A_k)$ as an
object of $J\D\beta^*(\br{\b n_1},\ldots,\br{\b n_k})$, we must
specify the maps  $\rho_{\br{T};s,U}$.  Setting $(i,j)=\beta^{-1}(s)$,
and writing $S_{ab}=\beta_{*}\br{T}_{ab}$, we define
$\rho_{\br{T};s,U,V}$ to be the composite
$$\align
&f(A_1\br{S_1},\ldots,A_i\br{S_i\subst{j}U},\ldots,A_k\br{S_k}) \oplus
f(A_1\br{S_1},\ldots,A_i\br{S_i\subst{j}V},\ldots,A_k\br{S_k})\\
&\qquad @>\hbox to 4em{\hss$\delta_i$\hss}>>
f(A_1\br{S_1},\ldots,A_i\br{S_i\subst{j}U}\oplus
A_i\br{S_i\subst{j}V},\ldots,A_k\br{S_k})\\
&\qquad @>\hbox to 4em{\hss$f(\rho^{A_i}_{\br{S_{i};j,U,V}})$\hss}>>
f(A_1\br{S_1},\ldots,A_i\br{S_i},\ldots,A_k\br{S_k}).
\endalign
$$
The coherence requirements on the $\delta_i$'s in the definition
of a $k$-linear map are just what is needed to ensure that this
composite is again a structure map for an object of $J\D$, as the
reader may check.

We extend $(Jf)_\beta$ to morphisms of $J\C_1\br{\b
n_1}\times\cdots\times J\C_k\br{\b n_k}$ as follows: Given
$\phi_i\colon A_i\to B_i$ in $J\C_i$ for $1\le i\le k$, the
$k$-linear map $f$, being a functor, provides us with a map
$$f(\phi)\colon f(A_1\br{S_1},\ldots,A_k\br{S_k})\to
f(B_1\br{S_1},\ldots,B_k\br{S_k}),
$$
and this gives the map on morphisms for $(Jf)_\beta$.

Unpacking the definitions shows that the maps $(Jf)_\beta$ form the
components of a $k$-map, and we have therefore specified the map $J$ on
objects of $\klin(\C_1,\ldots,\C_k;\D)$.  A morphism $\phi\colon f\to g$ in
$\klin(\C_1,\ldots,\C_k;\D)$ is a natural transformation commuting
with the structure maps $\delta_i$, while a morphism from $Jf$ to $Jg$
is a coherent choice of natural transformations from $(Jf)_\beta\br{\b
n}$ to $(Jg)_\beta\br{\b n}$, the functors we have just described.
Given $\phi$ in $\klin(\C_1,\ldots,\C_k;\D)$, the natural
transformation  $(Jf)_\beta\br{\b n}\to (Jg)_\beta\br{\b n}$
is induced by the components of $\phi$ via the assignment
$$ \align
&(Jf)_\beta(A_1,\ldots,A_k)\br{T} =
f(A_1\br{\beta_*\br{T}_1},\dotsc,A_k\br{\beta_*\br{T}_k})\\
&\qquad @> \hbox to 2em{\hss$\phi$\hss} >>
g(A_1\br{\beta_*\br{T}_1},\dotsc,A_k\br{\beta_*\br{T}_k})
=(Jg)_\beta(A_1,\ldots,A_k)\br{T}.
\endalign
$$
As above, the sequence $\br{\beta_*\br{T}_i}$ is formed from $\br{T}$
by $\beta_*\br{T}_{ij}=T_{\beta^{-1}(i,j)}$.  The reader may now
unpack these definitions to find that we do indeed have a functor
as needed.

It remains to check that $J$ preserves the symmetric group
actions, the units, and the multiproduct. These verifications are
entirely straightforward given the formulas above and are left to
the reader.
\enddemo

\section{From Symmetric Functors to Symmetric Spectra}\label{sftoss}

We turn next to the description of our multifunctor from symmetric
functors in $\Cat$ to symmetric spectra.  Again, to avoid the
confusion of the different levels of functors and natural
transformations, it is convenient to work as long as possible with
symmetric functors into an arbitrary category $\ACat$ (satisfying
certain hypotheses).  The construction in this context is a
multifunctor into the multicategory of symmetric spectra in $\SimpAC$.
We begin with a review of this multicategory.

The standard definition of the category of symmetric spectra in
$\SimpAC$ in the case when $\ACat$ is the category of sets is
usually phrased in terms of the smash product of based simplicial
sets.  We formulate the definition in terms of the cartesian
product with additional base point conditions, since this works
better when $\ACat$ is the category of small categories.  The
formulation of the category of symmetric spectra that follows
should therefore be thought of as a first step in the construction
of the functor $K$ rather than a generalization of the category of
symmetric spectra of \cite{7}.

\definition{Definition}\label{defsymsp}
Let $\ACat$ be a category with finite products and coproducts, and
assume that for all $X$ in $\ACat$, the functor $X\times(-)$
preserves coproducts (as, for example, when $\ACat$ is cartesian
closed). Let $\SimpAC$ denote the category of simplicial objects
in $\ACat$, i.e., contravariant functors from the category
$\Delta$ of simplices to the category $\ACat$.  We use $*$ to
denote both the final object of $\ACat$ and also the final object
in $\SimpAC$, the constant simplicial object on $*$.  For $X$ in
$\SimpAC$ and $K$ a finite simplicial set, write $X\times K$ for
the tensor of $X$ with $K$; concretely, $X \times K$ has
$n$-simplices
$$ (X\times K)_{n} = \coprod_{K_{n}} X_{n}. $$
A {\bf symmetric spectrum} in $\SimpAC$ consists of objects $X(p)$ in
$\SimpAC$ for all non-negative integers $p$, an action of the
symmetric group $\Sigma_{p}$ on $X(p)$, and maps
$$ *\to X(p)\qquad \text{and}\qquad
X(p)\times \SSp\to X(p+1),
$$
such that the map $*\to X(p)$ preserves the $\Sigma_{p}$-action
(with the trivial action on $*$), for each $q\geq 1$ the composite
$X(p)\times (\SSp)^{q}\to X(p+q)$ preserves the $(\Sigma_{p}\times
\Sigma_{q})$-action, and the composites
$$ X(p)\times *\to X(p)\times \SSp\to X(p+1)
\qquad \text{and}\qquad
* \times \SSp \to X(p)\times \SSp\to X(p+1)
$$
factor through the given map $*\to X(p+1)$.

A $k$-morphism in symmetric spectra in $\SimpAC$ from
$X_{1},\dotsc,X_{k}$ to $Y$ consists of maps
$$ X_{1}(p_{1}) \times \dotsb \times X_{k}(p_{k})\to
Y(p_{1}+\dotsb +p_{k}) $$
for all $p_{1},\dotsc,p_{k}$ that preserve the $\Sigma_{p_{1}}\times
\dotsb \times \Sigma_{p_{k}}$ action and that make the following diagrams
commute for all $1\leq i\leq k$:
$$ \vcenter{
\xymatrix@C-10pt{
X_1(p_{1})\times \dotsb
%\times X_{i-1}(p_{i-1})
\times *
%\times X_{i+1}(p_{i+1})
\times \dotsb \times X_k(p_{k})\ar[d]\ar[r]
&\relax*\ar[d]\\
X_1(p_{1})\times \dotsb
%\times X_{i-1}(p_{i-1})
\times X_i(p_{i})
%\times X_{i+1}(p_{i+1})
\times \dotsb \times X_k(p_{k})\ar[r]
&Y(p_{1}+\dotsb+p_{k})
}} \leqno{\text{\bf \refnum{defsymsp}(a)}}
$$
$$ \vcenter{\xymatrix{
(X_1(p_{1})\times \dotsb \times X_k(p_{k}))\times\SSp\ar[r]\ar[d]
&Y(p_{1}+\dotsb +p_{k})\times \SSp\ar[d]\\
X_1(p_{1})\times \dotsb \times (X_i(p_{i})\times\SSp)\times \dotsb
\times X_k(p_{k})
\ar[d]
&Y(p_{1}+\dotsb+p_{k}+1)\ar[d]^{c_{i}}\\
X_1(p_{1})\times \dotsb \times X_i(p_{i}+1)\times \dotsb \times
X_k(p_{k})
\ar[r]
&Y(p_{1}+\dotsb+p_{k}+1),
}}\leqno{\text{\bf \refnum{defsymsp}(b)}}
$$
where $c_{i}$ denotes the permutation in $\Sigma_{p_{1}+\dotsb
+p_{k}+1}$ that moves the last element to the
$(p_{1}+\dotsb+p_{i}+1)$-st position but otherwise preserves the
order, i.e., the cycle $(q+1,\dotsc,p,p+1)$ where
$q=p_{1}+\dotsb+p_{i}$ and $p=p_{1}+\dotsb+p_{k}$.  The
$\Sigma_{k}$ action on the $k$-morphisms is induced by permuting
the product factors and the symmetric group action on the target,
permuting blocks.  The identity $1$-morphisms are the
$1$-morphisms induced by the identity maps.  The multiproduct is
induced by products and compositions in $\ACat$.
\enddefinition

By the simplicial nature of the construction, the multicategory is
enriched over simplicial sets.  When $\ACat$ is enriched over
small categories or simplicial sets, the conditions in the
previous definition translate into limits on the categories or
simplicial sets of maps, and the multicategory of symmetric
spectra in $\SimpAC$ becomes enriched over simplicial categories
or bisimplicial sets.

\proclaim{Proposition}
The multicategory of symmetric spectra in simplicial sets as
defined above is isomorphic to the multicategory associated to the
symmetric monoidal category of symmetric spectra of \cite{7}.
\endproclaim

\demo{Proof}
This is an easy consequence of the definition of the smash product
of simplicial sets and the external formulation of the smash
product of symmetric spectra.  Technically, the paper \cite{7}
considers the category of ``left $S$-modules'' whereas the
(external) formulation above specifies the category of right
$S$-modules, but the identity isomorphism $S\iso S^{\op}$ induces
a strong symmetric monoidal isomorphism between these categories.
\enddemo

Now we describe the multifunctor from symmetric functors in
$\ACat$ to symmetric spectra in $\SimpAC$.  Recall from
\refname{secK} that we have defined our model of the circle $\SSp$
so that its based set of $n$-simplices is $\b n$, giving $\SSp$ as
a functor from $\Deltaop$ to $\F$.

\definition{Construction}
For $F$ a symmetric functor, let $IF(0)=F(\b1,\b1,\dotsc)$, let
$IF(1)$ be the simplicial object $F(\SSp)$ using the canonical
inclusion $\F\to\F^{\infty}$, and for $p>1$, let $IF(p)$ be the
diagonal simplicial object on the multisimplicial object
$F(\SSp,\dotsc,\SSp)$ using the canonical inclusion $\F^{p}\to
\F^{\infty}$.  We give $IF(p)$ the $\Sigma_{p}$ action arising
from the action of $\Sigma_{p}$ on $\F^{p}$, or more accurately,
its extension to the action of $\Aut(\N)$ on $\F^{\infty}$ fixing
the numbers greater than $p$. We have a canonical map
$*=F(\b0,\dotsc,\b0)\to IF(p)$ induced by the initial map in
$\F^{p}$, and this map preserves the $\Sigma_{p}$-action.  We have
maps
$$ IF(p)\times \SSp\to IF(p+1) $$
induced by the maps
$$(\b{n}_{1},\dotsc,\b{n}_{p},\b1,\b1,\dotsc )\times \b{n}_{p+1}
\to (\b{n}_{1},\dotsc,\b{n}_{p},\b{n}_{p+1},\b1,\dotsc )
$$
in $\F^{\infty}$ which for each $x$ in $\b{n}_{p+1}$ sends the $\b1$
in the $(p+1)$-st position to $\b{n}_{p+1}$ by the unique based map
that takes $1$ to $x$.  The composite map
$$ IF(p)\times (\SSp)^{q}\to IF(p+q) $$
has a similar description and so is easily seen to be
$\Sigma_{p}\times \Sigma_{q}$ equivariant.  Since $F\br{\b n}$ is
$*$ whenever any of the $\b{n}_{i}$ is $\b0$, the maps
$IF(p)\times
\SSp\to IF(p+1)$ restrict on $IF(p)\times *$ and $*\times \SSp$ to
the final map composed with the given map $*\to IF(p+1)$.  It
follows that these objects and maps assemble to a symmetric
spectrum; we denote this symmetric spectrum as $IF$.
\enddefinition

\proclaim{Theorem}
$I$ extends to a multifunctor from the multicategory of symmetric functors
in $\ACat$ to the multicategory of symmetric spectra in $\SimpAC$.
\endproclaim

\demo{Proof}
Let $F_{1},\dotsc,F_{k}$ and $G$ be symmetric functors and consider
a $k$-morphism $f$ from $F_{1},\dotsc,F_{k}$ to $G$.  We obtain a
$k$-morphism from $IF_{1},\dotsc,IF_{k}$ to $IG$ using the map
$$ f_{\beta}\colon
F_{1}(\underbrace{\SSp,\dotsc,\SSp}_{p_{1}}) \times \dotsb \times
F_{k}(\underbrace{\SSp,\dotsc,\SSp}_{p_{k}})
\to
G(\underbrace{\SSp,\dotsc,\SSp}_{p_{1}+\dotsb +p_{k}})
$$
where $\beta\colon \N\to k\cdot \N$ is any bijection that takes
$$ 1,\dotsc,p_{1}+\dotsb+p_{k}
\qquad \text{to}\qquad
(1,1),\dotsc,(1,p_{1}),\dotsc, (k,1),\dotsc,(k,p_{k}) $$ in
lexicographical order (all such $\beta$ have identical $f_{\beta}$
when restricted to $\F^{p_{1}}\times \dotsb \times \F^{p_{k}}$).
The equivariance condition follows from \refnum{sfkmap}(a). The
final object diagram
\refnum{defsymsp}(a) commutes because
$G\br{\b{n}}$ is $*$ whenever any $\b{n}_{i}$ is $\b0$.

Next we verify the suspension diagram \refnum{defsymsp}(b).  For this,
fix a bijection $\beta$ satisfing the condition above and the
additional condition $\beta(p_1+\cdots+p_k+1)=(i,p_i+1)$; this ensures
that $f_{\beta\circ c_i}$ defines the $k$-morphism after we suspend on
$F_i$.  Fix the simplicial degree $n$, and write $\brn{j}$ for
$$ (\underbrace{\b{n},\ldots,\b{n}}_{j},\b1,\b1,\ldots). $$
We write $\sigma$ for the map $F_{i}\brn{p_{i}}\times\b{n}\to
F_{i}\brn{p_{i}+1}$ and the map  $G\brn{p}\times \b n\to
G\brn{p+1}$.
Now we have the diagram
$$\xymatrix{
F_1\brn{p_1}\times\cdots\times F_k\brn{p_k}\times\b{n}
\ar[r]^-{f_\beta\times1}\ar[d]_-{1\times\tau}
&G\brn{p}\times\b{n}\ar[d]^-{\sigma}
\\F_1\brn{p_1}\times\cdots\times(F_i\brn{p_i}\times\b{n})
\times\cdots\times F_k\brn{p_k}\ar[d]_-{1\times\sigma\times1}
&G\brn{p+1}\ar[d]^-{(c_i)_!}
\\F_1\brn{p_1}\times\cdots\times F_i\brn{p_{i}+1}
\times\cdots\times F_k\brn{p_k}
\ar[ur]_-{f_\beta}\ar[r]_-{f_{\beta\circ c_i}}
&G\brn{p+1},}
$$
in which the upper part commutes by naturality of $f_\beta$, and
the lower part by \refnum{sfkmap}(b).  This diagram is precisely the suspension compatibility diagram
\refnum{defsymsp}(b) in simplicial degree $n$.

The description above therefore specifies a $k$-morphism
of symmetric spectra. We leave to the reader the exercise of
correlating definitions to check that this association preserves
the symmetric group action on the $k$-morphisms, the units, and
the multiproduct.

When we regard the $k$-morphisms of symmetric functors as discrete
simplicial sets, the multicategory of symmetric functors is
enriched over simplicial sets and the multifunctor described above
is enriched (for trivial reasons).  When $\ACat$ is enriched over
small categories or simplicial sets, we can regard the
multicategory of symmetric functors as enriched over simplicial
categories or bisimplicial sets by taking the (other) simplicial
direction to be discrete.  A straightforward check then shows that
the multifunctor described above is enriched over simplicial
categories or bisimplicial sets.
\enddemo

Composing the multifunctor $J$ from the previous section, the
multifunctor $I$, the nerve functor, and the diagonal functor (from
bisimplicial sets to simplicial sets), we obtain a multifunctor $K$
from the category of small permutative categories to the category of
symmetric spectra.  By inspection, the underlying functor is naturally
isomorphic to $\Knew$.  This completes the proof of \refname{num6}.

\section{Associative Categories, Bipermutative Categories, and
the Operads $\Sigma_{*}$ and $E\Sigma_{*}$}\label{secassoc}

This section is devoted to the proofs of Theorems \refnum{num10}
and \refnum{num12}.

\demo{Proof of \refname{num10}}
First, suppose we are given a small associative category $\A$; we
must produce a multifunctor $\Sigma_*\to\P$ sending the single
object of $\Sigma_*$ to $\A$. In this case, a multifunctor as
specified in the theorem is precisely a map of operads (in $\Cat$)
from $\Sigma_*$ to the endomorphism operad of $\A$ in $\P$, whose
component categories are the $k$-linear maps
$\P_k(\A,\dotsc,\A;\A)$. In other words, we must define a sequence
of functors $T_k\colon \Sigma_k\to\P_k(\A,\dotsc,\A;\A)$, and show
that they specify a map of operads.  Since $\Sigma_{k}$ is a
discrete category, specifying the functor $T_{k}$ is equivalent to
specifying a $k$-morphism $T_k\sigma$ for every element $\sigma$
in the group $\Sigma_{k}$. As per \refname{num8}, the $k$-morphism
$T\sigma$ consists of a functor $f^{\sigma}\colon \A^{k}\to \A$
and natural distributivity maps $\delta^{\sigma}_{i}$ for $1\leq
i\leq k$.

We define $f^{\sigma}$ by
$$f^\sigma(a_1,\ldots,a_k)=a_{\sigma^{-1}(1)}\otimes\cdots\otimes
a_{\sigma^{-1}(k)}.
$$
For notational convenience in defining $\delta^\sigma_i$, let
$P=a_{\sigma^{-1}(1)}\otimes\cdots\otimes
a_{\sigma^{-1}(\sigma(i)-1)}$, and
$Q=a_{\sigma^{-1}(\sigma(i)+1)}\otimes\cdots\otimes
a_{\sigma^{-1}(k)}$.  We then define $\delta^\sigma_i$  as the
common diagonal of the following square, which commutes by
\refname{num9}, condition (e):
$$\xymatrix{
(P\otimes a_i\otimes Q)\oplus(P\otimes a'_i\otimes
Q)\ar[r]^-{d_r}\ar[d]_-{d_l}
&P\otimes((a_i\otimes Q)\oplus(a'_i\otimes Q))\ar[d]^-{1\otimes d_l}\\
((P\otimes a_i)\oplus(P\otimes a'_i))\otimes
Q\ar[r]_-{d_r\otimes1}
&P\otimes(a_i\oplus a_i')\otimes Q.\\
}
$$
The reader may now verify that the requirements for distributivity
maps are satisfied.

We must verify that the $T_k$'s give a map of operads.
Equivariance is elementary; we check preservation of the
multiproduct. This follows as a consequence of the following
commutative diagram, where $\sigma\in\Sigma_k$ and
$\phi_i\in\Sigma_{j_i}$ for $1\le i\le k$:
$$\xymatrix@!C=65pt{
\A^{j_1}\times\cdots\times\A^{j_k}
\ar[dr]^-{{\phantom{mm}}f^{\phi_1}\times\cdots\times f^{\phi_k}}
\ar[d]_-{\phi_1\times\cdots\times\phi_k}\\
\A^{j_1}\times\cdots\times\A^{j_k}\ar[r]_-{\otimes^k}
\ar[d]_-{\sigma\langle j_1,\ldots,j_k\rangle}
&\A^k\ar[dr]^-{f^\sigma}\ar[d]^-{\sigma}\\
\A^{j_{\sigma^{-1}(1)}}\times\cdots\times\A^{j_{\sigma^{-1}(k)}}
\ar[r]_-{\otimes^k}
&\A^k\ar[r]_-{\otimes}&\A.
}
$$
We must also check that the distributivity maps of
$\multprod(T\sigma;T\phi_1,\ldots,T\phi_k)$ coincide with those of
$T\multprod(\sigma;\phi_1,\ldots,\phi_k)$.  However, both distribute
to the same ending point, which may be written
$$P_1\otimes P_2\otimes(a_i\oplus a_i')\otimes Q_2\otimes Q_1,
$$
where $P_1$ is the tensor product of blocks preceding the one in
which $a_i\oplus a_i'$ appears, and $P_2$ is the tensor product of
the terms in the same block which precede $a_i\oplus a_i'$.  $Q_1$
and $Q_2$ are described analogously.  Now
$\multprod(T\sigma;T\phi_1,\ldots,T\phi_k)$ distributes first
$P_1$ and $Q_1$, and then $P_2$ and $Q_2$, while
$T\multprod(\sigma;\phi_1,\ldots,\phi_k)$ does it all at once. The
resulting maps coincide by property (d) of the distributivity maps
in Definition \refnum{num9}. Therefore $T$ preserves the
multiproduct, and we get a map of operads, i.e., a multifunctor
$T\colon
\Sigma_*\to\P$.

Now suppose given a map of operads
$T\colon \Sigma_*\to\{\P_k(\A^k;\A)\}$; we must produce an associative
structure on $\A$.  First, the tensor product functor
$\otimes\colon \A^2\to\A$ is the functor part of the image of
$1\in\Sigma_2$, and the unit object is the image of the unique
element of $\Sigma_0$.  Write $1_n$ for the identity element of
$\Sigma_n$.  Then the strict associativity of $\otimes$ follows
from the fact that $\multprod(1_2;1_2,1_1)=1_3=\multprod(1_2;1_1,1_2)$,
and the unit condition follows from
$\multprod(1_2;1_1,1_0)=1_1=\multprod(1_2;1_0,1_1)$.

The distributivity maps $d_l$ and $d_r$ arise as part of the
structure of the target of $1_2\in\Sigma_2$. Properties (a), (b),
(c), and (f) follow immediately from requirements for
$k$-morphisms in $\P$.  Properties (d) and (e) follow from the
facts that $T$ is a map of operads, and also that
$\multprod(1_2;1_1,1_2)=\multprod(1_2;1_2,1_1)$.  The distributivity
maps for the images of these composites must therefore coincide,
and both (d) and (e) follow.  We therefore have an associative
structure whenever we have a map of operads
$\Sigma_*\to\{\P_k(\A^k;\A)\}$.

Finally, we must verify that these correspondences are inverse to each
other.  First suppose given an associative structure on $\A$, and
let $T\colon \Sigma_*\to\{\P_k(\A^k;\A)\}$ be the induced map of
operads. By definition, $T(1_2)$ is the tensor product on $\A$,
together with both distributivity maps, and the multiplicative
unit is given by $T(1_0)$. We therefore recover the original
structure from its induced map of operads.

Now suppose we start with a map of operads
$T\colon \Sigma_*\to\{\P_k(\A^k;\A)\}$, and give $\A$ the induced
associative structure.  By induction using the fact that
$\multprod(1_2;1_{k-1},1_1)=1_k$, we find that
$$f^{1_{k}}(a_1,\ldots,a_k)=a_1\otimes\cdots\otimes a_k,$$
and from equivariance it follows that, for $\sigma\in\Sigma_k$,
$$f^{\sigma}(a_1,\ldots,a_k)=a_{\sigma^{-1}(1)}\otimes\cdots\otimes
a_{\sigma^{-1}(k)}.
$$
We therefore recover the map of operads $T$ on underlying
functors $f$, and we are left with the recovery of the distributivity
maps. By equivariance, it suffices to recover the distributivity
maps $\delta_i^{1_k}$, which we do by induction on $k$.  This is
trivial if $k\le2$. Since $T$ is a map
of operads, we have
$$\multprod(T(1_2);T(1_i),T(1_{k-i}))=T(1_k).$$
If $i<k$, assume by induction that $\delta_i^{1_i}$ is given by
$$\xymatrix{
(P\otimes a_i)\oplus(P\otimes
a'_i)\ar[r]^-{d_r}&P\otimes(a_i\oplus a'_i). }
$$
Then by the definition of distributivity maps in the multiproduct
$\multprod(T(1_2);T(1_i),T(1_{k-i}))$, we have $\delta_i^{1_k}$ given
by the composite
$$\xymatrix{
(P\otimes a_i\otimes Q)\oplus(P\otimes a'_i\otimes Q)\ar[r]^-{d_l}
&((P\otimes a_i)\oplus(P\otimes a'_i))\otimes Q\ar[r]^-{d_r\otimes1}
&P\otimes(a_i\oplus a'_i)\otimes Q,}
$$
as required.  In the remaining case, where $i=k$, we use the fact
that the (single) distributivity map of $T(1_1)$ is the identity,
together with
$$\multprod(T(1_2);T(1_{k-1}),T(1_1))=T(1_k),$$
to exhibit $\delta_k^{1_k}$ as simply
$$\xymatrix{(P\otimes a_k)\oplus(P\otimes
a'_k)\ar[r]^-{d_r}&P\otimes(a_k\oplus a'_k), }
$$
as required.  This completes the proof.
\enddemo

\demo{Proof of \refname{num12}}
First suppose given a map of operads
$E\Sigma_*\to\{\P_k(\R^k;\R)\}$.  Then we have the composite
multifunctor
$$\xymatrix{\Sigma_*\ar[r]&E\Sigma_*\ar[r]^-{\R}&\P,}
$$
so by \refname{num10}, $\R$ is associative.  We therefore get all of
the bipermutative structure except for:
\roster
\item $\gamma^\otimes$,
\item The coherence diagram for $\gamma^{\otimes}$ from the
requirement that $(\R,\otimes,1)$ form
a permutative category, and
\item Diagram (e$'$).
\endroster
The symmetry
isomorphism $\gamma^\otimes$ is the image of the isomorphism
between the two objects of $E\Sigma_2$. The coherence diagram
$$\xymatrix{
a\otimes b\otimes c
\ar[rr]^-{\gamma^\otimes}\ar[dr]_-{1\otimes\gamma^\otimes}
&&c\otimes a\otimes b
\\&a\otimes c\otimes b
\ar[ur]_-{\gamma^\otimes\otimes1}}
$$
now follows as a consequence of there being exactly one
isomorphism in $E\Sigma_3$ between $1_3\in\Sigma_3$ and the
permutation sending $(abc)$ to $(cab)$.  Diagram (e$'$) is simply
the requirement that $\gamma^\otimes$, being the image of a
morphism in $E\Sigma_2$, must be a morphism in $\P_2(\R^2;\R)$. A
map of operads $E\Sigma_*\to\{\P_k(\R^k;\R)\}$ therefore
determines a bipermutative structure on~$\R$.

Suppose now that we are given that $\R$ is a small bipermutative
category; we need to construct the multifunctor
$T\colon E\Sigma_*\to\P$.  From \refname{num10}, we get the map of operads on
the objects $\Sigma_*$ once we know that $\R$ is an associative
category, and the only issue here is diagram (e) in \refname{num9},
which we have replaced with (e$'$). However, diagram (e) follows
as a consequence of the commutativity of the diagram in Figure~1%
\refonsamepage{figure}{}{ (see page~\refpage{figure})},
all of whose subdiagrams are instances of the coherence
requirements for a bipermutative category.
\midinsert
$$\xymatrix@C=-15pt{
(a\otimes b\otimes c)\oplus(a\otimes b'\otimes
c)\ar[rrr]^-{d_l}\ar[ddr]^-{\gamma\oplus\gamma}
\ar@/_.7pc/[ddddr]_-{(1\otimes\gamma)^{\oplus2}}\ar[ddddddd]_-{d_r}
&&&((a\otimes b)\oplus(a\otimes b'))\otimes
c\ar[ddddddd]^-{d_r\otimes1}\\
&&c\otimes((a\otimes b)\oplus (a\otimes
b'))\ar[ur]^-{\gamma}\ar[dd]^-{1\otimes d_r}\\
&(c\otimes a\otimes b)\oplus(c\otimes a\otimes
b')\ar[ur]^-{d_r}\ar[dr]^-{d_r}\ar[dd]^-{(\gamma\otimes1)^{\oplus2}}\\
&&c\otimes a\otimes(b\oplus
b')\ar[ddddr]^-{\gamma}
\ar[dd]_-{\gamma\otimes1}\\
&(a\otimes c\otimes b)\oplus(a\otimes c\otimes
b')\ar[dr]^-{d_r}\ar[dd]_-{d_r}\\
&&a\otimes c\otimes(b\oplus b')\ar[ddr]_-{1\otimes\gamma}\\
&a\otimes((c\otimes b)\oplus(c\otimes
b'))\ar[ur]_-{1\otimes d_r}\ar[dl]_-{1\otimes(\gamma\oplus\gamma)}\\
a\otimes((b\otimes c)\oplus(b'\otimes c))\ar[rrr]_-{1\otimes
d_l}&&&a\otimes(b\oplus
b')\otimes c\\
}
$$
\botcaption{Figure 1}
\endcaption
\labelpage{figure}
\endinsert
We therefore get a map of operads
$T\colon \Sigma_*\to\{\P_k(\R^k;\R)\}$, and it remains to extend this to
the morphisms in the $E\Sigma_k$'s.  These consist of one
isomorphism between each pair of objects.  Given any pair of
elements $\sigma$ and $\phi$ in $\Sigma_k$, the permutative
structure on $(\R,\otimes,1)$ gives a canonical isomorphism
$$\xymatrix{
a_{\sigma^{-1}(1)}\otimes\cdots\otimes
a_{\sigma^{-1}(k)}\ar[r]^{\cong}
&a_{\phi^{-1}(1)}\otimes\cdots\otimes a_{\phi^{-1}(k)},
}$$ as a composite of the maps $\gamma^{\otimes}$; we take this as
the image of the unique morphism from $\sigma$ to $\phi$. The
coherence condition for $\gamma^{\otimes}$ implies that any ways
of composing various instances of $\gamma^{\otimes}$ that lead to
the same permutation of the tensor factors give the same
isomorphism; we use this fact multiple times below, and refer to
it as ``uniqueness of the permutation isomorphisms''.
Compatibility of these permutation isomorphisms with the given
distributivity maps follows from coherence of the bipermutative
structure, specifically property (e$'$) using the fact that
$\Sigma_k$ is generated by transpositions. The uniqueness of the
permutation isomorphisms implies that $T_{k}$ is a functor
$E\Sigma_{k}\to
\P_{k}(\R^k;\R)$.  In order to see that $T$ defines a map
of operads on the morphisms, we apply a little more coherence
theory. Given objects $(\sigma;\phi_1,\ldots,\phi_k)$ and
$(\sigma';\phi_1',\ldots,\phi_k')$ of $E\Sigma_k\times
E\Sigma_{j_1}\times\cdots\times E\Sigma_{j_k}$, there is a unique
isomorphism from one to the other in $E\Sigma_k\times
E\Sigma_{j_1}\times\cdots\times E\Sigma_{j_k}$.  The target of
this morphism under $\multprod T$ first permutes within blocks,
and then permutes the blocks, while the target under $T\multprod$
does this all at once; these are the same isomorphism by the
uniqueness of the permutation isomorphisms. This concludes the
proof that $T$ is a map of operads, and consequently the given
data determine a multifunctor $E\Sigma_*\to\P$.  The proof that
the passages back and forth are inverse to each other is exactly
as in the proof of \refname{num10}.
\enddemo

\section{Modules and Algebras in Permutative Categories}
\label{secmodules}

In this section, we describe some of the module and algebra
structures in $\P$, the multicategory of permutative categories.
We first define each structure in terms of functors and natural
transformations; we then reinterpret the structure in terms of
parameter multicategories.  All of the parameter multicategories
we describe below have contractible components in their
$k$-morphism categories, so collapsing each component to a single
point gives a map of multicategories that is the identity on
objects and a weak equivalence on $k$-morphisms. From Theorems
\refnum{num5} and
\refnum{num5.5}, it follows that the structures we describe pass to
the associated strict structures on $K$-theory spectra.

\subsection{Modules}\nobreak

\definition{Definition}\label{num18} Let $\A$ be an associative category and
$\D$ a permutative category.  A {\bf left
$\A$-module} structure on $\D$ consists of a functor $\otimes\colon
\A\times\D\to\D$ that
is strictly associative in the sense that the diagram
$$\xymatrix{
\A\times\A\times\D\ar[r]^-{1\times\otimes}\ar[d]_-{\otimes\times1}
&\A\times\D\ar[d]^-{\otimes}
\\ \A\times\D\ar[r]_-{\otimes}&\D}
$$
commutes, strictly unital in the sense that the composite
$$\xymatrix{
\D\cong\{1\}\times\D\ar[r]&\A\times\D\ar[r]^-{\otimes}&\D}
$$
coincides with the identity, together with natural distributivity maps
$$d_l\colon (a\otimes d)\oplus(a'\otimes d)\to(a\oplus a')\otimes d
$$
and
$$d_r\colon (a\otimes d)\oplus(a\otimes d')\to a\otimes(d\oplus d')
$$
subject to the commutativity of all the diagrams in \refname{num9}.
\enddefinition

Since an associative category structure on a small permutative
category is equivalent to the structure of an algebra over an operad
$\Sigma_{*}$, we have the notion of a left module described in terms of the
parameter multicategory for left modules  discussed as the
third example following
\refname{num3}.  We repeat this here for convenience.

\definition{Definition}\label{num17} The multicategory $\Msig$ is the
following parameter multicategory for modules: It has two objects,
$A$ (the ``ring'') and $M$ (the ``module'').  In the case in which
all inputs and the output are $A$, we have
$\Msig_k(A^k;A)=\Sigma_k$, and if exactly one input is $M$ and the
output is also $M$, we set
$\Msig_k(A^{j-1},M,A^{k-j};M)=\{\sigma\in\Sigma_k:\sigma(j)=k\}$.
All other $k$-morphism sets are required to be empty. The
multiproduct and $\Sigma_*$-action are defined in exactly the same
way as in the operad $\Sigma_*$.
\enddefinition

Note that restricting our attention to the single
object $A$ gives a multifunctor
$$\Sigma_*\to\Msig,
$$
so if we have a multifunctor $\Msig\to\P$, the image of $A$ is an
associative category. The fundamental theorem about left module
structures on permutative categories is the following:

\proclaim{Theorem}\label{num19} Left $\A$-module structures on $\D$
determine and are determined by multifunctors $\Msig\to\P$ sending
$A$ to $\A$ and $M$ to $\D$ such that the restriction
$$\Sigma_*\to\Msig\to\P
$$
gives the structure map for $\A$ as an associative category.
\endproclaim

\demo{Proof}  First suppose given a left $\A$-module structure on
$\D$; we must produce a multifunctor $T\colon \Msig\to\P$.  The
associative structure on $\A$ gives us the multifunctor on the
$k$-morphisms of $\Msig$ involving only $A$, so consider
$\sigma\in\Msig_k(A^{j-1},M,A^{k-j};M)$, i.e., $\sigma\in\Sigma_k$
and $\sigma(j)=k$.  We define
$$T\sigma\colon \A^{j-1}\times\D\times\A^{k-j}\to\D
$$
by the formula
$$T\sigma(a_1,\ldots,a_{j-1},d,a_{j+1},\ldots,a_k)
=a_{\sigma^{-1}(1)}\otimes\cdots\otimes a_{\sigma^{-1}(k-1)}
\otimes d.
$$
Since $\sigma(j)=k$, all of the objects
$a_{\sigma^{-1}(1)},\ldots,a_{\sigma^{-1}(k-1)}$ are indeed
objects of $\A$, and this formula is simply a special instance of
the usual formula
$$T\sigma(b_1,\ldots,b_k)=b_{\sigma^{-1}(1)}\otimes\cdots\otimes
b_{\sigma^{-1}(k)}.
$$
The proof that this formula determines a multifunctor now proceeds
exactly as in the proof of \refname{num10}.

On the other hand, given a multifunctor $\Msig\to\P$ sending $A$ to $\A$
and $M$ to $\D$, and which restricts on $A$ to the associative
category structure map for $\A$, we must produce a left
$\A$-module structure on $\D$.  The tensor pairing
$\otimes\colon \A\times\D\to\D$ is the image of the single element of
$\Msig(A,M;M)$, and the distributivity maps are part of the
structure of the target of this element.  The rest of the proof
now follows exactly as in the proof of \refname{num10}.
\enddemo

Since the multicategories $\Sigma_*$ and $\Msig$ are discrete, we
do not need to apply \refname{num5.5}, and we have the following
result.

\proclaim{Corollary}\label{num19.5} If $\D$ is a left $\A$-module,
then $K\D$ is a left $K\A$ module.
\endproclaim

When $\A$ is not just associative but actually bipermutative, we
can describe a parameter multicategory that captures this further
structure using the translation category construction $E$ applied
to $\Msig$: for a multicategory of sets $\M$, let $E\M$ denote the
multicategory enriched over small categories for which
$E\M_k(B_1,\ldots,B_k;C)$ is the category obtained by applying $E$
to $\M_k(B_1,\ldots,B_k;C)$.  There is an obvious inclusion of
multicategories $\M\to E\M$, where we consider $\M$ enriched over
small categories with all the categories discrete.

\proclaim{Lemma}\label{num20}  Let $\Sigma_*\to\Msig$ be the
inclusion of the $k$-morphisms of $\Msig$ involving only $A$.
Then the diagram of multicategories
$$\xymatrix{
\Sigma_*\ar[r]\ar[d]&\Msig\ar[d]
\\ E\Sigma_*\ar[r]&E\Msig}
$$
is a pushout.  In other words, making the $k$-morphisms in
$\Sigma_*$ all canonically isomorphic forces all the other
$k$-morphisms in $\Msig$ to be canonically isomorphic as well.
\endproclaim

\demo{Proof} Let $\Q$ be another multicategory, and suppose we
have a commutative diagram
$$\xymatrix{
\Sigma_*\ar[r]\ar[d]&\Msig\ar[d]
\\ E\Sigma_*\ar[r]&\Q}
$$
of multicategories.  We must show that there is a unique dashed
arrow making the diagram of multicategories
$$\xymatrix{
\Sigma_*\ar[r]\ar[d]&\Msig\ar[d]\ar[ddr]
\\ E\Sigma_*\ar[r]\ar[drr]&E\Msig\ar@{-->}[dr]
\\ &&\Q}
$$
commute.  Certainly there is no choice about the values on the
objects of the $k$-morphism category $E\Msig_k(B_1,\ldots,B_k;C)$,
since the objects are the same as the objects of $\Msig$. The
values on morphisms of $E\Sigma_*$ are also determined. We show
that whenever $\sigma_{1}$ and $\sigma_{2}$ are objects in
$E\Msig_k(A^{j-1},M,A^{k-j};M)$, the image of the map from
$\sigma_{1}$ to $\sigma_{2}$ is also determined. Since
$\sigma_{2}\circ \sigma_1^{-1}$ fixes $k$, we can think of it as
an element of $\Sigma_{k-1}$, and let $\phi$ be the unique map in
$E\Sigma_{k-1}$ from the identity permutation to $\sigma_{2}\circ
\sigma_1^{-1}$.  Then we can express the unique map from
$\sigma_{1}$ to $\sigma_{2}$ in $E\Msig_k(A^{j-1},M,A^{k-j};M)$ by
the formula
$$ \multprod(\id_{\xi};\phi,1_{M})\cdot \sigma_{1} $$
where $\xi$ is the single object of $\Msig_2(A,M;M)$. This
establishes uniqueness of such a multifunctor, and it remains to
show existence. Using the formula above to define the functors, it
is straightforward to show that they preserve the symmetric group
action and the multiproduct and therefore define a multifunctor
$E\Msig\to \Q$.
\enddemo

\proclaim{Corollary}\label{num21}  Let $\R$ be a small bipermutative category,
$\D$ a small permutative category.  Then left $\R$-module structures on
$\D$ determine and are determined by multifunctors $E\Msig\to\P$
sending $M$ to $\D$ and restricting on $A$ to the bipermutative
structure map $E\Sigma_*\to\P$ for $\R$.
\endproclaim

\demo{Proof} This follows immediately from \refname{num20} with $\Q$
replaced by $\P$.
\enddemo

\proclaim{Corollary}\label{num21.5} If $\D$ is a left module over
a bipermutative category $\R$, then $K\D$ is equivalent to a
strict module over a strictly commutative ring spectrum equivalent
to $K\R$.
\endproclaim

For right modules, the relevant definitions are as follows.

\definition{Definition}\label{num22}  Let $\A$ be a associative category,
$\D$ a permutative category.  Then the structure of a {\bf right
$\A$-module} on $\D$ consists of a functor $\otimes\colon \D\times\A\to\D$ that
is strictly associative and unital in the analogous sense as in
\refname{num18}, together with distributivity maps again defined
analogously and satisfying the corresponding diagrams.
\enddefinition

\definition{Definition}\label{num23} The multicategory $\rMsig$ is the
following parameter multicategory for modules: It has two objects, $A$ and
$M$, with $k$-morphism sets being empty unless all inputs are $A$
and the output is $A$ or exactly one input is $M$ and the
output is $M$.  In the first case, the $k$-morphisms are
$\Sigma_k$, so the endomorphism operad of $A$ is $\Sigma_*$ (as in
$\Msig$), but we set
$$\rMsig_k(A^{j-1},M,A^{k-j};M)
=\{\sigma\in\Sigma_k:\sigma(j)=1\}.
$$
The $\Sigma_*$-action and multiproduct are defined exactly as in
$\Sigma_*$.
\enddefinition

\proclaim{Theorem}\label{num24}  Let $\A$ be a small associative category and
$\D$ a small permutative category.  Then right $\A$-module
structures on $\D$ determine and are determined by multifunctors
$\rMsig\to\P$ sending $M$ to $\D$ and restricting on $A$ to the
structure map for $\A$ as an associative category.
\endproclaim

The proof is safely left to the reader, given the proof of
\refname{num19}.  The obvious analog to Corollaries \refnum{num19.5}, \refnum{num21},
and \refnum{num21.5} also hold.

Just as in ordinary algebra,
a right module over $\A$ is the same thing as a left module over
the opposite structure ``$\A^\op$'', which we now define.

\definition{Definition}\label{num14} The {\bf opposite} map is the
particular map of operads $\op\colon \Sigma_*\to\Sigma_*$ defined as
follows.  For $k\ge0$, define $r_k\in\Sigma_k$ by $r_k(j)=k+1-j$,
so $r_k$ reverses order.  We then define
$$\op\colon \Sigma_k\to\Sigma_k
$$
by $\op(\sigma)=r_k\circ\sigma$.
\enddefinition

We leave to the reader the check that $\op$ defines a map of
operads.

\definition{Definition}\label{num15} Let $\A$ be an associative category.
The {\bf opposite} of $\A$, written $\A^\op$, is the associative
category given by the composite
$$\xymatrix{\Sigma_*\ar[r]^-{\op}&\Sigma_*\ar[r]^-{\A}&\P.}
$$
\enddefinition

\proclaim{Corollary}\label{num25} Right $\A$-module structures on a
small permutative category $\D$ determine and are determined by left
$\A^\op$-module structures on $\D$.
\endproclaim

\demo{Proof}  The automorphism $\Sigma_*\overto{\op}
\Sigma_*$ extends to an isomorphism $\Msig
\overto{\op}\rMsig$ for which the diagram
$$\xymatrix{\Sigma_*\ar[r]^-{\op}\ar[d]&\Sigma_*\ar[d]
\\ \Msig\ar[r]_-{\op}&\rMsig}
$$
commutes.  The extension is given by exactly the same formula:
using the elements $r_k\in\Sigma_k$ defined by $r_k(j)=k+1-j$, we
define $\op(\sigma)=r_k\circ\sigma$, and clearly if $\sigma(j)=k$,
then $\op(\sigma)(j)=1$.  The result now follows immediately.
\enddemo

\proclaim{Corollary}\label{num16} If $\R$ is bipermutative, so is $\R^\op$.
\endproclaim

\demo{Proof} The map ``op'' of operads extends to the map of operads
$$E(\op)\colon E\Sigma_*\to E\Sigma_*.
$$
\enddemo

\subsection{Bimodules}

The following is the explicit definition of a bimodule in the
context of permutative categories.

\definition{Definition}\label{num27}  Let $\A$ and $\B$ be associative
categories, and $\D$ a permutative category. We say that $\D$ is
an {\bf $\A$-$\B$ bimodule} if $\D$ is a left $\A$-module and a
right $\B$-module, the associativity diagram
$$\xymatrix{\A\times\D\times\B\ar[r]^-{\otimes\times1}
\ar[d]_-{1\times\otimes}
&\D\times\B\ar[d]^-{\otimes}
\\ \A\times\D\ar[r]_-{\otimes}&\D}
$$
commutes, and diagrams (e) and (f) from \refname{num9} commute in
all situations in which the maps are defined.
\enddefinition

For bimodule structures, the fundamental
parameter multicategory is as follows.

\definition{Definition}\label{num26}  The bimodule parameter multicategory
$\Bi$ has objects $A$, $B$ (the ``rings'', with $A$ acting on the
left and $B$ on the right) and $M$ (the ``module''). All sets of
$k$-maps are empty with the exception of those in which $M$
appears exactly once in the input and is the output, those
where all inputs and the output are $A$, and those where all inputs
and the output are $B$.  In the latter two cases the set of $k$-maps
is $\Sigma_k$.  In the case of
$\Bi_k(C_1,\ldots,C_k;D)$ with $C_j=D=M$ and all other entries
either $A$ or $B$, we set
$\Bi_k=\{\sigma\in\Sigma_k:\sigma(i)<\sigma(j)\Leftrightarrow
C_i=A\}$.  These are precisely the $\sigma$'s for which the list
$C_{\sigma^{-1}(1)},\ldots,C_{\sigma^{-1}(k)}$ is the list
$A^{\sigma(j)-1},M,B^{k-\sigma(j)}$.  In particular, $\sigma(j)$
must always be one plus the number of $A$'s occurring in the
input.  The $\Sigma_k$ action and the multiproduct are defined
exactly as for the operad $\Sigma_*$.
\enddefinition

Note in particular that
restriction to either of the single objects $A$ or $B$ determines
a multifunctor $\Sigma_*\to\Bi$.

\proclaim{Theorem}\label{num28} Let $\A$ and $\B$ be small associative categories.
Then an $\A$-$\B$ bimodule structure on a small permutative
category $\D$ determines and is determined by a multifunctor
$\Bi\to\P$ sending $M$ to $\D$, restricting on the single object
$A$ to the structure multifunctor $\Sigma_*\to\P$ for $\A$ and on
the single object $B$ to the structure multifunctor for $\B$.
\endproclaim

\demo{Proof} Given a bimodule structure on $\D$ and an element
$\sigma\in\Bi_k(C_1,\ldots,C_k;D)$, we need to define a functor
$T\sigma$, and we use the usual formula
$$T\sigma(c_1,\ldots,c_k)
=c_{\sigma^{-1}(1)}\otimes\cdots\otimes c_{\sigma^{-1}(k)}.
$$
The proof that this gives a multifunctor $\Bi\to\P$ now proceeds
in exactly the same way as in the proof of
\refname{num10}. Conversely, suppose
we are given a multifunctor $T\colon \Bi\to\P$ satisfying the
conditions in the theorem.  Restricting to pairs of objects
$(A,M)$ or $(B,M)$ gives us restriction multifunctors
$\Msig\to\Bi$ and $\rMsig\to\Bi$, and we immediately obtain a left
$\A$-module structure on $\D$ and a right $\B$-module structure on
$\D$.  The associativity diagram commutes because $\Bi_3(A,M,B;M)$
has only one element, and diagrams (e) and (f) commute exactly as
in the proof of \refname{num10}.  This concludes the proof.
\enddemo

\proclaim{Corollary}\label{num28.5} If $\D$ is an $\A$-$\B$ bimodule for
associative categories $\A$ and $\B$, then $K\D$ is a $K\A$-$K\B$
bimodule in symmetric spectra.
\endproclaim

In the case where $\A=\B$, we can collapse the parameter
multicategory further using the parameter multicategory in the
second example after \refname{num3}:

\definition{Definition}\label{num29} The parameter multicategory
$\bMsig$ has two objects, $A$ and $M$, and is a parameter
multicategory for modules, so there are no $k$-morphisms unless
$M$ is the output and appears exactly once in the input, or else
$A$ is the output and only $A$ appears in the input.  In these
cases the $k$-morphisms are $\Sigma_k$, with the multiproduct
defined as in $\Sigma_*$.
\enddefinition

To compare this multicategory with the previous one, we use
the following lemma:

\proclaim{Lemma}\label{num30}  Consider the diagram of multicategories
$$\xymatrix{\Sigma_*\vphantom{\Msig}\ar@<.5ex>[r]\ar@<-.5ex>[r]&\Bi\ar[r]
&\bMsig}
$$
where the two arrows on the left are the inclusions of the
endomorphism operads of the objects $A$ and $B$, and the arrow on
the right sends both $A$ and $B$ to $A$, and sends permutations in
$\Bi$ to corresponding ones in $\bMsig$.  This is a coequalizer
diagram of multicategories.
\endproclaim

\demo{Proof} The key point here is that each permutation in
$\bMsig_k(A^{j-1},M,A^{k-j};M)$ has exactly one preimage in $\Bi$.
Once we realize this, extending an equalizing multifunctor to
$\bMsig$ is simply a matter of sending all permutations to their
images under the multifunctor.
\enddemo

The characterization of $\A$-$\A$ bimodules in terms of a
parameter multicategory now follows immediately.

\proclaim{Corollary}\label{num31} If $\A$ is a small associative
category and $\D$ is a small permutative category, then an
$\A$-$\A$ bimodule structure on $\D$ determines and is determined
by a multifunctor $\bMsig\to\P$ sending $M$ to $\D$ and
restricting on $A$ to the associative category structure
multifunctor $\Sigma_*\to\P$ for $\A$.
\endproclaim

The analog of \refname{num28.5} now follows as well.

If one or both of $\A$ and $\B$ are bipermutative, one can also
describe $\A$-$\B$ bimodules with this extra structure in terms of
parameter multicategories.  We leave this to the interested
reader.

We can also ask for an analogous characterization of $\A$-$\A$
bimodules as in \refname{num31} in the case where $\A$ is
bipermutative.  The answer is NOT to apply $E$ to all the
multicategories in the diagram in \refname{num30}.  (This illustrates
the fact that $E$ does not preserve coequalizers).  Instead, we get
a multicategory described as follows.

\definition{Definition}\label{num36} The multicategory $\bEMsig$ is
a parameter multicategory for modules, so has objects $A$ and $M$,
with the $k$-morphisms empty except in the cases where $M$ appears
exactly once in the input and is the output, or else all inputs
and the output are $A$.  We set $\bEMsig_k(A^k;A)=E\Sigma_k$. The
objects of $\bEMsig_k(A^{j-1},M,A^{k-j};M)$ are the elements of
$\Sigma_k$, but the objects are not all isomorphic.  Instead, we
look at the equivalence relation on $\Sigma_k$ in which
$\sigma\sim\sigma'$ if and only if $\sigma(j)=\sigma'(j)$ and
$\sigma$ and $\sigma'$ are in the same coset of the left action of
$\Sigma_{\sigma(j)-1}\times\Sigma_{k-\sigma(j)}$ on $\Sigma_k$.
Equivalently, we could say that $\sigma\sim\sigma'$ means that
$\sigma(i)<\sigma(j)\Leftrightarrow\sigma'(i)<\sigma'(j)$ whenever
$1\le i\le k$. There is exactly one morphism from $\sigma$ to
$\sigma'$ when $\sigma$ and $\sigma'$ are equivalent and no
morphisms when they are not equivalent. We leave it to the reader
to check that the same formula for the multiproduct in $\Sigma_*$
extends to give multicategory structure on $\bEMsig$.
\enddefinition

\proclaim{Lemma}\label{num37} Consider the diagram of multicategories
$$\xymatrix{E\Sigma_*\vphantom{\Msig}\ar@<.5ex>[r]\ar@<-.5ex>[r]&E\Bi\ar[r]
&\bEMsig}
$$
where the two arrows on the left are the inclusions of the
endomorphism operads of the objects $A$ and $B$, and the arrow on
the right sends both $A$ and $B$ to $A$, and sends permutations to
themselves.  This is a coequalizer diagram of
multicategories.
\endproclaim

\demo{Proof} Given \refname{num30}, the only issue is the
morphisms.  However, the definition of the morphisms in $\bEMsig$
is precisely the requirement that two $k$-morphisms are isomorphic
in $\bEMsig$ if and only if they come from isomorphic
$k$-morphisms in $E\Bi$.  The result follows.
\enddemo

\proclaim{Corollary}\label{num38}  Let $\R$ be a small bipermutative category.
Then $\R$-$\R$ bimodule structures on a small permutative category
$\D$ determine and are determined by multifunctors $\bEMsig\to\P$
sending $A$ to $\R$ and $M$ to $\D$, and which restrict on $A$ to
the bipermutative structure map $E\Sigma_*\to\P$ for $\R$.
Consequently, the $K$-theory spectrum $K\D$ is equivalent to a
bimodule over a strictly commutative ring spectrum equivalent to
$K\R$.
\endproclaim

This still leaves the question of what sort of bimodule structure
is parameterized by $E\bMsig$.  The relevant definition is as
follows.

\definition{Definition}\label{num39} Let $\R$ be a bipermutative category.
The structure of a {\bf symmetric bimodule} over
$\R$ on a permutative category $\D$ consists of an $\R$-$\R$ bimodule
structure together with a natural isomorphism
$$\gamma\colon r\otimes d\cong d\otimes r
$$
for $r$ an object of $\R$ and $d$ an object of $\D$.  The
isomorphism $\gamma$ must be compatible with the multiplicative
symmetry isomorphism $\gamma^{\otimes}$ for $\R$, in the
sense that all possible diagrams of the form given in part 3 of
\refname{num7} must commute (with the $\oplus$'s
replaced with $\otimes$'s).  We also require diagram
(e$'$) given in \refname{num11} to commute.
\enddefinition

\proclaim{Theorem}\label{num40}  Let $\R$ be a small bipermutative category and
$\D$ a small permutative category.  Then symmetric bimodule
structures for $\D$ over $\R$ determine and are determined by
multifunctors $E\bMsig\to\P$ sending $M$ to $\D$ and restricting
on $A$ to the structure map $E\Sigma_*\to\P$ for $\R$ as a
bipermutative category.  Consequently, the $K$-theory spectrum
$K\D$ is equivalent to a module over a strictly commutative ring
spectrum equivalent to $K\R$.
\endproclaim

The proof is the same as the proof of \refname{num12} with
$\bMsig$ in place of $\Sigma_{*}$.

%\demo{Proof}  First suppose that we are given a multifunctor
%$E\bMsig\to\P$ as specified.  Applying $E$ to the diagram in
%\refname{num30} and using \refname{num37}, we get a canonical
%multifunctor
%$\bEMsig\to E\bMsig$, and therefore by \refname{num38}, we get an
%$\R$-$\R$ bimodule structure on $\D$.  The symmetry
%isomorphism $\gamma\colon  r\otimes d\cong d\otimes r$ is the image of
%the isomorphism between the two objects of $E\bMsig_2(A,M;M)$
%which are the two elements of $\Sigma_2$.  The rest of the proof
%that this makes $\D$ into a symmetric $\R$-bimodule is exactly as
%in the proof of \refname{num12}.
%
%Next, suppose $\D$ is a symmetric $\R$-bimodule, and we wish to
%construct the multifunctor $E\bMsig\to\P$.  The map on
%$k$-morphism objects is as in the proof of \refname{num12}, leaving us
%with the need for the map on the morphisms in $E\Sigma_k$ for each
%$k$.  However, our requirements for a symmetric bimodule are
%designed to allow the rest of the proof of \refname{num12} to go
%through, as the reader may easily check.  This completes the
%proof.
%\enddemo

\subsection{Algebras}

We turn our attention next to algebras.
The definition of a central algebra over a
bipermutative category depends on the notion of a central map from
a bipermutative category to an associative category, which we
define first.

\definition{Definition}\label{num42}  Let $\R$ be a bipermutative category
and $\A$ an associative category.  A {\bf central map} from $\R$
to $\A$ is a lax map $\phi\colon \R\to\A$ (i.e.,
$(\phi,\lambda)\in\obj{\P_1(\R;\A)}$)
and a natural isomorphism
$\gamma\colon \phi(r)\otimes a\cong a\otimes\phi(r)$ for $r$ an
object of $\R$ and $a$  an object of $\A$, satisfying the following
conditions:
\roster
\item $\phi$ preserves the tensor product in the sense that the diagram
$$\xymatrix{
\R\times\R\ar[r]^-{\phi\times\phi}\ar[d]_-{\otimes}
&\A\times\A\ar[d]^-{\otimes}
\\ \R\ar[r]_-{\phi}&\A}
$$
commutes strictly and $\phi(1)=1$.
\item The lax structure map
$\lambda$ preserves the distributivity maps in the sense that the
diagram
$$\xymatrix{
(\phi r_1\otimes\phi r_2)\oplus(\phi r_1\otimes\phi r_3)
\ar[r]^-{d_r}\ar[d]_-{=}
&\phi r_1\otimes(\phi r_2\oplus\phi r_3)\ar[d]^-{1\otimes\lambda}
\\ \phi(r_1\otimes r_2)\oplus\phi(r_1\otimes r_3)
\ar[d]_-{\lambda}
&\phi r_1\otimes\phi(r_2\oplus r_3)\ar[d]^-{=}
\\ \phi[(r_1\otimes r_2)\oplus(r_1\otimes r_3)]
\ar[r]_-{\phi(d_r)}
&\phi(r_1\otimes(r_2\oplus r_3))}
$$
and a similar diagram involving $d_l$ commute.
\item $\gamma$ must be consistent with the symmetry
isomorphism $\gamma^{\otimes}$ in $\R$ in the sense for all
objects $r_1$, $r_2$ of $\R$, the diagram
$$\xymatrix{
\phi(r_1)\otimes\phi(r_2)\ar[r]^-{\gamma}\ar[d]_-{=}
&\phi(r_2)\otimes\phi(r_1)\ar[d]^-{=}
\\ \phi(r_1\otimes r_2)\ar[r]_-{\phi(\gamma^{\otimes})}
&\phi(r_2\otimes r_1)}
$$
commutes.
\item $\gamma$ satisfies all
instances of the diagrams in part (3) of \refname{num7}, and
diagram (e$'$) of \refname{num11}.
\endroster
An {\bf $\R$-algebra} structure on $\A$ consists of a central map from
$\R$ to $\A$.
\enddefinition

\definition{Definition}\label{num41} Let $\AA$ be the
multicategory with two objects, $R$ (the ground ring) and $A$ (the
algebra).  The category $\AA_k(B_1,\ldots,B_k;C)$ is empty if $C=R$
and one or more of the $B_{j}$'s are $A$.  Otherwise,
$\AA_k(B_1,\ldots,B_k;C)$ has $\Sigma_{k}$ as its set of objects, and
has morphisms as follows.
Let $S=\{j:B_j=A\}$ and consider the equivalence relation
on the elements of $\Sigma_k$ where $\sigma\sim\sigma'$ means
that for all $i$ and $j$ in $S$,
$\sigma(i)<\sigma(j)\Leftrightarrow\sigma'(i)<\sigma'(j)$.  We have
precisely one morphism from $\sigma$ to $\sigma'$ when
$\sigma\sim\sigma'$, and no morphisms between inequivalent elements.
\enddefinition

In the previous definition, if we restrict our attention to the
object $R$, we get $E\Sigma_*$, while if we restrict our attention
to the object $A$, we get $\Sigma_*$. We wish to show that
$\R$-algebra structures on a small associative category $\A$
correspond to multifunctors from $\AA$ to $\P$ extending the
structure multifunctors for both $\R$ and $\A$.  To do this, we
need the following combinatorial lemma about permutations.

\proclaim{Lemma}\label{num43} Suppose $T\subset\u{k}=\{1,\ldots,k\}$ and
that $\rho\in\Sigma_k$ is order-preserving on $T$ in the sense
that if $i$ and $j$ are elements of $T$ with $i<j$, then
$\rho(i)<\rho(j)$.  Then $\rho$ can be written as a product of
transpositions of consecutive integers in $\u{k}$, say
$\rho=t_1\cdots t_m$, in such a way that for $1\le n\le m$, $t_n$
does not transpose two elements of $t_{n+1}\cdots t_mT$.
\endproclaim

\demo{Proof} Let the elements of $T$ be written in order as
$\{a_1,\ldots,a_q\}$.  First, we use transpositions of the
required form to map $T$ to $\{1,\ldots,q\}$; we do this by first
transposing $a_1$ with its predecessors, in order, and then
repeating the process with $a_2$ through $a_q$. Then use
transpositions of adjacent elements of $\{q+1,\ldots,k\}$ to
rearrange this set in the same order that $\rho$ rearranges
$\u{k}\setminus T$. Finally, start with $q$ and transpose it with
its successors, in order, until it reaches $\rho(a_q)$, and repeat
the process with $q-1$ back through 1. The result is $\rho$, with
the transpositions involved having the required property.
\enddemo

\proclaim{Theorem}\label{num44} Let $\R$ be a small bipermutative
category and $\A$ a small associative category.  Then $\R$-algebra
structures on $\A$ determine and are determined by multifunctors
from $\AA$ to $\P$ restricting on the object $R$ to the structure
multifunctor for $\R$ as a bipermutative category and on the
object $A$ to the structure multifunctor for $\A$ as an
associative category.  Consequently, $K\A$ is equivalent to a
central algebra over a strictly commutative ring spectrum
equivalent to $K\R$.
\endproclaim

\demo{Proof} Suppose we are given a multifunctor from
$\AA$ restricting as required.  Then we obtain a functor
$\phi\colon \R\to\A$ as the image of the unique element $1_1$ of
$\AA_1(R;A)$; we claim that this functor is a central map.  First,
we have the formula
$\multprod(1_1;1_2)=\multprod(1_2;1_1,1_1)=1_2$ in $\AA$, which we
can express by saying that the diagram in $\AA$
$$\xymatrix{
(R,R)\ar[r]^-{(1_1,1_1)}\ar[d]_-{1_2}
&(A,A)\ar[d]^-{1_2}
\\ R\ar[r]_-{1_1}&A}
$$
commutes, and consequently its image in $\P$
$$\xymatrix{
\R\times\R\ar[r]^-{\phi\times\phi}\ar[d]_-{\otimes}
&\A\times\A\ar[d]^-{\otimes}
\\ \R\ar[r]_-{\phi}&\A}
$$
commutes as well.  A similar argument shows that $\phi(1)=1$.
Since the commutativity of this diagram in $\P$ also requires that
the distributivity maps coincide, we get the diagrams showing that
$\lambda$ preserves the distributivity maps.  The natural
isomorphism $\gamma\colon \phi(r)\otimes a\cong a\otimes\phi(r)$ is the
image of the isomorphism between the two elements of
$\AA_2(R,A;A)=\Sigma_2$.  Because the diagram
$$\xymatrix{
(R,R)\ar[r]^-{(1_1,1_1)}\ar[d]
&(R,A)\ar[d]
\\ R\ar[r]_-{1_1}&A}
$$
in $\AA$ commutes with the downward arrows being either of the two
elements of $\Sigma_2$, the isomorphism between the two possible
elements on the left gets taken by $\phi$ to the isomorphism
between the two possible elements on the right, i.e.,
$\gamma=\phi(\gamma^\otimes)$, as required.  Further, diagram
(e$'$) of
\refname{num11} is satisfied because $\gamma$ is a morphism in
$\P_2(\R,\A;\A)$.  We therefore get a central map $\phi\colon \R\to\A$
given a multifunctor $\AA\to\P$ restricting to the structure
multifunctors of $\R$ and $\A$ on the objects $R$ and $A$,
respectively.

Now suppose we are given a central map $\phi\colon \R\to\A$; we must
show that this extends uniquely to a multifunctor $\AA\to\P$ by
requiring the multifunctor to restrict to the structure
multifunctors for $\R$ and $\A$ and also by requiring the single
element of $\AA_1(R;A)$ to map to $\phi$.
The functor on $\AA_k(B_1,\ldots,B_k;C)$ is
already determined when $C=R$ or when $C=A$ and all the $B_{j}$'s are
$A$.  In the other cases, set
$S=\{i:B_i=A\}$ as in the definition.  It remains to
determine the images of the categories $\AA_k(B_1,\ldots,B_k;A)$
with $S\ne\emptyset$ and $S\ne\{1,\dotsc,k\}$.  By equivariance, it
suffices to consider the special case $S=\{1,\ldots,q\}$ for
$q<k$.  The objects are the elements of $\Sigma_k$, and it is
clear that the image of $1_k$ is the composite
$$\xymatrix@C+6pt{
\A^q\times\R^{k-q}\ar[r]^-{1\times\phi^{k-q}}
&\A^k\ar[r]^-{\otimes}&\A},
$$
and the images of the rest of the objects are determined by
equivariance.  We must also determine the images of the
isomorphisms in $\AA_k(B_1,\ldots,B_k;A)$. For this, note that
when $\sigma \sim \sigma'$ as in the definition,
$\sigma'\sigma^{-1}$ is order-preserving on $\sigma S$, so by
\refname{num43}, can be written as a product of transpositions of
adjacent integers which are not both elements of $\sigma S$.  Now
the image of a typical $k$-tuple $(b_1,\ldots,b_k)$ under the
element $\sigma$ is $b_{\sigma^{-1}(1)}\otimes\cdots\otimes
b_{\sigma^{-1}(k)}$, and we need to produce an isomorphism between
this and the image under $\sigma'$.  Write $\sigma'\sigma^{-1}$ as
$t_1\cdots t_m$, where $t_j$ is a transposition of adjacent
integers not both in $t_{j+1}\cdots t_m\sigma S$, and say $t_m$
transposes $i$ and $i+1$. Then the term $b_{\sigma^{-1}(i)}\otimes
b_{\sigma^{-1}(i+1)}$ appears as part of the image under $\sigma$,
and since $\sigma^{-1}(i)$ and $\sigma^{-1}(i+1)$ are not both
elements of $S$, the two $b$'s are not both objects of $\A$, so
they can be transposed using $\gamma$.  We get an isomorphism
between a tensor product of elements of the form
$$b_{\sigma^{-1}(i)}=b_{\sigma'{}^{-1}\sigma'\sigma^{-1}(i)}
=b_{\sigma'{}^{-1}t_1\cdots t_m(i)}
$$
and elements of the form
$$b_{\sigma'{}^{-1}t_1\cdots t_{m-1}(i)}.
$$
By iterating the process $m$ times, we get an isomorphism between
the image under $\sigma$ and the image under $\sigma'$.  The
isomorphism is uniquely determined by $\sigma'\sigma^{-1}$ and not
its presentation, because the $\gamma$'s satisfy the relations
among transpositions in $\Sigma_k$.  This completes the proof.
\enddemo

In the special case where $\A$ is also a bipermutative category
and the symmetry isomorphism is given by the isomorphism
already present in $\A$, we can give a somewhat simpler
description.

\definition{Definition}\label{num45} Let $\R$ and $\A$ be bipermutative
categories.  A {\bf map} of bipermutative categories
$\phi\colon \R\to\A$ is a lax map that preserves the tensor product,
distributivity maps, and multiplicative unit in the same sense
that a central map does, and for which also
$\phi(\gamma^{\otimes}_\R)=\gamma^{\otimes}_\A$.
\enddefinition

The corresponding definition in terms of a parameter multicategory
is as follows.

\definition{Definition}\label{num46} The multicategory $\AE$ is a parameter
multicategory for algebras, so by \refname{num4} has two objects,
$A$ and $R$, and with $\AE_k(B_1,\ldots,B_k;C)=\emptyset$ if
$S\ne\emptyset$ and $C=R$, where $S=\{i:B_i=A\}$.  Otherwise, we
set $\AE_k(B_1,\ldots,B_k;C)=E\Sigma_k$, so this is an example of
the sort discussed as the third example following \refname{num4}.
\enddefinition

The proof of the following theorem can now be safely left to the
reader.

\proclaim{Theorem}\label{num47} Let $\R$ and $\A$ be small bipermutative
categories.  Then a map of bipermutative categories $\phi\colon
\R\to\A$ determines and is determined by a multifunctor $\AA\to\P$
which restricts on the object $R$ to the structure multifunctor
for $\R$ and on the object $A$ to the structure multifunctor for
$\A$.  Consequently, $K\phi$ is equivalent to a map of strictly
commutative ring spectra.
\endproclaim

\section{Free Permutative Categories}\label{freepcs}

This section is devoted to the construction of additional examples
of both associative and bipermutative categories via the ``free
permutative category'' construction.  This associates to any
small category $\C$ a small permutative category $\PP\C$ as follows.  Let
$E\Sigma_k$ be the translation category of $\Sigma_k$.  Then we
define
$$\PP\C=\coprod_{k\ge0}E\Sigma_k\times_{\Sigma_k}\C^k.
$$
The objects of $\PP\C$ are the elements of the free monoid on the
objects of $\C$, with 0 given by the empty string and the direct
sum given by concatenation, which is the monoid operation.  The
symmetry isomorphism arises from the isomorphism in $E\Sigma_2$
between the two elements of $\Sigma_2$.  Dunn \cite{3} apparently
first observed that $\PP$ defines a monad in $\Cat$ whose algebras
are precisely the small permutative categories. The resulting
morphisms are called the {\bf strict} morphisms and are even more
restrictive than the strong morphisms.  In fact, they are too
restrictive to form a multicategory.

The following theorem shows how additional structure on
$\C$ gives rise to additional structure on $\PP\C$.

\proclaim{Theorem}\label{num13.5} Let $\C$ be a small strict monoidal category
(i.e., one equipped with a strictly associative and unital
``tensor product'' operation).  Then $\PP\C$ can be made into an
associative category.  If $\C$ is permutative, then $\PP\C$
becomes a bipermutative category.
\endproclaim

\demo{Proof} There are actually uncountably many different ways of
constructing such structure, depending on one's choice of what we
call a {\bf priority order}.  Let $\u m$ denote the set
$\{1,\ldots,m\}$ for positive integers $m$.  Then a priority order
is a choice of bijection $\omega_{m,n}\colon \u{mn}\to\u m\times\u n$
for each $m$ and $n$ that is coherent in the sense that all
diagrams of the form
$$\xymatrix{
\u{mnp}\ar[r]^-{\omega_{mn,p}}\ar[d]_-{\omega_{m,np}}
&\u{mn}\times\u p\ar[d]^-{\omega_{m,n}\times1}
\\ \u m\times\u{np}\ar[r]_-{1\times\omega_{n,p}}
&\u m\times\u n\times\u p}
$$
commute.  By ordering $\u m\times\u n$ using lexicographic order
and taking the inverse of the resulting bijection, we get a
priority order, as we do using reverse lexicographic order, but
there are uncountably many other choices as well.  For example, we
can use lexicographic order to define a bijection $\u
m\to\u{2^{\nu(m)}}\times\u{\hat{m}}$, where $\hat{m}$ is odd, and
then for any $m$ and $n$, use the inverse of the bijection
$$\xymatrix@R=8pt{
\u m\times\u n\ar[r]&
\u{2^{\nu(m)}}\times\u{\hat{m}}\times\u{2^{\nu(n)}}\times\u{\hat{n}}
\\ \relax\ar[r]^-{1\times\tau\times1}&
\u{2^{\nu(m)}}\times\u{2^{\nu(n)}}\times\u{\hat{m}}\times\u{\hat{n}}
\ar[r]&\u{2^{\nu(m)}2^{\nu(n)}\hat{m}\hat{n}}=\u{mn},}
$$
where the unlabelled arrows are given by lexicographic order or
its inverse.  We can use the same sort of trick for any set of
primes, not just 2, to get uncountably many additional priority
orders.  In any case, pick one, and call it $\omega$.  Let
$\omega_1$ and $\omega_2$ denote $\omega$ followed by projection
onto the first or second factor, respectively. Then we define an
associative structure on $\PP\C$ as follows. Write a typical
object $(a_1,\ldots,a_m)$ of $\PP\C$ as $\oplus_{i=1}^m(a_i)$, and
write the monoidal operation in $\C$ as $\otimes$. Then we define
the tensor product on $\PP\C$ by the formula
$$\bigoplus_{i=1}^m(a_i)\otimes\bigoplus_{j=1}^n(b_j)
:=\bigoplus_{k=1}^{mn}(a_{\omega_1(k)}\otimes b_{\omega_2(k)}).
$$
In the case where $\C$ is permutative, we can then use the
symmetry isomorphism in $\C$ to map this to
$$\bigoplus_{k=1}^{mn}(b_{\omega_2(k)}\otimes a_{\omega_1(k)}),
$$
and then shuffle inside of $\PP\C$ to map this to
$$\bigoplus_{k=1}^{mn}(b_{\omega_1(k)}\otimes a_{\omega_2(k)}),
$$
defining the multiplicative symmetry isomorphism necessary
for a bipermutative category.  The reader can check that one needs
only the associativity condition on a priority order to show that
these definitions satisfy the requirements for an associative or a
bipermutative category, respectively.
\enddemo

An example of particular
importance of this form is the free permutative category $\PP(*)$
on a one point category, which becomes a bipermutative category
via this construction.  The reader should be aware, however, that
modules over $\PP(*)$ depend strongly on the priority order
chosen.  We leave as an exercise to the reader that if we use
lexicographic order, then any permutative category is a left
module over $\PP(*)$, while if we use reverse lexicographic order,
every permutative category is a right module over $\PP(*)$.  Of
course, the two orders give opposite bipermutative structures on
$\PP(*)$, so the duality is to be expected.  Other choices of
priority order seem to give far fewer modules over $\PP(*)$.

\section{Model Categories of Rings, Modules, and Algebras in Symmetric
Spectra}
\label{secpmssmodel}

In this section we prove \refname{num5}.  Fix a small
multicategory $\M$ enriched over simplicial sets, and let $\OB$
denote its set of objects.  Let $\S^{\OB}$ denote the category
obtained as the product of copies of the category $\S$ of
symmetric spectra indexed on the set $\OB$.  As a product
category,
$\S^{\OB}$ inherits a simplicial closed model structure for each
simplicial closed model structure on $\S$, precisely, one with its
fibrations, cofibrations, and weak equivalences formed objectwise
(i.e., coordinatewise). Our goal is to prove that the category
$\S^{\M}$ of simplicial multifunctors from $\M$ to $\S$ has a
simplicial closed model structure with the fibrations and weak
equivalences the maps that are fibrations and weak equivalences
respectively in $\S^{\OB}$ for the positive stable model structure
on $\S$. Throughout this section, we use the terminology {\bf stable
equivalence}, {\bf positive stable fibration}, and {\bf acyclic
positive stable fibration} in $\S^{\M}$ to indicate those maps in
$\S^{\M}$ whose underlying maps in $\S^{\OB}$ are weak
equivalences, fibrations, and acyclic fibrations in the positive
stable model structure.

The first step is to show that the category $\S^{\M}$ has limits
and colimits.  For this, it is convenient to observe that
$\S^{\M}$ is the category of algebras over a monad $\MM$ on
$\S^{\OB}$.

\definition{Definition}\label{defmm}
For $b\in \OB$, and $T$ in $\S^{\OB}$, let
$$ (\MM T)_{b}=
\bigvee_{n\geq 0}\left(\bigvee_{a_{1},\dotsc,a_{n}\in \OB}
\M(a_{1},\dotsc,a_{n};b)_{+}\sma
(T_{a_{1}}\sma \dotsb \sma T_{a_{n}})\right)/\Sigma_{n}, $$ let
$\eta \colon T\to \MM T$ be the map
$$ T_{b}\iso \{\id_{b}\}_{+}\sma T_{b}\to  \M(b;b)_{+}\sma T_{b}
\to (\MM T)_{b}, $$
and $\mu \colon \MM \MM T\to \MM T$ the map induced by the
multiproduct of $\MM$.
\enddefinition

The proof of the following theorem in the special case of operads
\cite{13} easily generalizes to multicategories.

\proclaim{Theorem}\label{monad}
$\MM$ is a simplicial monad on the category $\S^{\OB}$.  An
$\MM$-algebra structure on an object of $\S^{\OB}$ is equivalent
to an $\M$-multifunctor structure, and the simplicial category of
$\MM$-algebras is isomorphic to $\S^{\M}$.  \endproclaim

\proclaim{Corollary}\label{moadj}
$\MM$, viewed as a functor $\S^{\OB}\to \S^{\M}$, is left adjoint
to the forgetful functor $\S^{\M}\to \S^{\OB}$.
\endproclaim

\proclaim{Corollary}\label{complete}
The category $\S^{\M}$ is complete and cocomplete (has all small
limits and colimits), and is tensored and cotensored over
simplicial sets.
\endproclaim

\demo{Proof} As a category of algebras over a monad on a complete
category, $\S^{\M}$ is complete, with limits and cotensors formed
in $\S^{\OB}$. Since $\MM$ preserves reflexive coequalizers (by
the argument of
\cite{5} Proposition~II.7.2), $\S^{\M}$ is cocomplete with
reflexive coequalizers created in $\S^{\OB}$ by \cite{5}
Proposition~II.7.4. General colimits are formed by rewriting the
colimit as a reflexive coequalizer, and the tensor of an object
$A$ of $\S^{\M}$ and a simplicial set $X$ is formed as a
(reflexive) coequalizer of the form
$$ \xymatrix{
\MM((\MM A)\sma X_{+})\ar@<.5ex>[r]\ar@<-.5ex>[r]
&\MM(A\sma X_{+})\ar[r]
&A \otimes X.
} $$
\enddemo

In order to prove the required factorization and lifting
properties, we need to review briefly the positive stable model
structure on $\S$. Recall that in any category $\ACat$ with small
colimits, for any set $I$ of maps, a relative $I$-complex
(\cite{12} Definition~5.4) is a map $X\to Y$ in $\ACat$ where
$Y=\Colim X_{k}$, with $X_{0}=X$, and $X_{k+1}$ is formed from
$X_{k}$ as a pushout of a coproduct of maps in $I$.  In this
terminology, a map of symmetric spectra is a cofibration in the
positive stable model structure if and only if it is a retract of
a relative $\IP$-complex, where
$$ \IP =\{F_{m}\partial \Delta[n]_{+}\to F_{m}\Delta[n]_+\mid m>0,n\geq
0  \},
$$
and $F_{m}$ is the functor from simplicial sets to symmetric spectra
left adjoint to the $m$-th space functor. A map is an acyclic
cofibration if and only if it is a retract of a relative $\KP$-complex
for a certain set of maps $\KP$ (q.v. \cite{7} Definition~3.4.9 and
\cite{12} Section~14).  A complete description of the maps in $\KP$
is not difficult but would require an unnecessary digression; all we
need to know about the maps is that the domain and codomain are small,
meaning that the sets of maps out of them commute with sequential
colimits.

For $a\in \OB$, let $\IO{a}$ denote the functor $\S\to \S^{\OB}$
that is left adjoint to the projection functor $\pi_{a}\colon
\S^{\OB}\to
\S$.  For a symmetric spectrum $T$, the object $\IO{a}T$ of $\S^{\OB}$
satisfies
$$ (\IO{a} T)_{b}=
\cases
T&b=a\\
*&b\neq a.
\endcases $$
The positive stable model structure on $\S^{\OB}$ then has a
similar description of its cofibrations and acyclic cofibrations:
Let
$$ \align{
\IO*\IP&=\{\IO{a}f\mid f\in \IP,a\in \OB\}\cr
\IO*\KP&=\{\IO{a}f\mid f\in \KP,a\in \OB\}.
}\endalign
$$
A map in $\S^{\OB}$ is cofibration if and only if it
is the retract of a relative $\IO*\IP$-complex and is an acyclic
cofibration if and only if it is a retract of a relative
$\IO*\KP$-complex. Let
$$ \align{
\MI&=\MM\IO*\IP=\{\MM\IO{a}f\mid f\in \IP,a\in \OB\}=\{\MM f\mid f\in \IO*\IP\}\cr
\MK&=\MM\IO*\KP=\{\MM\IO{a}f\mid f\in \KP,a\in \OB\}=\{\MM f\mid f\in \IO*\KP\}.
}\endalign
$$
The adjunction of \refname{moadj} and the lifting
properties in $\S^{\OB}$ then imply the following.

\proclaim{Proposition}\label{rlp}
A map in $\S^{\M}$ is an acyclic positive stable fibration if and
only if it has the right lifting property with respect to $\MI$,
if and only if it has the right lifting property with respect to
retracts of relative $\MI$-complexes. It is a positive stable
fibration if and only if it has the right lifting property with
respect to $\MK$, if and only if it has the right lifting property
with respect to retracts of relative $\MK$-complexes.
\endproclaim

Because the domains and codomains of the maps in $\IP$ and $\KP$
are small in symmetric spectra, the domains and codomains of the
maps in $\MI$ and $\MK$ are small in $\S^{\M}$.  The Quillen small
object argument then gives the following.

\proclaim{Proposition}\label{factor}
A map in $\S^{\M}$ can be factored as a relative $\MI$-complex
followed by an acyclic positive stable fibration or as a relative
$\MK$-complex followed by a positive stable fibration.
\endproclaim

The proof of the following lemma is complicated but similar to the
analogous lemma in the case of commutative ring symmetric spectra.
%(\cite{M} Theorem~7.5).
Since we need some specifics of the argument in the next section, we
provide the proof at the end of that section.

\proclaim{Lemma}\label{mklemma}
A relative $\MK$-complex is a stable equivalence.
\endproclaim

The usual lifting and retract argument then gives the following.

\proclaim{Proposition}\label{llp}
A map in $\S^{\M}$ has the left lifting property with respect to
the acyclic positive stable fibrations if and only if it is a
retract of a relative $\MI$-complex. A map in $\S^{\M}$ has the
left lifting property with respect to the positive stable
fibrations if and only if it is a retract of a relative
$\MK$-complex.
\endproclaim

We have now collected all the facts we need to prove
\refname{num5}.

\demo{Proof of \refname{num5}}
We have shown (in \refname{complete}) that $\S^{\M}$ has all
finite limits and colimits.  It is clear by their definition that
weak equivalences (the stable equivalences) are closed under
retracts and have the two-out-of-three property.  Also clear from
the definition is that the fibrations (the positive stable
fibrations) are closed under retracts, and if we define the
cofibrations in terms of the left lifting property, then it is
clear that these are closed under retracts. The lifting properties
follows from \refname{rlp} and \refname{llp}, and the
factorization properties follow from \refname{factor}.  Thus, all
that remain is SM7.

We need to show that when $i\colon T\to U$ is a cofibration and
$p\colon X\to Y$ is a fibration, the map of simplicial sets
$$ \S^{\M}(U,X)\longrightarrow\S^{\M}(U,Y)\times_{\S^{\M}(T,Y)}\S^{\M}(T,X) $$
is a fibration, and a weak equivalence if either $i$ or $p$ is.
Using the characterization in \refname{llp} of cofibrations and
acyclic cofibrations as the maps that are retracts of relative
$\MI$- and $\MK$-complexes respectively, this easily reduces to
the case when $i$ is a map in $\MI$ or a map in $\MK$.  Using the
adjunction of
\refname{moadj}, this reduces to SM7 in $\S^{\OB}$, which reduces
to SM7 in $\S$, proved in \cite{7}.
\enddemo

\section{Multifunctors and Quillen Adjunctions}
\label{secpmssfunc}

In this section we prove \refname{num5.5}.

Let $f\colon \M\to \Mp$ be a simplicial multifunctor between small
multicategories enriched over simplicial sets. Let $\OB$ denote
the set of objects of $\M$ and $\OBp$ the set of objects of $\Mp$.
The multifunctor $f$ in particular induces a projection functor
$\pi_{f}
\colon \S^{\OBp}\to \S^{\OB}$. Let $\IO{f}\colon \S^{\OB}\to\S^{\OBp}$  be
the left adjoint: For $T$ an object in $\S^{\OB}$ and $b$ in
$\OBp$,
$$ (\IO{f}T)_{b}= \bigvee_{a\in f^{-1}(b)} T_{a}. $$
The multifunctor $f$ induces a natural transformation
$$ \IO{f}\MM\to \MMp\IO{f}, $$
where $\MMp$ is the monad on $\S^{\OBp}$ from \refname{defmm}. For
an object $A$ of $\S^{\M}$, we use this natural transformation and
the structure map $\MM A\to A$ to construct $f_{*}A$ in $\S^{\Mp}$
by the (reflexive) coequalizer diagram
$$ \xymatrix{
\MM'\IO{f}\MM A \ar@<-.5ex>[r]\ar@<.5ex>[r]
&\MM'\IO{f}A \ar[r]
&f_{*}A.
} $$ Unwinding the universal property and the adjunctions, we
obtain the following result.

\proclaim{Proposition}\label{consfls}
$f_{*}\colon \S^{\M}\to \S^{\Mp}$ is left adjoint to the pullback
functor $f^{*}\colon \S^{\Mp}\to \S^{\M}$.
\endproclaim

Since the functor $f^{*}$ clearly preserves weak equivalences and
fibrations, the first statement of \refname{num5.5} is an immediate
consequence of the previous proposition. For the rest of
\refname{num5.5}, we need the full definition of weak
equivalence.  We begin by reviewing the definition from \cite{4}
of a weak equivalence of categories enriched over simplicial sets,
and for this, we need to recall the category of components.  When
$\ACat$ is a category enriched over simplicial sets, the sets of
components $\pi_{0}\ACat(x,y)$ for objects $x,y$ have the
composition
$$ \pi_{0}\ACat(y,z)\times \pi_{0}\ACat(x,y)\to \pi_{0}\ACat(x,z) $$
induced by the composition in $\ACat$.  This composition and the
identity components make $\pi_{0}\ACat$ into a category, called
the category of components.  Recall that a simplicial functor
$f\colon \ACat\to
\ACat'$ is a weak equivalence when the induced functor
$\pi_{0}f$ is an equivalence of categories of components and for
all objects $x,y$ in $\ACat$, the map of simplicial sets
$\ACat(x,y)\to \ACat'(fx,fy)$ is a weak equivalence. In the
following definition, we understand the category of components of
a enriched multicategory to be the category of components of its
underlying enriched category.

\definition{Definition}\label{wkequiv}
A simplicial multifunctor $f\colon \M\to \Mp$ is a weak
equivalence when the induced functor $\pi_{0}f$ is an equivalence
of categories of components and for all $a_{1},\dotsc,a_{n},b$ in
$\OB$, the map of simplicial sets $\M(a_{1},\dotsc,a_{n};b)\to
\Mp(fa_{1},\dotsc,fa_{n};fb)$ is a weak equivalence.
\enddefinition

For the rest of the section, we assume that $f$ is a weak
equivalence.  We need to show that $(f_{*},f^{*})$ is a Quillen
equivalence.  The following lemma is the first step.

\proclaim{Lemma}\label{psel}
A map $\phi \colon T\to U$ is a stable equivalence in $\S^{\Mp}$
if and only if $f^{*}\phi$ is a stable equivalence in $\S^{\M}$.
\endproclaim

\demo{Proof}
By definition, $f^{*}\phi$ is a stable equivalence in $\S^{\M}$ if
and only if it is a stable equivalence in $\S^{\OB}$, i.e., if and
only if $\pi_{f}\phi$ is a stable equivalence.  Since $\phi$ is a
stable equivalence in $\S^{\Mp}$ if and only if it is a stable
equivalence in $\S^{\OBp}$, it follows that $f^{*}$ takes stable
equivalences in $\S^{\Mp}$ to stable equivalences in $\S^{\M}$.
Thus, it remains to show that $\phi$ is a stable equivalence when
$f^{*}\phi$ is.

Assume that $f^{*}\phi$ is a stable equivalence.  Then for any $a$
in $\OBp$ in the image of $f$, $\phi_{a}\colon T_{a}\to U_{a}$ is
a stable equivalence.  If $b$ is an arbitrary element of $\OBp$,
then the hypothesis that $f$ is a weak equivalence implies that we
can find an $a$ in the image of $f$ and an isomorphism from $a$ to
$b$ in the category of components of $\Mp$.  Choosing maps in
$\Mp(a,b)$ and $\Mp(b,a)$ in the components giving such an
isomorphism and its inverse, there are generalized simplicial
intervals connecting the composites with the appropriate identity
map (on $a$ and on $b$). Using the naturality of $\phi$, it
follows that $\phi_{b}$ is (levelwise) weakly equivalent to
$\phi_{a}$, and is therefore a positive stable equivalence.
\enddemo

We spend much of the rest of the section proving the following
theorem.

\proclaim{Theorem}\label{unitthm}
If $A$ is a cofibrant object of $\S^{\M}$, then the unit $A\to
f^{*}f_{*}A$ of the $(f_*,f^*)$ adjunction is a stable
equivalence.
\endproclaim

Assuming the previous theorem for the moment, we have all we need
to prove \refname{num5.5}.

\demo{Proof of \refname{num5.5}}
It remains to show that when $f$ is a weak equivalence, the
Quillen adjunction $(f_{*}, f^{*})$ is a Quillen equivalence.  Let
$A$ be a cofibrant object of $\S^{\M}$ and $B$ a fibrant object of
$\S^{\Mp}$; we need to show that a map $\phi \colon f_{*}A\to B$
is a stable equivalence if and only if the adjoint map $\psi
\colon A\to f^{*}B$ is a stable equivalence.  By \refname{psel},
we know that $\phi$ is a stable equivalence if and only if
$f^{*}\phi$ is a stable equivalence.  Since $\psi$ is the
composite
$$ A \longrightarrow f^{*}f_{*}A
\buildrel{f^{*}\phi}\over{\longrightarrow} f^{*}B, $$
\refname{unitthm} implies that $\psi$ is a stable equivalence
if and only if $f^{*}\phi$ is.  This concludes the proof.
\enddemo

We now move on to the proof of \refname{unitthm}.  The proof
requires an analysis of the pushouts in $\S^{\M}$ of the form
$B\amalg_{\MM
\IO{x}X}\MM\IO{x}Y$ for a map of symmetric spectra $X\to Y$ and a map
$\IO{x}X\to B$ in $\S^{\OB}$.  For this we need to set up two
constructions.  For the first, for each $x_{1},\dotsc,x_{k}$ in
$\OB$, construct $\U{x_{1},\dotsc,x_{k}}B$ as the coequalizer in
$\S^{\OB}$
$$ \align{
&\displaystyle\bigvee_{n\geq 0}
\left(\displaystyle\bigvee_{a_{1},\dotsc,a_{n}}
\M(a_{1},\dotsc,a_{n},x_{1},\dotsc,x_{k};-)_{+}\sma
(\MM B)_{a_{1},\dotsc,a_{n}}
\right)/\Sigma_{n}\\
{\xymatrix@C=20pt{\relax\ar@<.5ex>[r]\ar@<-.5ex>[r]&\relax}}
&\displaystyle
\bigvee_{n\geq 0} \left(\bigvee_{a_{1},\dotsc,a_{n}}
\M(a_{1},\dotsc,a_{n},x_{1},\dotsc,x_{k};-)_{+}\sma
B_{a_{1},\dotsc,a_{n}}
\right)/\Sigma_{n}\\
{\xymatrix@C=20pt{\relax\ar[r]&\relax}}
&\;\;\U{x_{1},\dotsc,x_{k}}B.
}\endalign
$$
where $B_{a_{1},\dots,a_{n}}$ is shorthand for
$B_{a_{1}}\sma\dotsb
\sma B_{a_{n}}$ and similarly for $\MM B$.  (One map is induced
by the action map $\MM B\to B$ and the other by the multiproduct.)
The purpose of introducing $\U{*}B$ is that for any $T$ in
$\S^{\OB}$, the underlying object in $\S^{\OB}$ of the coproduct
$B\amalg \MM T$ in $\S^{\M}$ is
$$ \bigvee_{k}\left(\bigvee_{x_{1},\dotsc,x_{k}}
\U{x_{1},\dotsc,x_{k}}B\sma
T_{x_{1}}\sma \dotsb \sma T_{x_{k}}\right)/\Sigma_{k}. $$ When
$x_{1}=\dotsb=x_{k}=x$ and $x$ is understood, we write $\U{k}B$
for $\U{x_{1},\dotsc,x_{k}}B$.

The second construction is defined for maps of symmetric spectra
$g\colon X\to Y$. We construct symmetric spectra $Q^{k}_{i}(g)$
(or $Q^{k}_{i}$ when $g$ is understood) for $k\geq 0$, $0\leq
i\leq k$ inductively as follows: $Q^{k}_{0}=X^{(k)}$,
$Q^{k}_{k}=Y^{(k)}$ (the $k$-th smash power of $X$ and $Y$), and
for $0<i<k$, we define $Q^{k}_{i}$ by the pushout square:
$$ \xymatrix{
\Sigma_{k+}\sma_{\Sigma_{k-i}\times \Sigma_{i}}X^{(k-i)}\sma Q^{i}_{i-1}
\ar[r]\ar[d]
&\Sigma_{k+}\sma_{\Sigma_{k-i}\times \Sigma_{i}}X^{(k-i)}\sma Y^{(i)}
\ar[d]\\
Q^{k}_{i-1}\ar[r]&Q^{k}_{i} } $$ Essentially, $Q^{k}_{i}$ is the
$\Sigma_{k}$-sub-spectrum of $Y^{(k)}$ of with $i$ factors of $Y$
and $k-i$ factors of $X$: The quotient $Y^{(k)}/Q^{k}_{k-1}$ is
naturally isomorphic to $(Y/X)^{(k)}$.  When $g$ is $F_{m}$ of an
injection of simplicial sets $X\to Y$, $Q^{k}_{i}$ is precisely $F_{mk}$
of the subspace of $Y^{k}$ where at most $i$ factors are in $Y
\setminus X$.

Combining these constructions, we get a filtration on
$B\amalg_{\MM\IO{x}X}\MM\IO{x}Y$ as follows.  Let $B_{0}=B$, and
let $B_{k}$ be the pushout in $\S^{O}$
$$ \xymatrix{
\U{k}B\sma_{\Sigma_{k}}Q^{k}_{k-1}\ar[r]\ar[d]
&\U{k}B\sma_{\Sigma_{k}}\IO{x}Y^{(k)}\ar[d]\\
B_{k-1}\ar[r]&B_{k}, }
$$
where the map $\U{k}B\sma_{\Sigma_{k}}Q^{k}_{k-1}\to B_{k-1}$ is
induced by the map $\IO{x}X\to B$.  Let $B_{\infty}=\Colim B_{k}$.

\proclaim{Proposition}\label{pushfilt}
With notation above, $B_{\infty}$ is isomorphic to the underlying
object of $B\amalg_{\MM\IO{x}X}\MM\IO{x}Y$ in $\S^{\OB}$.
\endproclaim

In order to use this below, we need to know that the map
$B_{k-1}\to B_{k}$ is objectwise a level cofibration of symmetric
spectra.

\proclaim{Lemma}\label{coflemma}
Let $T$ be any right $\Sigma_{k}$ object in symmetric spectra. If
$g\colon X\to Y$ is a cofibration, then $T\sma_{\Sigma_{k}}
Q^{k}_{k-1}(g)\to T\sma_{\Sigma_{k}}Y^{(k)}$ is a level
cofibration, i.e., level injection.
\endproclaim

\demo{Proof}
It suffices to consider the case when $X\to Y$ is a relative
$\IP$-complex, and a filtered colimit argument reduces to the case
when $X\to Y$ is formed by attaching a single cell, i.e., is the
pushout over a map
$$ F_{m}i\colon F_{m}\partial \Delta[n]_{+}\to F_{m}\Delta[n]_{+} $$
in $\IP$. Then the map in the statement is the pushout over the
map
$$ T\sma_{\Sigma_{k}} Q^{k}_{k-1}(F_{m}i)\to
T\sma_{\Sigma_{k}}(F_{m}\Delta[n]_{+})^{(k)}. $$ We can identify
this as $T\sma_{\Sigma_{k}}(-)$ applied to the map
$$ F_{mk}\partial (\Delta[n]^{k})_{+}\to F_{mk}\Delta[n]^{k}_{+}.  $$
It is easy to check explicitly that this is a level cofibration.
\enddemo

\demo{Proof of \refname{unitthm}}
It suffices to consider the case when $A$ is an $\MI$-complex,
i.e., the map from the initial object $\M(;-)_{+}\sma S$ to $A$ is
a relative $\MI$-complex.  Then $A=\Colim A_{n}$ where
$A_{0}=\M(;-)_{+}\sma S$, and $A_{n+1}$ is formed from $A_{n}$ as
a pushout over a coproduct of maps in $\IP$.  Since $f^{*}f_{*}A =
\Colim f^{*}f_{*}A_{n}$, it suffices to show that $A_{n}\to
f^{*}f_{*}A_{n}$ is a weak equivalence for all $n$.

We prove this by induction on $n$ for all $A_{n}$.  Specifically,
we say that an $\MI$-complex $B$ can be {\bf built in $n$ stages}
if, starting with $B_{0}=\M(;-)_{+}\sma S$, we can construct $B$
as a sequence of $n$ pushouts over coproducts of maps in $\MI$,
$B_{0}\to B_{1}\to \dotsb \to B_{n}=B$.  Our inductive hypothesis
is that for any $\MI$-complex $B$ that can be built in $n$ stages,
$B\to f^{*}f_{*}B$ is a stable equivalence.  Since $f$ is a weak
equivalence, $\M(;-)_{+}\sma S\to\Mp(;-)_{+}\sma S$ is a stable
equivalence, and this gives the base case $n=0$.  Our argument
also needs the base case $n=1$, where we are looking at a map of
the form $\MM T\to f^{*}\MMp \IO{f}T$ for some $T$ in $\S^{\OB}$
that is objectwise cofibrant.  Using the explicit formula for
$\MM$ and $\MMp$ in \refname{defmm}, we see that this is a stable
equivalence.

For the inductive step from $n$ to $n+1$, a filtered colimit
argument reduces to the case of $C=B\amalg_{\MM\IO{x}X}\MM\IO{x}Y$
for $X\to Y$ in $\IP$, where $B$ can be built in $n$ stages.  We
have the filtration preceding \refname{pushfilt},
$$ B=B_{0}\to B_{1}\to \dotsb, \qquad C=B_{\infty}=\Colim B_{k}, $$
whose associated graded is
$$ \bigvee_{k} \U{k}B\sma_{\Sigma_{k}}(Y/X)^{(k)},
$$
which is isomorphic in $\S^{\OB}$ to $B\amalg \MM\IO{x}(Y/X)$,
with the coproduct in $\S^{\M}$.  Let $B'=f_*B$ and $C'=f_*C$.
Since $C'=B'\amalg_{\MMp\IO{fx}X}\MMp\IO{fx}Y$, we have the
analogous filtration
$$ B'=B'_{0}\to B'_{1}\to\dotsb, \qquad C'=B'_{\infty}=\Colim B'_{k},
$$
whose associated graded is isomorphic in $\S^{\OBp}$ to $B'\amalg
\MMp\IO{fx}(Y/X)$.  The map $C\to f^{*}C'=\pi_{f}C'$ preserves the
filtrations, and the map of associated gradeds
$$ B\amalg \MM\IO{x}(Y/X) \to \pi_{f}(B'\amalg \MMp\IO{fx}(Y/X) \iso
f^{*}f_{*}(B\amalg\MM\IO{x}(Y/X))
$$
is a stable equivalence, because
$B\amalg \MM\IO{x}(Y/X)$ can be built in $n$ stages (since $n\geq
1$). By \refname{coflemma}, the maps in the filtration are
objectwise level cofibrations, and it follows that each map
$B_{k}\to \pi_{f}B_{k}$ is a stable equivalence.  The map $C\to
\pi_{f}C'=f^{*}f_{*}C$ is therefore a stable equivalence.
\enddemo

The constructions in this section also provide what is needed for
the proof of Lem\-ma~\refnum{mklemma}.

\demo{Proof of \refname{mklemma}}
A filtered colimit argument reduces to showing that the map $B\to
B\amalg_{\MM\IO{x}X}\MM\IO{x}Y$ is a stable equivalence for $X\to
Y$ in $\KP$.  Let $B=B_{0}\to B_{1}\to\dotsb$ be as above
\refname{pushfilt};
it suffices to show that each $B_{k-1}\to B_{k}$ is a stable
equivalence.  The quotient $B_{k}/B_{k-1}$ is naturally isomorphic
to $\U{k}B\sma_{\Sigma_{k}}(Y/X)^{(k)}$.  Moreover, $Y/X$ is
positive cofibrant and stably equivalent to the trivial symmetric
spectrum $*$, and so $B_{k}/B_{k-1}$ is stably equivalent to the
trivial object $*$ in $\S^{\OB}$. Since the map $B_{k-1}\to B_{k}$
is objectwise a level cofibration, it follows that it is a stable
equivalence.
\enddemo

\enddocument
\bye